\documentclass[letterpaper]{article} % DO NOT CHANGE THIS
\usepackage{aaai2026}  % DO NOT CHANGE THIS
\usepackage{times}  % DO NOT CHANGE THIS
\usepackage{helvet}  % DO NOT CHANGE THIS
\usepackage{courier}  % DO NOT CHANGE THIS
\usepackage[hyphens]{url}  % DO NOT CHANGE THIS
\usepackage{graphicx} % DO NOT CHANGE THIS
\usepackage{natbib}  % DO NOT CHANGE THIS AND DO NOT ADD ANY OPTIONS TO IT
\usepackage{caption} % DO NOT CHANGE THIS AND DO NOT ADD ANY OPTIONS TO IT
\usepackage{algorithm}
\usepackage{algorithmic}

\usepackage{newfloat}
\usepackage{listings}
\DeclareCaptionStyle{ruled}{labelfont=normalfont,labelsep=colon,strut=off} % DO NOT CHANGE THIS
\floatstyle{ruled}
\newfloat{listing}{tb}{lst}{}
\floatname{listing}{Listing}
\usepackage{subfig}	
\usepackage[shortlabels]{enumitem}	
\usepackage{booktabs}	
\usepackage{multirow}	
\usepackage{amsmath}	
\usepackage{subeqnarray}
\usepackage{cases}      
\usepackage{amsthm}     
\usepackage{amsfonts,amssymb}
\usepackage{color}
\usepackage{bm}
\usepackage[T1]{fontenc}
\usepackage{makecell} % \Xhline

\newtheorem{assumption}{Assumption}
\newtheorem{lemma}{Lemma}
\newtheorem{theorem}{Theorem}
\newtheorem{definition}{Definition}
\newtheorem{remark}{Remark}

\newcommand{\rmnum}[1]{\romannumeral #1}
\newcommand{\Rmnum}[1]{\uppercase\expandafter{\romannumeral #1}}

\title{Parallelizable Riemannian Alternating Direction Method of
	Multipliers for Non-convex Pose Graph Optimization}

\author {
	Xin Chen\textsuperscript{\rm 1},
	Chunfeng Cui\textsuperscript{\rm 1}\thanks{Corresponding author.},
	Deren Han\textsuperscript{\rm 1},
	Liqun Qi\textsuperscript{\rm 2},
}
\affiliations {
	\textsuperscript{\rm 1}School of Mathematical Sciences, Beihang University, Beijing, China\\
	\textsuperscript{\rm 2}Department of Applied Mathematics, The Hong Kong Polytechnic University, HKSAR, China\\
	\{chenxin2020, chunfengcui, handr\}@buaa.edu.cn, maqilq@polyu.edu.hk
}

\usepackage{bibentry}
\begin{document}
	
	\maketitle
	
	\begin{abstract}
		Pose graph optimization (PGO) is fundamental to robot perception and navigation systems, serving as the mathematical backbone for solving simultaneous localization and mapping (SLAM). Existing solvers suffer from polynomial growth in computational complexity with graph size, hindering real-time deployment in large-scale scenarios. In this paper, by duplicating variables and introducing equality constraints, we reformulate the problem and propose a Parallelizable Riemannian Alternating Direction Method of Multipliers (PRADMM) to solve it efficiently. Compared with the state-of-the-art methods that usually exhibit polynomial time complexity growth with graph size,  PRADMM  enables efficient parallel computation across vertices regardless of graph size. Crucially, all subproblems admit closed-form solutions, 
		ensuring  PRADMM maintains exceptionally stable performance. 
		Furthermore, by carefully exploiting the structures of the coefficient matrices in the constraints, we establish the global convergence of PRADMM under mild conditions, enabling larger relaxation step sizes within the interval  $(0,2)$.   Extensive empirical validation on two synthetic datasets and multiple real-world 3D SLAM benchmarks confirms the superior computational performance of PRADMM. 
	\end{abstract}
	
	\begin{links}
		\link{Code}{https://github.com/HeartsHorizon/PRADMM}
	\end{links}

	\section{Introduction} \label{sec-introduction}
	%ccf: I have rewritten this paragraph, begining with PGO, not SLAM
	Pose graph optimization (PGO) \cite{lu1997globally} is a fundamental technique for trajectory estimation from noisy sensor data. This  optimization framework corrects cumulative errors in motion systems, serving as the computational core for  sensor networks \cite{so2007theory} and simultaneous localization and mapping (SLAM) \cite{smith1990estimating}. In robotics and autonomous driving, PGO enables SLAM by optimally reconciling pose estimates with environmental observations \cite{montemerlo2002fastslam}, while also supporting structure-from-motion \cite{martinec2007robust} and bundle adjustment \cite{bender2013graph}. 
	
	However, scalability remains a critical bottleneck: existing solvers suffer from polynomial growth in computational complexity with graph size, hindering real-time deployment in large-scale scenarios. To address this limitation, we propose PRADMM---a Riemannian ADMM framework with  convergence guarantees, enabling parallel computation across graph vertices. 
	
	%with closed-form subproblem solutions.
	
	\subsection{Related Work}
	\subsubsection{Pose Graph Optimization}
	Early first-order methods such as stochastic gradient descent \cite{olson2006fast,grisetti2009nonlinear} reduced computational complexity, while second-order techniques such as Gauss-Newton \cite{carlone2015duality}, Levenberg-Marquardt \cite{kummerle2011g2o} and trust-region methods \cite{rosen2012incremental} achieved faster local convergence. Nevertheless, these approaches remain susceptible to convergence at local minima and have high computational complexity in Hessian construction. Efforts to construct sparse Hessians based on graph connectivity \cite{grisetti2010tutorial} ultimately fail to address this issue in dense graph models. To address non-convexity, initialization schemes such as chordal \cite{carlone2015initialization} and rotation synchronization \cite{nasiri2018linear} have been proposed. 
	Unlike prior methods, convex relaxation techniques guarantee convergence to certifiably global optima regardless of initialization \cite{carlone2015lagrangian,rosen2019se}. However, solving the resultant semi-definite programming or its low-rank approximations remains computationally prohibitive. Recent work further accelerated convergence via problem-aware first-order methods \cite{fan2023majorization} and trigonometric parameterizations \cite{nasiri2020novel}, balancing speed and accuracy. However, these methods scale poorly for large-scale or densely connected graphs, demanding the development of algorithms with near-constant complexity relative to graph scale.
	%with weak dependency on graph scale.

	%  {\bf Pose graph optimization}
	% Early first-order methods (e.g., stochastic gradient descent \cite{olson2006fast,grisetti2009nonlinear}) reduced computational complexity, while second-order techniques (Gauss-Newton \cite{carlone2015duality}, Levenberg-Marquardt \cite{kummerle2011g2o} and trust-region methods \cite{rosen2012incremental}) achieved faster local convergence--though often at the cost of local minima and high Hessian complexity.   To address non-convexity, initialization schemes like chordal \cite{carlone2015initialization} and rotation synchronization \cite{nasiri2018linear} were proposed, while convex relaxations \cite{rosen2019se} enabled certifiably global solutions under noise thresholds. Recent work further accelerated convergence via problem-aware first-order methods \cite{fan2023majorization} and trigonometric parameterizations \cite{nasiri2020novel}, balancing speed and accuracy. 
	%  \red{However, xx}
	\begin{table*}[tbp]
		\centering
		\begin{tabular}{lccccc}
			\toprule[1pt]
			Algorithm& Manifold & $A_{N}$ & Stepsize &Conv.& G.C. \\
			\midrule
			ADMM \cite{wang2019global}  &  $\mathbb{R}^{n}$ & $\operatorname{Im}(A_{-N}) \subseteq \operatorname{Im}(A_{N})$& $1$ & S.C.& \checkmark\\
			Prox-ADMM \cite{boct2020proximal}  &  $\mathbb{R}^{n}$ &  surjective& $(0,2)$ & S.C.& \checkmark\\
			iADMM \cite{hien2022inertial}  &  $\mathbb{R}^{n}$ &  surjective& $(0,2)$ & S.C.& $\checkmark_1$\\
			\midrule
			SOC \cite{lai2014splitting} &   Stiefel& $I_{n}$& $1$ & -& -\\
			MADMM \cite{kovnatsky2016madmm}  & Riemannian& $I_{n}$& $1$ & -& -\\
			ADMM-NSSC \cite{lu2018nonconvex}  & Orthogonal & $I_{n}$& $1$ & S.C.& -\\
			Prox-ADMM \cite{zhang2020primal}  & Riemannian& $I_{n}$& $1$ & E.C.& - \\
			RADMM \cite{li2022riemannian}  & Riemannian& $I_{n}$& $1$ & E.C.& - \\
			PieADMM \cite{chen2024pieadmm}  & Spherical& $I_{n}$& $1$ & E.C.& - \\
			PRADMM (Algorithm \ref{alg1}) & Riemannian & injective& $(0,2)$ & S.C. & \checkmark \\
			\bottomrule[1pt]
		\end{tabular}
		{\raggedright 
			``S.C.'', ``E.C.'' and `` G.C.'' represent ``subsequence convergence'', `` ergodic complexity'' and ``global convergence''.  \\} 
		{\raggedright 
			$A_{-N}=[A_{1},\dots,A_{N-1}]$, and relaxed stepsize $\tau$ comes from the update of dual variables (Algorithm \ref{alg1}).\\} 
		{\raggedright The symbol $\checkmark_1$  denotes validity solely for $\tau=1$.\\}
		\caption{Comparison of recent developments of Riemannian ADMM and part of Euclidean ADMM for non-convex problems.}
		\label{tab:related work of PRADMM}
	\end{table*}
	
	\subsubsection{Riemannian ADMM}  
	ADMM is an efficient algorithm widely applied in large-scale optimization across diverse domains 
	\cite{gabay1976admm,Boyd2011}.
	The emergence of optimization problems with manifold constraints has prompted increasing research interest in extending vanilla ADMM (see Table \ref{tab:related work of PRADMM}). Early approaches such as SOC \cite{lai2014splitting} and MADMM \cite{kovnatsky2016madmm} demonstrated empirical effectiveness but lacked theoretical convergence guarantees. Subsequent advancements, including ADMM-NSSC \cite{lu2018nonconvex}, Prox-ADMM \cite{zhang2020primal}, RADMM \cite{li2022riemannian}, and PieADMM \cite{chen2024pieadmm},  established rigorous convergence foundations for Riemannian ADMM variants. Notably, \citeauthor{li2022riemannian} \shortcite{li2022riemannian} enhanced the applicability of PRADMM to nonsmooth objectives via the Moreau smoothing technique. However, two fundamental limitations persist: (i) Existing theories \cite{lu2018nonconvex,zhang2020primal,li2022riemannian,chen2024pieadmm} uniformly require the last block coefficient matrix $A_{N}$ in the linear equality constraints $\sum_{i=1}^{N}A_{i}x_{i}=b$ to be the identity matrix. This restriction persists even in Euclidean nonconvex ADMM theory \cite{wang2019global,boct2020proximal,hien2022inertial}, where $A_{N}$ must typically be bijective or surjective, conditions that fail to hold in specific application scenarios. (ii) All these manifold-constrained ADMM analyses establish only subsequence convergence or ergodic complexity, and exclusively for the relaxed dual stepsize $\tau=1$. Proving convergence for a wider range of step sizes and establishing global convergence remain significant theoretical challenges.
	
	% The extension of vanilla ADMM (Table \ref{tab:related work of PRADMM}) to manifold-constrained optimization has attracted growing research interest. Early approaches such as SOC \cite{lai2014splitting} and MADMM \cite{kovnatsky2016madmm} demonstrated empirical effectiveness but lacked theoretical convergence guarantees. Subsequent work in  \cite{zhang2020primal} and  \cite{li2022riemannian} established rigorous convergence analysis for Riemannian ADMM variants, advancing the theoretical foundations of manifold-constrained optimization. 
	% \red{However, xx} 

	\subsection{Main Contributions}  
	Our specific contributions are summarized as follows.
	\begin{itemize}
		\item We decouple vertex correlations through variable splitting and constraint propagation and develop PRADMM with closed-form subproblem solutions. Given the sparsity measurement $s$, each subproblem complexity is $\mathcal{O}(s)$, which is near-constant relative to graph scale. 
		
		%We investigate the pose graph optimization problem and achieve decoupling of vertex correlations through mathematical reformulation. Then we develop PRADMM with closed-form subproblem solutions, enabling fully parallelizable vertex updates and stable computation for large-scale PGO.
		
		\item We establish global convergence under row rank-deficient coefficient matrices, significantly weaker conditions than existing requirements. Concurrently, we permit extended relaxed dual steps $\tau \in (0,2)$, surpassing conventional limits ($\tau \in (0,(1+\sqrt{5})/2)$ for convex and $\tau = 1$ for non-convex cases).

		\item Our algorithm is evaluated on two synthetic datasets at varying scales and five widely-used 3D SLAM benchmark datasets. Experimental results demonstrate that PRADMM consistently outperforms state-of-the-art graph-SLAM methods. As the problem scale increases, iPRADMM exhibits superior computational scalability with significantly less performance degradation.
	\end{itemize}
	
	% First, we develop PRADMM with closed-form subproblem solutions, enabling fully parallelizable vertex updates and stable computation for large-scale PGO. Second, we prove global convergence permitting relaxed dual steps ($\tau \in (0,2)$), surpassing conventional limits ($\tau \in (0,(1+\sqrt{5})/2)$ for convex and $\tau = 1$ for non-convex cases). Third, experiments show PRADMM consistently outperforms state-of-the-art graph-SLAM methods. 
	
	\section{Model} \label{sec-model}
	
\begin{figure}[!thbp]
		\centering
		\includegraphics[width=0.92\linewidth]{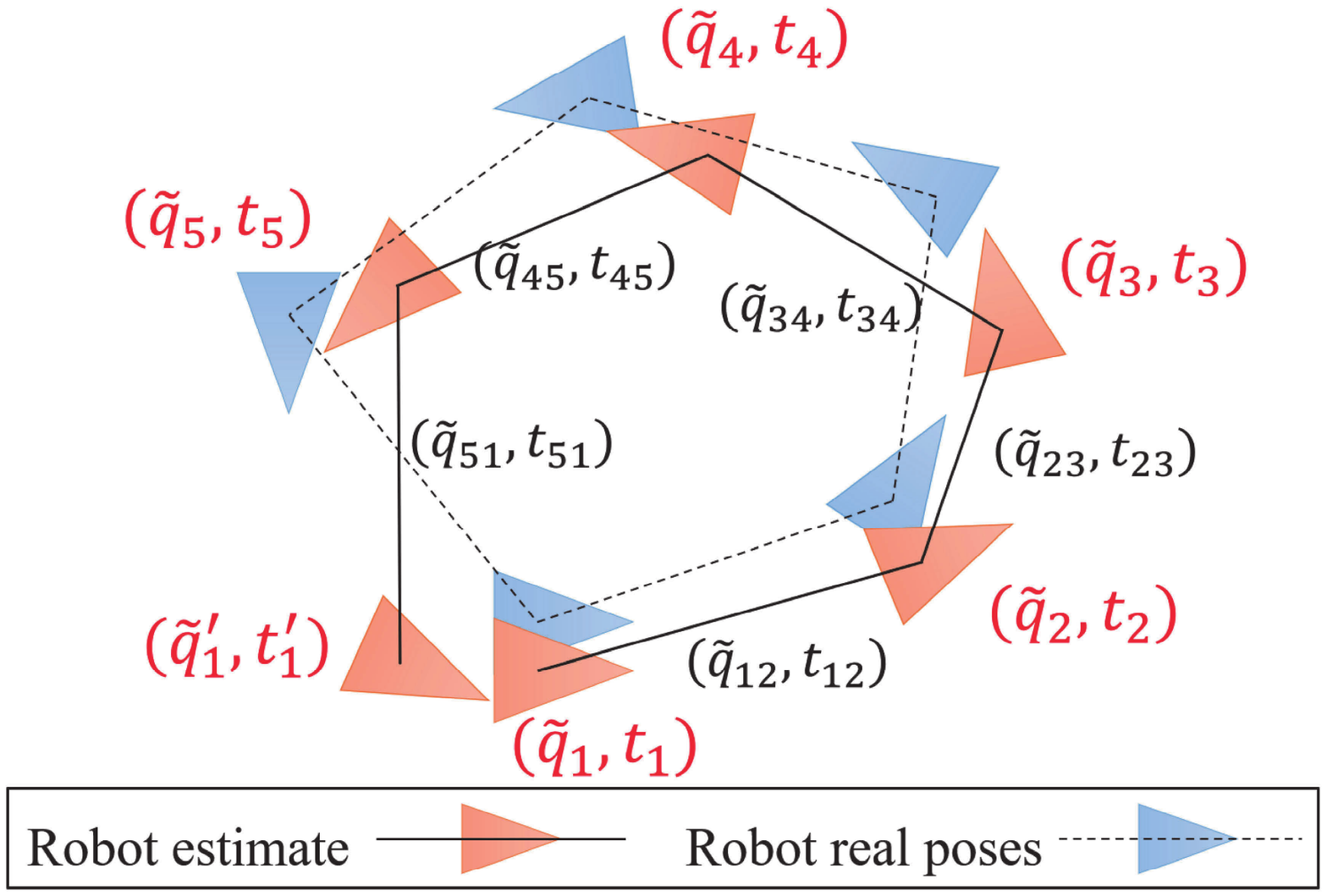}
		\caption{A pose-graph representation of the  SLAM process.}
		\label{fig-PGO}
	\end{figure}

	PGO is mathematically represented by a directed graph $\mathcal{G} = (\mathcal{V}, \mathcal{E})$ (Figure \ref{fig-PGO}). The vertex set $\mathcal{V}$ contains $n = |\mathcal{V}|$ nodes, each corresponding to an unknown robot pose $(\tilde{q}_{i},\bm{t}_{i})$. The edge set $\mathcal{E}$ comprises $m = |\mathcal{E}|$ directed arcs, where each edge $(i,j)$ represents a relative measurement $(\tilde{q}_{ij},\bm{t}_{ij})$ between poses. The objective is to recover the $n$ unknown poses from these $m$ noisy relative observations. 
	Here, $\tilde{q}_{i},\tilde{q}_{ij}\in\mathbb U$ are unit quaternions, and $\bm{t}_{i},\bm{t}_{ij}\in\mathbb R^n$ are 3-dimensional vectors for any $i,j=1,\dots,n$.
	
	Given that the prior noise distribution is typically unknown in real-world scenarios, we leverage an augmented unit quaternion formulation with the following generative model \cite{chen2024pieadmm}:
	\begin{equation*}
		\begin{aligned}
			\bm{t}_{ij}&=R_{i}^{\top}(\bm{t}_{j}-\bm{t}_{i})+\bm{t}_{\epsilon},    \quad \text{where}~ \bm{t}_{\epsilon} \sim  ~ \mathcal{N}(0,\Omega_{1}),\\
			\tilde{q}_{ij}&=\tilde{q}_{i}^{*}\tilde{q}_{j}\tilde{q}_{\epsilon}, \quad\quad\quad\quad\quad~ \text{where}~\tilde{q}_{\epsilon}\sim  ~ \text{vMF}([1,0,0,0],\kappa),
		\end{aligned}
	\end{equation*}
	where $R_{i}$ and $\tilde{q}_{i}$ are the rotation representation of vertex $i$ in $SO(3)$ and unit quaternion, respectively. 
	%``$\mathcal{N}(\bm{\mu},\Omega)$'' denotes a Gaussian distribution with mean $\bm{\mu}$ and covariance matrix $\Omega$. 
	``$\text{vMF}(\bm{\mu},\kappa)$'' denotes a $d$-dimensional von Mises-Fisher distribution where $\bm{\mu} \in \mathbb{S}^{d-1}$ and $\kappa \geq 0$ are mean direction and concentration parameters, respectively. It is one of the most commonly used distributions to model data distributed on the surface
	of the unit hypersphere \cite{fisher1953dispersion,sra2012short,hornik2014maximum}. 
	%and can be considered a circular analogue of the normal distribution. 
	In fact, as  $\kappa$ increases, the vMF distribution becomes increasingly concentrated at the mean direction $\bm{\mu}$. When $\kappa=0$, it corresponds to the uniform distribution on $\mathbb{S}^{d-1}$. When $\kappa\rightarrow+\infty$, the distribution approximates a Gaussian distribution with mean $\bm{\mu}$ and covariance  $1/\kappa$.
	The MLE of the graph vertices yields the optimization problem
	\begin{equation}
		\begin{aligned}
			\underset{\{\tilde{q}_i, \bm{t}_{i}\}}{\min}	&
			\sum_{(i,j)\in \mathcal{E}} 
			\|\tilde{t}_j - \tilde{t}_i-\tilde{q}_{i}\tilde{t}_{ij}\tilde{q}_{i}^{*}\|_{\Sigma_{1}}^2
			\!+\!\|\tilde{q}_{j}^{*}\tilde{q}_{i}\tilde{q}_{ij}-1\|_{\Sigma_{2}}^2 \\
			\operatorname{s.t.} ~ ~&
			\tilde{q}_i \in \mathbb{U},~ \bm{t}_{i} \in \mathbb{R}^3,~ i=1,2,\dots,n, 
		\end{aligned}
		\label{PGO-model1}
	\end{equation}
	where $\Sigma_{1} = \left( \begin{array}{cc}
		c & 0\\
		0 & \Omega_{1}^{-1}
	\end{array}\right)$ and $\Sigma_{2}=\kappa I$ are positive semidefinite matrices, and $c>0$ is an arbitrary scalar. $\tilde{t}_{i}=[0,\bm{t}_{i}]$ is a pure quaternion.  
	
	\section{Algorithm}\label{sec-algorithm}
	\subsection{Reformulation via Structured Splitting}
	Scalability remains a critical bottleneck in PGO: existing solvers such as 
	the proximal Riemannian gradient method \cite{gabay1982minimizing}, Riemannian conjugate gradient \cite{smith1994optimization}, and Riemannian Newton method \cite{hu2018adaptive} 
	suffer from polynomial growth in computational complexity with graph size.
	This bottleneck arises because computing gradients requires global pose information, resulting in high complexity $\mathcal{O}(m)$ and $\mathcal{O}(n^{2})$ for rotations and translations, respectively. 
	
	To overcome this limitation, we employ a splitting strategy that decouples these vertices. The existing splitting strategy, SOC \cite{lai2014splitting}, in Riemannian optimization separates manifold constraints from variables. 
	However, the absence of closed-form solutions in subproblems fundamentally prevents computational complexity reduction via this approach. The cost of SOC is  $\mathcal{O}(km+n^{2})$ per iteration, where $k$ is the number of internal cycles. Another splitting strategy, PieADMM \cite{chen2024pieadmm}, decouples rotation variables across vertices, achieving the computational complexity $\mathcal{O}(m/n)$. However, by overlooking the iteration cost for translation subproblems, the cost of $\bm{t}$ is still $\mathcal{O}(n^{2})$. By introducing auxiliary variables $\tilde{p}_i$ and $\bm{s}_{i}$, we reformulate \eqref{PGO-model1} into an equivalent model: 
	\begin{equation}\label{PGO-model-splitting}
		\begin{aligned}
			\underset{\{\tilde{p}_i,\tilde{q}_i,\bm{t}_{i},\bm{s}_{i}\}}{\min}	&
			\!\sum_{(i,j)\in \mathcal{E}} \!
			\underbrace{\|\tilde{t}_j - \tilde{s}_i-\tilde{q}_{i}\tilde{t}_{ij}\tilde{p}_{i}^{*}\|_{\Sigma_{1}}^2}
			_{f(\tilde{\bm{p}},\tilde{\bm{q}},\bm{t},\bm{s})}
			\!+\!
			\underbrace{\|\tilde{p}_{j}^{*}\tilde{q}_{i}\tilde{q}_{ij}-1\|_{\Sigma_{2}}^2}_{g(\tilde{\bm{p}},\tilde{\bm{q}})} \\
			\operatorname{s.t.} ~~~~& 
			\tilde{p}_i=\tilde{q}_i,~
			\bm{t}_{i} = \bm{s}_{i},~
			\tilde{p}_i \in \mathbb{U},~
			\tilde{q}_i \in \mathbb{R}^{4},\\
			&\bm{t}_{i} \in \mathbb{R}^3,~ 
			\bm{s}_{i} \in \mathbb{R}^3,~
			i=1,\dots,n.
		\end{aligned}
	\end{equation}
	This reformulation yields two computational advantages, as demonstrated in the following:
	\begin{itemize}
		\item[(i)] The quartic polynomial in \eqref{PGO-model1} transforms into a multi-linear model, enabling closed-form solutions for all subproblems.  
		\item[(ii)] Vertex interdependencies are decoupled, facilitating massively parallel computation. Specifically, our PRADMM algorithm achieves $\mathcal{O}(m/n)$ computational complexity for both rotations and translations.
	\end{itemize}
	Note that $m$ depends on the density of the graph, i.e., $n \leq m \leq n^{2}$. Therefore, with fixed graph sparsity (i.e., constant ratio $m/n$), the theoretical iteration complexity becomes scale-invariant.
	
	\subsection{PRADMM}
	Let $\tilde{\bm{p}}=(\tilde{p}_{1},\tilde{p}_{2},\dots,\tilde{p}_{n})$, $\tilde{\bm{q}}=(\tilde{q}_{1},\tilde{q}_{2},\dots,\tilde{q}_{n})$,  $\bm{t}=(\bm{t}_{1},\bm{t}_{2},\dots,\bm{t}_{n})$, and $\bm{s}=(\bm{s}_{1},\bm{s}_{2},\dots,\bm{s}_{n})$. 
	% We define
	% \begin{align}
	% 	f(\tilde{\bm{p}},\tilde{\bm{q}},\bm{t},\bm{s})
	% 	&=\sum_{(i,j)\in \mathcal{E}} 
	% 	\|\tilde{t}_j - \tilde{s}_i-\tilde{q}_{i}\tilde{t}_{ij}\tilde{p}_{i}^{*}\|_{\Sigma_{1}}^2,
	% 	\label{function_f1}\\
	% 	g(\tilde{\bm{p}},\tilde{\bm{q}})
	% 	&=\sum_{(i,j)\in \mathcal{E}} 
	% 	\|\tilde{p}_{j}^{*}\tilde{q}_{i}\tilde{q}_{ij}-1\|_{\Sigma_{2}}^2.
	% 	\label{function_g1}
	% \end{align}
	The augmented Lagrangian function of \eqref{PGO-model-splitting} is 
	\begin{equation*}
		\begin{aligned}
			\mathcal{L}_{\beta}(\tilde{\bm{p}},\tilde{\bm{q}},&\bm{t},\bm{s},\bm{\lambda},\bm{z})=
			f(\tilde{\bm{p}},\tilde{\bm{q}},\bm{t},\bm{s})
			+g(\tilde{\bm{p}},\tilde{\bm{q}})\\
			&+\sum_{i=1}^{n} \left\lbrace 
			-\left\langle \bm{\lambda}_{i}, \tilde{p}_{i}-\tilde{q}_{i} \right\rangle
			+\frac{\beta_1}{2}\|\tilde{p}_{i}-\tilde{q}_{i}\| ^2\right\rbrace \\
			& +\sum_{i=1}^{n} \left\lbrace 
			-\left\langle \bm{z}_{i}, \bm{t}_{i}-\bm{s}_{i} \right\rangle
			+\frac{\beta_2}{2}\|\bm{t}_{i}-\bm{s}_{i}\| ^2\right\rbrace,
		\end{aligned}
	\end{equation*}
	where $\bm{\lambda} \in \mathbb{R}^{4n}$ and $\bm{z} \in \mathbb{R}^{3n}$ are the Lagrange multipliers and $\beta_1,\beta_2>0$ are penalty parameters. In terms of algorithm design,  we use the over-relaxation technique, which allows a larger step size for dual variables. The accelerated algorithm may reduce the number of iterations and achieve faster convergence. The iterative scheme of Parallelizable Riemannian ADMM is given by
	\begin{small}
		\begin{subequations}
			\begin{numcases}{}
				\tilde{\bm{p}}^{k+1}\!\in
				\mathop{\arg\min}\limits_{\tilde{\bm{p}}\in\mathbb U^n}
				\mathcal{L}_{\beta}(\tilde{\bm{p}},\tilde{\bm{q}}^{k},\bm{t}^{k},\bm{s}^{k},\bm{\lambda}^{k},\bm{z}^{k})
				+\mathcal{Q}_{1}(\tilde{\bm{p}})
				\label{subproblem_p_1}\\
				\tilde{\bm{q}}^{k+1}\!\in
				\mathop{\arg\min}\limits_{\tilde{\bm{q}}}
				\mathcal{L}_{\beta}(\tilde{\bm{p}}^{k+1},\tilde{\bm{q}},\bm{t}^{k},\bm{s}^{k},\bm{\lambda}^{k},\bm{z}^{k})
				+\mathcal{Q}_{2}(\tilde{\bm{q}})
				\label{subproblem_q_1}\\
				\bm{t}^{k+1}\!\in
				\mathop{\arg\min}\limits_{\bm{t}}
				\mathcal{L}_{\beta}(\tilde{\bm{p}}^{k+1},\tilde{\bm{q}}^{k+1},\bm{t},\bm{s}^{k},\bm{\lambda}^{k},\bm{z}^{k})
				+\mathcal{Q}_{3}(\bm{t})
				\label{subproblem_t_1}\\
				\bm{s}^{k+1}\!\in
				\mathop{\arg\min}\limits_{\bm{s}}
				\mathcal{L}_{\beta}(\tilde{\bm{p}}^{k+1}\!,\tilde{\bm{q}}^{k+1}\!,\bm{t}^{k+1}\!,\bm{s},\bm{\lambda}^{k}\!,\bm{z}^{k})
				\!+\!\mathcal{Q}_{4}(\bm{s})
				\label{subproblem_s_1}\\
				\bm{\lambda}^{k+1}=\bm{\lambda}^{k}-\tau\beta_1(\tilde{\bm{p}}^{k+1}-\tilde{\bm{q}}^{k+1})
				\label{subproblem_lambda_1}\\
				\bm{z}^{k+1} = \bm{z}^{k}-\tau\beta_2(\bm{t}^{k+1}-\bm{s}^{k+1})\label{subproblem_z_1}
			\end{numcases}
		\end{subequations}
	\end{small}
	where $\mathcal{Q}_{i}(x)=\frac{1}{2}\|x-x^{k}\|_{H_{i}}^{2}$; $H_i \succ 0$, $i=1,\dots,4$, are positive definite matrices with diagonal structures. Here, $\frac{1}{2}\|\tilde{\bm{p}}-\tilde{\bm{p}}^{k}\|_{H_1}^2$ can be split as 
	$\frac{1}{2}\sum_{i=1}^{n}\|\tilde{p}_{i}-\tilde{p}_{i}^{k}\|_{H_{1,i}}^{2}$. 
	%and the same as $\tilde{\bm{q}}$. 
	In other words, each block of $H_{1}$ is still a diagonal matrix with   $H_{1,i} = \gamma_{1,i}I_{4}$. 
	Assume $H_2, H_3, H_4$ have similar structures. 
	
	\subsection{Subproblems}
	We partition the given directed graph $\mathcal{G}=(\mathcal{V},\mathcal{E})$ according to the vertices.
	We define $\mathcal{E}_{i}^{in}=\{(l,i)\in\mathcal{E} \}$ for all $l \in \mathcal{V}$, and $\mathcal{E}_{i}^{out}=\{(i,j)\in\mathcal{E} \}$ for all $j \in \mathcal{V}$. In other words, $\mathcal{E}_{i}^{in}$ represents all directed edges that pointing to vertex $i$, while $\mathcal{E}_{i}^{out}$ is the opposite. Then we have the properties that 
	\begin{align*}
		&\mathcal{E}=\bigcup_{i\in\mathcal{V}}\left( \mathcal{E}_{i}^{in}\cup \mathcal{E}_{i}^{out} \right),\\
		&\mathcal{E}_{i}^{in} \cap \mathcal{E}_{i}^{out} = \varnothing \text{ for all } i\in\mathcal{V},\\
		&\mathcal{E}_{i}^{in} \cap \mathcal{E}_{j}^{in} = \varnothing ,\text{ for all } i\neq j.
	\end{align*}
	Next, we analyze the subproblems \eqref{subproblem_p_1}--\eqref{subproblem_s_1}. All detailed derivations can be found in the appendix.
	\subsubsection{$\tilde{\bm{p}}$--subproblem:}
	The problem \eqref{subproblem_p_1} can be reformulated as 
	\begin{align}
		&\mathop{\arg\min}\limits_{\tilde{\bm{p}} \in \mathbb{U}^{n}}  
		\sum_{i=1}^{n}  \left\lbrace  \sum_{(i,j)\in \mathcal{E}_{i}^{out}}
		\!\!\!\|M(\tilde{q}_{i}^{k})M(\tilde{t}_{ij})D\tilde{p}_{i}-(\tilde{t}_j^{k} - \tilde{s}_i^{k}) \|_{\Sigma_{1}}^2\right. \nonumber\\
		&\quad  + \sum_{(l,i)\in \mathcal{E}_{i}^{in}}\|W(\tilde{q}_{li})W(\tilde{q}_{l}^{k})\tilde{p}_{i}-1\|_{\Sigma_{2}}^2\nonumber\\
		&\quad\left.+ \frac{\beta_1}{2}\|\tilde{p}_{i}-(\tilde{q}_{i}^{k}+\frac{1}{\beta_1}\bm{\lambda}_{i}^{k})\|^{2}
		+\frac{1}{2}\|\tilde{p}_{i}-\tilde{p}_{i}^{k}\|_{H_{1,i}}^2\right\rbrace. 
		\label{subproblem_p_2-1}
	\end{align}
	where $M(\cdot)$, $W(\cdot) \in \mathbb{R}^{4\times 4}$ are quaternion-generated matrices designed to simplify quaternion multiplication. The matrix $D=\text{diag}(1,-1,-1,-1)$ is a diagonal matrix. Since $\tilde{p}_{i}$ are fully separable, we can update them in parallel, i.e.,
	\begin{align}
		\tilde{p}_{i}^{k+1} 
		= & \mathop{\arg\min}\limits_{\tilde{p}_{i} \in \mathbb{U} } ~ \frac{1}{2}\tilde{p}_{i}^{\top}A_{1,i}^{k}\tilde{p}_{i} + (b_{1,i}^{k})^{\top}\tilde{p}_{i} \label{subproblem_p_2}
	\end{align}
	% \red{where $G_{1,ij}^k = M(\tilde{q}_{i}^{k})M(\tilde{t}_{ij})D$,  $G_{2,ij}^k =W(\tilde{q}_{li})W(\tilde{q}_{l}^{k})$, and
	% \begin{align*}
	%     A_{1,i}^{k} = &
	% 2\sum_{(i,j)\in \mathcal{E}_{i}^{out}}
	% \left\lbrace(G_{1,ij}^k)^{\top}\Sigma_{1}G_{1,ij}^k\right\rbrace \\
	% &+2\sum_{(l,i)\in \mathcal{E}_{i}^{in}}
	% \left\lbrace(G_{2,ij}^k)^{\top}\Sigma_{2}G_{2,ij}^k\right\rbrace
	% +\beta_{1}I_{4}+H_{1,i},\\
	% b_{1,i}^{k} = &
	% -2\sum_{(i,j)\in \mathcal{E}_{i}^{out}}
	% (G_{1,ij}^k)^{\top}\Sigma_{1}(\tilde{t}_j^{k} - \tilde{s}_i^{k})\\
	% &-2\sum_{(l,i)\in \mathcal{E}_{i}^{in}}
	% (G_{2,ij}^k)^{\top}\Sigma_{2}\tilde{1}
	% -\beta_1\tilde{q}_{i}^{k}-\bm{\lambda}_{i}^{k}-H_{1,i}\tilde{p}_{i}^{k}.
	% \end{align*}}
	where $A_{1,i}^{k}\in\mathbb{R}^{4\times 4}$ and $b_{1,i}^{k}\in\mathbb{R}^{4}$ are obtained by rearranging \eqref{subproblem_p_2-1}. When $\Sigma_{1}$ and $\Sigma_{2}$ are both scalar matrices, $A_{1,i}^{k}$ is also a scalar matrix. The solution is $\tilde{p}_{i}^{k+1}=-b_{1,i}^{k}/\|b_{1,i}^{k}\|$ if $b_{1,i}^{k}$ is non-zero. When $\Sigma_{1}$ and $\Sigma_{2}$ are general matrices, the $\tilde{\bm{p}}$-subproblem \eqref{subproblem_p_2} is a special quadratic constraint quadratic programming with spherical constraint. 
	We show the results in Lemma \ref{TRS-subproblem3} that \eqref{subproblem_p_2} admits an eigenvalue problem.
	\begin{lemma}\label{TRS-subproblem3} \cite{adachi2017solving}
		Consider the spherical constrained problem
		\begin{equation}
			\min_{x\in\mathbb{R}^{n}} m(x)=g^{\top}x+\frac{1}{2}x^{\top}Ax,  \qquad 
			\operatorname{s.t.}~\|x\|= \Delta, 
			\label{eq-lemma4.2-1}
		\end{equation}
		where $A \in \mathbb{S}^{n \times n}$.
		For a solution $(x^{*},\lambda^{*})$ of the problem \eqref{eq-lemma4.2-1}, the multiplier $\lambda^{*}$ is equal to the largest real eigenvalue of $\tilde{Q}(\lambda)$, where
		$$
		\tilde{Q}(\lambda)=\left(\begin{array}{cc}
			-I & A+\lambda I\\
			A+\lambda I & \frac{-gg^{\top}}{\Delta^{2}}.
		\end{array}\right).
		$$
	\end{lemma}
	From the above lemma, problem \eqref{subproblem_p_2} aims at finding the largest real eigenvalue $\lambda^{*}$ such that $\operatorname{det}(\tilde{Q}(\lambda^{*}))=0$.
	This can be reformulated as a generalized eigenvalue problem, in which we need to compute the rightmost eigenvalue $\lambda_{*}$ such that $Q_1 y = -\lambda_{*}Q_2y$, where $$ Q_1 = \left(\begin{array}{cc}
		-I_4 & A_{1,i}^{k} \\
		A_{1,i}^{k} & - b_{1,i}^{k}(b_{1,i}^{k})^{\top}
	\end{array}\right),\qquad
	Q_2 = \left(\begin{array}{cc}
		0 & I_4 \\
		I_4 & 0
	\end{array}\right).
	$$
	%Note that $Q_{2}^{\top}Q_{2}=I$, and 
	Then we can derive an eigenvalue problem (Algorithm \ref{alg2}) for solving problem \eqref{subproblem_p_2}.
	\begin{algorithm}[t]
		\caption{Eigenvalue problem for solving $\tilde{p}$-subproblem \eqref{subproblem_p_1} .}
		\label{alg2}
		\hspace*{0.02in} 
		{ \textbf{Input:} $A_{1,i}^{k}$ and $b_{1,i}^{k}$, which are defined in \eqref{subproblem_p_2}.}
		
		\begin{algorithmic}[1]
			\STATE Compute the rightmost eigenvalue $\lambda_{*}$ such that $Q y = -\lambda_{*}y$, where $$ Q = \left(\begin{array}{cc}	
				A_{1,i}^{k} & - b_{1,i}^{k}(b_{1,i}^{k})^{\top}\\
				-I_4 & A_{1,i}^{k} 
			\end{array}\right).
			$$
			\STATE $\tilde{p}_{i}^{k+1} = -(A_{1,i}^{k} + \lambda_{*}I_4)^{-1}b_{1,i}^{k}$.
		\end{algorithmic}
		\hspace*{0.02in} {\bf Output: $\tilde{p}_{i}^{k+1}$.}
	\end{algorithm}
	
	\subsubsection{$\tilde{q}_{i}$, $\bm{t}_{i}$, and $\bm{s}_{i}$ subproblems}
	
	Similar as $\eqref{subproblem_p_2}$, we can update $\tilde{q}_{i}$, $\bm{t}_{i}$, and $\bm{s}_{i}$ in parallel. Since these subproblems can be formulated as least squares problems without manifold constraints, they admit directly closed-form solutions. 
	
	The pseudocode for solving \eqref{PGO-model-splitting} is summarized in Algorithm \ref{alg1}. All analytical expressions about $A_{s,i}^{k}$ and $b_{s,i}^{k}$ can be found in the appendix.
	\begin{algorithm}[t]
		\caption{PRADMM for solving PGO Model \eqref{PGO-model-splitting}.}
		\label{alg1}
		\hspace*{0.02in} 
		{ \textbf{Input:} \textbf{Graph} $\mathcal{G}=(\mathcal{V},\mathcal{E})$, 
			$\tilde{q}_{ij} \in \mathbb{U}$, $\bm{t}_{ij} \in \mathbb{R}^{3}$ for all $(i,j) \in \mathcal{E}$. \\
			\textbf{Initial points} $\tilde{\bm{p}}^{0} \in \mathbb{U}^{n}$, $\tilde{\bm{q}}^{0}\in \mathbb{U}^{n}$, $\bm{t}^{0} \in \mathbb{R}^{3n}$, $\bm{s}^{0}\in \mathbb{R}^{3n}$, ${\bm{\lambda}}^{0} \in \mathbb{R}^{4n}$, ${\bm{z}}^{0} \in \mathbb{R}^{4n}$}.\\
		\textbf{Parameters} $\beta_1,\beta_2>0$, $\tau \in \left(0,2\right) $, $H_{i}\succ0$, $i = 1,\dots,4$. (For the choice of these parameters, see Theorem \ref{thm-decrease}).
		
		\begin{algorithmic}[1]
			\FOR{$k=1,2,\dots$}
			\STATE Update $\tilde{p}_{i}^{k+1}$, $i\in[n]$ in parallel as follows 
			\begin{equation*}
				\tilde{p}_{i}^{k+1}=
				\left\lbrace 
				\begin{aligned}
					&-\frac{b_{1,i}^{k}}{\|b_{1,i}^{k}\|},   \quad \text{if}~ \Sigma_{i}~ \text{are scalar matrices},\\
					&\text{run} ~~\textbf{Algorithm} ~\textbf{\ref{alg2}}, \quad \text{otherwise}.
				\end{aligned}
				\right. 
			\end{equation*}
			
			\STATE Update 
			$\tilde{q}_{i}^{k+1}$, $i\in[n]$ in parallel by 
			$$\tilde{q}_{i}^{k+1}=(A_{2,i}^{k})^{-1}b_{2,i}^{k}.$$ 
			
			\STATE Update $\bm{t}_{i}^{k+1}$, $i\in[n]$ in parallel by $$\bm{t}_{i}^{k+1}=(A_{3,i}^{k})^{-1}b_{3,i}^{k}.$$
			
			\STATE Update $\bm{s}_{i}^{k+1}$, $i\in[n]$ in parallel by $$\bm{s}_{i}^{k+1}=(A_{4,i}^{k})^{-1}b_{4,i}^{k}.$$
			
			\STATE Update $\bm{\lambda},\bm{z}$ as follows
			$$\bm{\lambda}^{k+1}=\bm{\lambda}^{k}-\tau\beta_1(\tilde{\bm{p}}^{k+1}-\tilde{\bm{q}}^{k+1}),$$
			$$\bm{z}^{k+1} = \bm{z}^{k}-\tau\beta_2(\bm{t}^{k+1}-\bm{s}^{k+1}).$$
			\ENDFOR
		\end{algorithmic}
		\hspace*{0.02in} {\bf Output: $\tilde{\bm{p}}^{k+1} \in \mathbb{U}^{n}$ and $\bm{t}^{k+1} \in \mathbb{R}^{3n}$.}
	\end{algorithm}
	
	\section{Convergence Analysis} \label{sec-Convergence Analysis}
	For the ease of analysis, let us simplify the  notations in model \eqref{PGO-model-splitting} as 
	\begin{equation}\label{sec5-problem}
		\begin{aligned}
			\underset{x_{i}}{\min}	~~~~&
			f(x_{1},x_{2},x_{3},x_{4})+g(x_{1},x_{2}), \\
			\operatorname{s.t.} ~~~~& 
			\sum_{i=1}^{4}A_{i}x_{i}=0, ~x_{1} \in \mathcal{M},
		\end{aligned}	
	\end{equation}
	where $A_{1}=[I_{4n\times 4n},O_{3n\times 4n}]^{\top}$, $A_2=[-I_{4n\times 4n},$ $O_{3n\times 4n}]^{\top}$, $A_{3}=[O_{4n\times 3n},I_{3n\times 3n}]^{\top}$, $A_{4}=[O_{4n\times 3n},-I_{3n\times 3n}]^{\top}$. $O$ is a zero matrix.
	$\mathcal{M}\subseteq \mathbb{R}^{n_1}$ is a smooth Riemannian submanifold embedded in $n_1$-dimensional Euclidean space. For PGO model \eqref{PGO-model-splitting},  $\mathcal{M}$ is the Cartesian product of spheres. Denote $\bm{x} = (x_{1},x_{2},x_{3},x_{4})$, $\bm{x}_{-i} = (x_{1},\dots,x_{i-1},x_{i+1},x_{4})$,
	%$\bm{x}_{[i:j]} = (x_{i},\dots,x_{j})$, 
	and $\bm{x}^{k,i} = (\dots,x_{i-1}^{k+1},x_{i}^{k+1},x_{i+1}^{k},\dots)$. Let $\bm{x}^{k,0}=\bm{x}^{k}$ and $\bm{x}^{k,4}=\bm{x}^{k+1}$. 
	When we choose the same $\beta_1$ and $\beta_2$, Algorithm \ref{alg1} can be rewritten as 
	\begin{equation}\label{sec5-ADMM}
		\left\lbrace \begin{aligned}
			& \text{for}~ i=1,2,3,4,\\
			&\qquad x_{i}^{k+1}=\underset{x_{i}\in \mathcal{X}_{i}}{\arg\min}~
			\mathcal{L}_{\beta}^{k}(x_{i})
			+\frac{1}{2}\|x_{i}-x_{i}^{k}\|_{H_{i}}^{2},
			\\
			&\lambda^{k+1} = \lambda^{k} - \tau\beta\sum_{i=1}^{4}A_{i}x_{i}^{k+1},
		\end{aligned}\right. 	
	\end{equation}
	where $\mathcal{L}_{\beta}(\bm{x},\lambda)$ is the augmented Lagrangian function of \eqref{sec5-problem} and $\mathcal{L}_{\beta}^{k}(x_{i}) = \mathcal{L}_{\beta}(\dots,x_{i-1}^{k+1},x_{i},x_{i+1}^{k},\dots,\lambda^{k})$. $\mathcal{X}_{i} = \mathcal{M}$ for $i=1$, and $\mathcal{X}_{i} = \mathbb{R}^{n_{i}}$ for others.
	
	First, we characterize the geometries of objective functions \eqref{PGO-model-splitting} in the following assumption.
	\begin{assumption}\label{assumption1} 
		We assume that
		\begin{itemize}
			\item[(i)] $f$ and $g$ are both proper and lower semicontinuous functions and bounded from below in the feasible region, i.e., $f^{*}=\inf_{x_{i}\in \mathcal{X}_{i}} f>-\infty $ and $ g^{*}=\inf _{x_{i}\in \mathcal{X}_{i}} g>-\infty$.
			\item[(ii)] The gradient of $f(\bm{x})$ is Lipschitz continuous on bounded subset of $\mathbb{R}^{n_1} \times \mathbb{R}^{n_2} \times \mathbb{R}^{n_3} \times \mathbb{R}^{n_4}$ with Lipschitz constant $L_{f}>0$, i.e., for any $\bm{x}$ and $\bm{y}\in \mathbb{R}^{n_1} \times \mathbb{R}^{n_2} \times \mathbb{R}^{n_3} \times \mathbb{R}^{n_4}$, it holds that 
			\begin{align*}
				\|\nabla f(\bm{x})-\nabla f(\bm{y})\|
				\leq L_{f}\|\bm{x}-\bm{y}\|^{2}.
			\end{align*}  
			Similarly, the gradient of $g$ is Lipschitz continuous on bounded subset of $\mathbb{R}^{n_3} \times \mathbb{R}^{n_4}$ with Lipschitz constant $L_{g}>0$.
			
			\item[(iii)]  
			For any fixed $\bm{x}_{-i}$, $i=1,\dots,4$, the partial gradient $\nabla_{i}f(\bm{x})$ is globally Lipschitz with constant $L_{f,i}(\bm{x}_{-i})>0$, that is, for any $x_{i},y_{i} \in \mathbb{R}^{n_{i}}$,
			\begin{align*}
				&\|\nabla_{i}f(x_{i},\bm{x}_{-i})
				-\nabla_{i}f(y_{i},\bm{x}_{-i})\|
				\leq L_{f,i}(\bm{x}_{-i})\|x_{i}-y_{i}\|.
			\end{align*}
			Similarly,  $\nabla_{i}g(x_{1},x_{2})$, $i=1,2$, are also globally Lipschitz with constants $L_{g,1}(x_{2})>0$ and $L_{g,2}(x_{1})>0$.
			
			\item[(iv)] 
			If $\bm{x}_{-i}$ lies in a bounded subset, then the Lipschitz constants of the partial gradient $\nabla_{i}f(\bm{x})$ and $\nabla_{i}g(\bm{x})$
			have uniform upper bounds, respectively, i.e.,
			\begin{align*}
				\sup_{x_{j}, j\neq i} L_{f,i}(\bm{x}_{-i})
				\leq L_{f,i}, \text{and} 
				\sup_{x_{j}, j\neq i} L_{g,i}(x_{-i})\leq L_{g,i}.
			\end{align*}
			
			\item[(v)]
			$f$ and $g$ are block multi-convex functions, i.e., for each $i$ with fixed $n-1$ blocks $\bm{x}_{-i}$, and for any $x_{i},y_{i} \in \mathbb{R}^{n_{i}}$,  
			\begin{align*}
				f(y_{i},\bm{x}_{-i}) \geq f(x_{i},\bm{x}_{-i}) + 
				\left\langle \nabla_{i}f(x_{i},\bm{x}_{-i}), y_{i}-x_{i}\right\rangle.  
			\end{align*}
		\end{itemize}
	\end{assumption}
	\begin{remark}
		Assumption 1(ii) implies (iii), but (iii) may admit a smaller constant. Since $f$ and $g$ in PGO model \eqref{PGO-model-splitting} are polynomial functions, it is straightforward to verify that they satisfy this assumption \cite{chen2024pieadmm}.
	\end{remark}
	
	\subsection{Subsequential Convergence}
	Our analysis circumvents a fundamental limitation in classical nonconvex ADMM convergence analysis. Existing works typically require that $A_4$ is surjective or $Im([A_1,A_2,A_3]) \subset Im(A_4)$ to bound the dual variable residuals (see \cite{boct2020proximal,hien2022inertial,wang2019global}). In contrast, we define $A_4$ in \eqref{sec5-problem} as injective which violates this standard assumption, necessitating a novel convergence analysis. Our key innovation lies in exploiting the invertibility of the combined matrix $[A_2,A_4]$---designed to facilitate alternating updates of $x_{2}$ and $x_{4}$ rather than conventional simultaneous updates. This novel framework enables us to bound the iterative residual $\|\lambda^{k+1}-\lambda^{k}\|$. We denote $\Delta x_{i}^{k}=x_{i}^{k}-x_{i}^{k-1}$, $\Delta \bm{x}^{k}=\bm{x}^{k}-\bm{x}^{k-1}$, and $\Delta \lambda^{k}=\lambda^{k}-\lambda^{k-1}$. 
	This analysis framework follows the  ``bounded dual by primal'' technique. 
	\begin{lemma}\label{lem-5.1}
		Suppose that Assumption \ref{assumption1} holds. Let $\{(\bm{x}^{k},\lambda^{k})\}$ be the sequence generated by \eqref{sec5-ADMM} which is assumed to be bounded, then
		\begin{align}
			&\|\Delta \lambda^{k+1}\|^{2}
			\leq ~\alpha_1\left( \|\Delta \lambda^{k}\|^{2} -\|\Delta \lambda^{k+1}\|^{2}\right)
			\nonumber\\
			&+\tau \alpha_2
			\left(
			\|\Delta x_{1}^{k+1}\|_{7L_{f}^{2}+4L_{g}^{2}}^2
			+\|\Delta x_{3}^{k+1}\|_{3L_{f}^{2}}^2
			+\|\Delta x_{3}^{k}\|_{4L_{f}^{2}}^2
			\right. 
			\nonumber\\	 
			&+\left.
			\|\Delta x_{2}^{k+1}\|_{7L_{f}^{2}I+4L_{g}^{2}I+4H_{2}^{\top}H_{2}}^2
			+\|\Delta x_{2}^{k}\|_{4H_{2}^{\top}H_{2}}^2
			\right. 
			\nonumber\\	 
			&+\left.
			\|\Delta x_{4}^{k+1}\|_{3L_{f}^{2}I+3H_{4}^{\top}H_{4}}^2
			+\|\Delta x_{4}^{k}\|_{4L_{f}^{2}I+3H_{4}^{\top}H_{4}}^2
			\right).
			\label{eq-bounded lambda}
		\end{align}
		where 
		$$\alpha_1 = \frac{|1-\tau|}{1-|1-\tau|},
		\quad
		\alpha_2 = \frac{\tau}{(1-|1-\tau|)^2}. 
		$$
	\end{lemma}
	Next, we construct a Lyapunov function 
	\begin{align}\label{def-Psi}
		\Psi^{k}:
		&=\Psi(\bm{x}^{k},\bm{x}^{k-1},\lambda^{k},\lambda^{k-1})\nonumber\\
		&=\mathcal{L}_{\beta}(\bm{x}^{k},\lambda^{k})
		+\frac{\alpha_{1}}{\tau\beta}\|\Delta\lambda^{k}\|^{2}
		+\sum_{i=1}^4 \|\Delta x_{i}^{k}\|_{M_{i}}^{2},
	\end{align}
	to prove the subsequential convergence.
	\begin{theorem}\label{thm-decrease}
		Suppose that Assumption \ref{assumption1} holds. Let $\{(\bm{x}^{k},\lambda^{k})\}$ be the sequence generated by \eqref{sec5-ADMM} which is assumed to be bounded. There exists a constant $\underline{\beta}>0$, such that for all $\beta>\underline{\beta}$, we have
		\begin{itemize}
			\item[(i)]  $\{\Psi^{k}\}$ is nonincreasing, i.e.,
			\begin{align}
				\Psi^{k}-\Psi^{k+1}\geq \sum_{i=1}^4 (1-\delta_{i})\|\Delta x_{i}^{k}\|_{M_{i}}^{2},
				\label{eq-thm5.2-L1}
			\end{align}
			with $0<\delta_{i}<1$, and the right-hand side is nonnegative.
			\item[(ii)] the sequences $\{\Delta x_{i}^{k}\}$ and $\{\Delta \lambda^{k}\}$ converge to $0$.
			\item[(iii)] every limit point of the generated sequence $\{(\bm{x}^{k},\lambda^{k})\}$ is a stationary point of $\mathcal{L}_{\beta}$.
		\end{itemize}
	\end{theorem}
	
	\begin{table*}[t]
		\centering
		\begin{tabular}{l|ccc|ccc|ccc}
			\toprule[1pt]
			\multicolumn{1}{c}{}& 
			\multicolumn{3}{c}{$\sigma_{r}=0.01$, $\sigma_{t}=0.01$} & 
			\multicolumn{3}{c}{$\sigma_{r}=0.03$, $\sigma_{t}=0.05$} & 
			\multicolumn{3}{c}{$\sigma_{r}=0.05$, $\sigma_{t}=0.1$} \\
			\cmidrule(lr){2-4}\cmidrule(lr){5-7}\cmidrule(lr){8-10}
			\multicolumn{1}{c}{Algorithm} & Rel. Err. & NRMSE & \multicolumn{1}{c}{Time (s)} & Rel. Err. & NRMSE & \multicolumn{1}{c}{Time (s)} & Rel. Err. & NRMSE & Time (s) \\
			\midrule[0.5pt]
			mG-N & 0.0711 & 0.0354 & 0.407 & 0.3683 & 0.1834 & 0.401 & 0.5024 & 0.2502 & 0.407 \\
			SE-sync & 0.0711 & 0.0354 & 0.179 & 0.3683 & 0.1834 & 0.166 & 0.5025 & 0.2502 & 0.161 \\
			RS+PS & 0.0699 & 0.0348 & 0.069 & 0.3457 & 0.1721 & 0.065 & 0.4861 & 0.2420 & 0.072 \\
			RGD & 0.0692 & 0.0344 & 0.322 & 0.3087 & 0.1537 & 0.355 & 0.4729 & 0.2355 & 0.354 \\
			SOC & 0.0691 & 0.0344 & 0.647 & 0.3085 & 0.1536 & 0.425 & 0.4729 & 0.2354 & 0.436 \\
			PieADMM & 0.0689 & 0.0343 & 0.123 & 0.3085 & 0.1536 & 0.131 & 0.4729 & 0.2354 & 0.215 \\
			PRADMM (ours) & \textbf{0.0689} & \textbf{0.0343} & \textbf{0.065} & \textbf{0.3085} & \textbf{0.1536} & \textbf{0.034} & \textbf{0.4724} & \textbf{0.2352} & \textbf{0.042} \\
			\bottomrule[1pt]
		\end{tabular}
		\caption{Numerical results of different noise levels of circular ring datasets with $m=n=100$.}
		\label{tab:circular ring_results}
	\end{table*}
	
	\subsection{Global Convergence}
	Since the multiplication of two quaternion variables is a polynomial, the proposed model \eqref{PGO-model-splitting} satisfies the Riemannian Kurdyka--\L ojasiewicz property (Theorem 3 in  \cite{huang2022riemannian}), which helps us establish the global convergence of \eqref{sec5-ADMM}. First, we prove the following bound for $\operatorname{grad} \Psi^{k+1}$.
	\begin{lemma}\label{lem-subgradient-bound}
		Suppose that Assumption \ref{assumption1} holds. Let $\{(\bm{x}^{k},\lambda^{k})\}$ be the sequence generated by \eqref{sec5-ADMM} which is assumed to be bounded. If $\beta>\max\{\beta^{\prime},\beta^{\prime\prime}\}$, then there exist some constants $\varrho_1,\varrho_2,\varrho_3>0$ such that 
		\begin{align*}
			\|\operatorname{grad} \Psi^{k+1}\|\leq \varrho_1\|\Delta \bm{x}^{k+1}\|+\varrho_2\|\Delta \bm{x}^{k}\|+\varrho_3\|\Delta \lambda^{k+1}\|.
		\end{align*} 
	\end{lemma}
	
	Next, we show the whole sequence convergence. Let $\hat{z}:=(\hat{\bm{x}},\hat{\bm{x}},\hat{\lambda},\hat{\lambda})$ be a limit point of the sequence $\{z^{k}:=(\bm{x}^{k},\bm{x}^{k-1},\lambda^{k},\lambda^{k-1})\}$.  The set of all limit points of $\{z^{k}\}$ is denoted by
	\begin{align*}
		\Omega^{*}:=
		\left\lbrace \hat{z}\mid 
		\exists~\{z^{k_j}\}~\text{such that}~
		z^{k_j}\rightarrow \hat{z}~
		\text{as}~j \rightarrow +\infty \right\rbrace. 
	\end{align*}
	
	\begin{theorem}\label{thm-global-convergence}
		(A finite length property)
		Suppose that Assumption \ref{assumption1} holds. Let $\{(\bm{x}^{k},\lambda^{k})\}$ be the sequence generated by \eqref{sec5-ADMM} which is assumed to be bounded. If $\Psi$ satisfies the Riemannian K\L~  property at every point in  $\Omega^{*}$,  and $\beta>\max\{\beta^{\prime},\beta^{\prime\prime}\}$, then
		\begin{itemize}
			\item[(i)] The sequence $\{(\bm{x}^{k},\lambda^{k})\}$ has finite length, that is 
			\begin{align}
				\sum_{k=1}^{+\infty} \|\bm{x}^{k+1}-\bm{x}^{k}\|+\|\lambda^{k+1}-\lambda^{k}\|<+\infty.
				\label{eq-thm-gc-1}
			\end{align}
			\item[(ii)] The sequence $\{(\bm{x}^{k},\lambda^{k})\}$ converges to a stationary point of $\mathcal{L}_{\beta}$.
		\end{itemize}	
	\end{theorem}

	\begin{remark}\label{remark-KL}
		While the global convergence of nonconvex optimization via K\L~  functions has been extensively studied in \cite{bolte2014proximal,pock2016inertial}, these results cannot be directly extended to our algorithm due to their exclusion of dual variables $\lambda$. Prior work \cite{hien2022inertial,hien2024inertial} established global convergence for ADMM with over-relaxation at $\tau=1$, but convergence guarantees for $0<\tau<2$ remained an open challenge. The core difficulty lies in bounding $\|\Delta \lambda^{k+1}\|$ rather than its squared counterpart $\|\Delta \lambda^{k+1}\|^{2}$. 
		
		Although our analysis in Lemma \ref{lem-5.1} demonstrates that bounding $\|\Delta \lambda^{k+1}\|^{2}$ inherently includes the term $\|\Delta \lambda^{k}\|^{2}-\|\Delta \lambda^{k+1}\|^{2}$. This term vanishes when $\tau=1$ but becomes intractable in our proof of global convergence with $0<\tau<2$. To overcome this persistent challenge, we develop a novel upper bound for $\|\Delta \lambda^{k+1}\|$ which enables term-by-term cancellation throughout the summation analysis. This theorem extends prior ADMM frameworks and resolves a key limitation in relaxation parameter selection,  while preserving validity on Riemannian manifolds.
	\end{remark}
	
	\section{Numerical Experiments} \label{sec-numerical experiments}
	We test the effectiveness of Algorithm \ref{alg1} on different 3D pose graph datasets. As a basis for comparison, we also evaluate the performance of several state-of-the-art pose-graph SLAM approaches, including the manifold-based Gauss-Newton (mG-N) method \cite{wagner2011rapid}, SE-sync method \cite{rosen2019se}, pose synchronization  (RS+PS) algorithm \cite{nasiri2020novel}, and PieADMM \cite{chen2024pieadmm}. The CPU time of SE-Sync is computed without the time spent on the optimality check. In addition, we complement our test experimental evaluation encompassing both the standard Riemannian gradient descent \cite{gabay1982minimizing} with line search and specialized manifold splitting methods, such as SOC \cite{lai2014splitting}. All experiments are performed using MATLAB 2020a 
	on an Intel i7-10700F CPU desktop computer with 16GB of memory. More details of all experiments can be found in the appendix.
	
	\subsection{Synthetic Datasets}
	\textbf{Data Settings.} We evaluate the performance on two synthetic datasets:
	\begin{itemize}
		\item[(a).] Circular ring: this is a single loop with a radius of 2 and odometric edges. The closed-loop constraint requires the first point to coincide with the last point. Recovering the trajectory poses is particularly challenging because of the sparse observations, i.e., $m=n$.
		\item[(b).] Cube dataset: assume that a robot travels on a  $2 \times 2\times 2$ grid world and random loop closures are added between nearby nodes with probability $p_{cube}$. The total number of vertices is $n=\hat{n}^{3}$, where $\hat{n}$ is the number of nodes on each side of the cube, and
		the expectation of the number of edges is $\mathbb{E}(m)=2(2\hat{n}^{3}-3\hat{n}^2+1)p_{cube} + \hat{n}^3 - 1$.
	\end{itemize}
	The noisy relative pose measurements are generated by
	\begin{equation*}
		\begin{aligned}
			\bm{t}_{ij}&=R_{i}^{\top}(\bm{t}_{j}-\bm{t}_{i})+\bm{t}_{\epsilon},  ~ \text{where}~ \bm{t}_{\epsilon} \sim   \mathcal{N}(0,\sigma_{t}^2I_{3}),\\
			\tilde{q}_{ij}&=\tilde{q}_{i}^{*}\tilde{q}_{j}\tilde{q}_{\epsilon}, \quad\quad\quad\quad~~ \text{where}~\tilde{q}_{\epsilon}\sim   \text{vMF}([1,0,0,0],\frac12\sigma_{r}^2).
		\end{aligned}
	\end{equation*}
	
	\begin{figure*}[thbp]
		\centering
		% \subfloat{\label{cube_03_4_trajectory}	\includegraphics[width=1\linewidth]{fig//cube_03_4_trajectory.eps}}
		
		\subfloat{\label{cube_03_6_trajectory}	\includegraphics[width=1\linewidth]{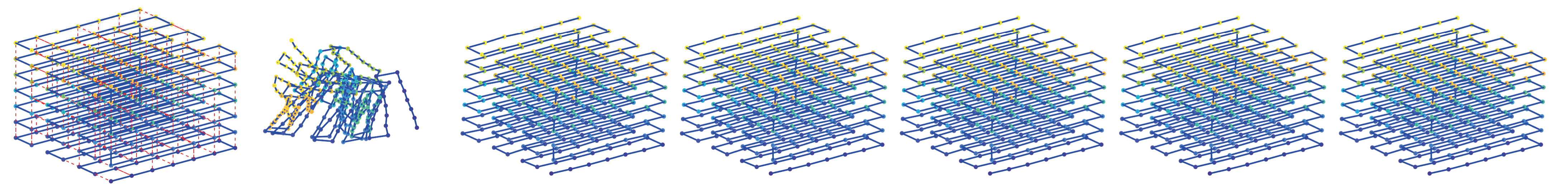}}
		
		{\raggedright \quad~~ Ground truth\qquad\quad Noisy  ~\quad\qquad RE: 0.1782 \qquad~ RE: 0.1783 \qquad~ RE: 0.1487 \qquad~ RE: 0.1502 \qquad~ RE: 0.1478 \qquad\\}
		% Ground truth \quad\quad 9.169 \quad mG-N:0.1782 \quad SE-sync:0.1783 \quad RS+PS:0.1487 \quad PieADMM:0.1502 \quad PRADMM:0.1478
		\caption{The comparison of cube trajectory with $\hat{n}=7$. From left to right are the real trajectory, the corrupted trajectory, and the recovered results by mG-N, SE-sync, RS+PS, PieADMM, and PRADMM, respectively.}
		\label{fig:cube_trajectory}
	\end{figure*}
	
	\noindent \textbf{Experiment Settings.} In our experiments, we set the noise level of the translation part $\sigma_{t} \in (0.01,0.1)$, and the magnitude of rotation noise with $\sigma_{r} \in (0.01,0.1)$. The relaxation stepsize $\tau=1.4$.  We estimate the upper bounds of $L_{f}$ and $L_{g}$ via the principal block in the Hessian matrices of $f$ and $g$. In addition, we find an upper bound of $L_{g,\tilde{\bm{q}}}$ using $\|\tilde{p}_{i}\|=1$. Empirically, we set $H_{1}\in\{0.1,1,10\}$ and $H_{i}=0.001$ for $i=2,3,4$. Then we adjust $\beta$ around the range and check the descent property of the merit function defined in Theorem \ref{thm-decrease}. We also set $tol=10^{-4}$ and $MaxIter=300$ for PRADMM. All algorithms use  chordal to find an initial point. Results are averaged over 5 runs.
	
	\noindent \textbf{Experiment Results.} The numerical results of circular ring datasets with different noise levels are listed in Table \ref{tab:circular ring_results}, which demonstrate that PRADMM achieves both faster computation and higher accuracy on small-scale circular ring datasets. 
	%Notably, RGD and SOC exhibit prohibitively growing computation times versus other algorithms with scale, due to non-parallelization or ill-fitting decomposition methods. 
	This validates the significance of developing PRADMM for model \eqref{PGO-model-splitting}. Due to non-parallelization or ill-fitting decomposition issues, we exclude RGD and SOC from our subsequent large-scale experiments. 
	%due to their computational inefficiency.
	
\begin{figure}[!tbp]
		\centering
		\subfloat{\label{cube_03_time}	\includegraphics[width=0.5\linewidth]{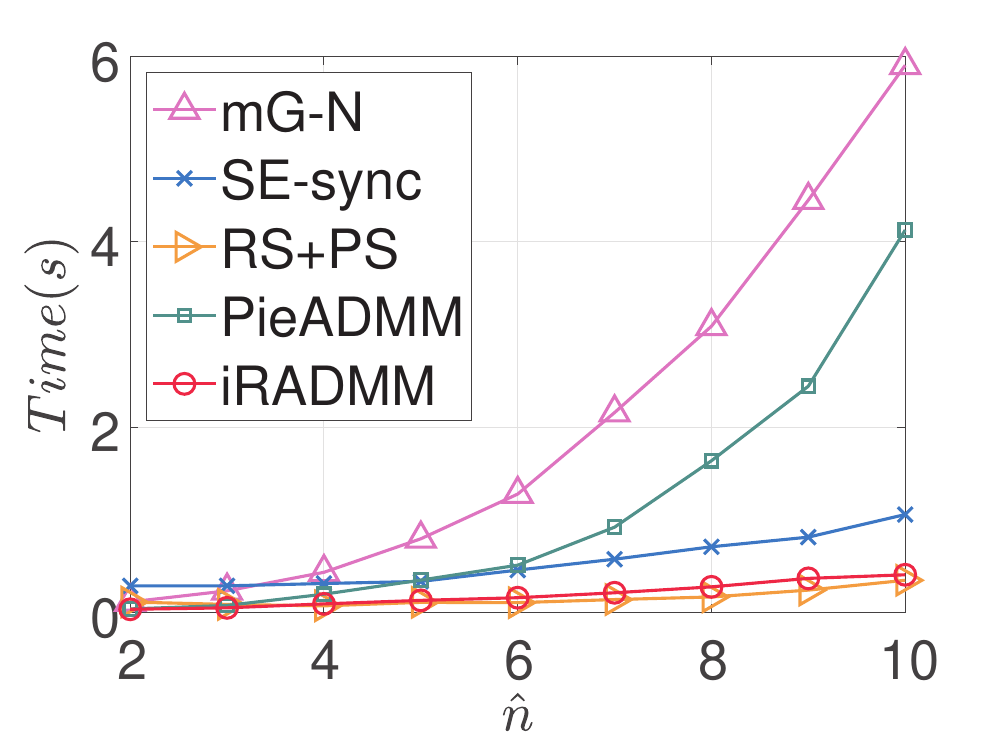}}
		\subfloat{\label{cube_09_time}	\includegraphics[width=0.5\linewidth]{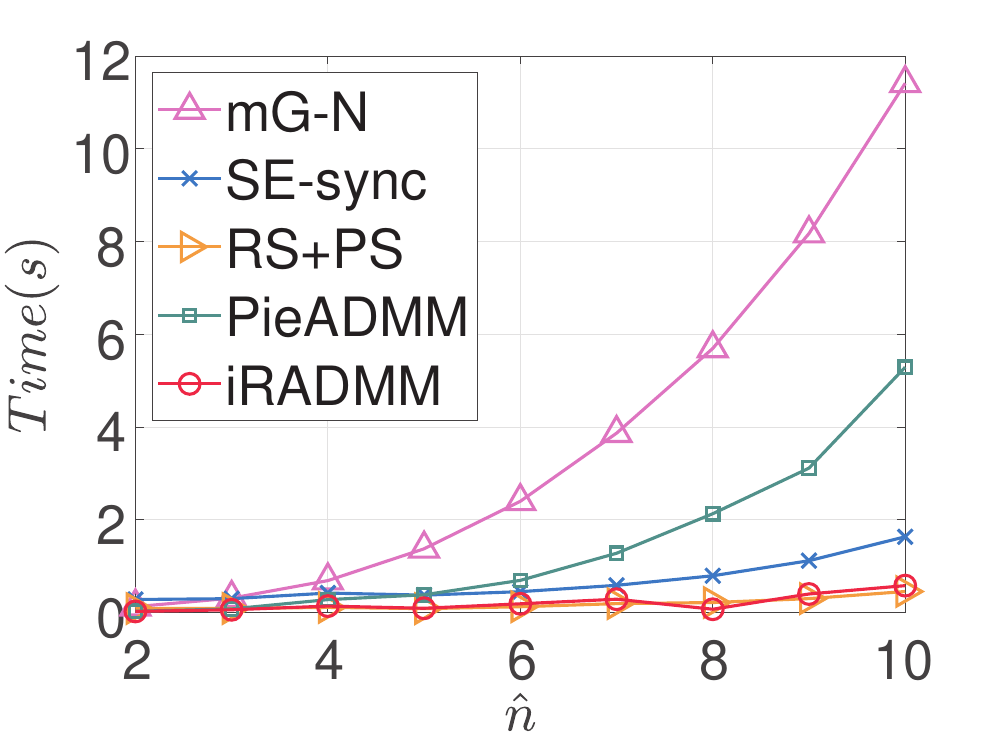}}
		
		\caption{The trend of CPU time along with different $\hat{n}$ for cube datasets. Left: $p_{cube}=0.3$. Right: $p_{cube}=0.9$.
			%The last three rows correspond to $p_{cube}=\{0.3,0.6,0.9\}$, respectively.
		}
		\label{fig:cube_time}
	\end{figure}

	For cube datasets, we evaluate $\hat{n}$ values ranging from $2$ to $10$ under relative noise conditions $\sigma_{t}=0.1$ and $\sigma_{r}=0.1$. Figure \ref{fig:cube_trajectory} visually compares reconstructed trajectories for $\hat{n}=7$ with $p_{cube}=0.3$, demonstrating comparable reconstruction quality across the tested algorithms. The results under the parameters $\hat{n}=7$ and $p_{cube}=0.3$ are listed in Table \ref{tab:cube_results}.   PRADMM achieves superior reconstruction accuracy compared to alternative algorithms. For large-scale PGO problems, PRADMM maintains competitive efficiency -- operating marginally slower than RS+PS yet substantially outperforming all other methods in computation speed. Figure \ref{fig:cube_time} demonstrates PRADMM's stable efficiency with increasing graph size, indicating stronger scalability for large-scale applications.
	
	\begin{table}[tbp]
		\centering
		\begin{tabular}{l|ccc}
			\toprule[1pt]
			Algorithm & Rel. Err. & NRMSE & \multicolumn{1}{c}{Time (s)} \\
			\midrule[0.5pt]
			mG-N & 0.1782 & 0.2057 & 2.157 \\
			SE-sync & 0.1783 & 0.2059 & 0.575 \\
			RS+PS  & 0.1487 & 0.1717 & 0.139 \\
			PieADMM & 0.1502 & 0.1735 & 0.920 \\
			PRADMM (ours) & 0.1478 & 0.1707 & 0.213 \\
			\bottomrule[1pt]
		\end{tabular}
		\caption{Numerical results of $\hat{n}=7$ and $p_{cube}=0.3$ for cube datasets with $m=697$.}
		\label{tab:cube_results}
	\end{table}
	
	% \begin{figure*}[!tbp]
	% 	\begin{minipage}{.26\linewidth}
	% 		\centering
	% 		\subfloat[$n$ and $m$]{\label{cube_numedge}	
	% 			\includegraphics[width=1\linewidth]{fig//cube_numedge.eps}}
	% 	\end{minipage}
	% 	\begin{minipage}{.74\linewidth}
	% 		\centering
	% 		\subfloat[$p_{cube}=0.3$]{\label{cube_03_re}	\includegraphics[width=0.33\linewidth]{fig//cube_03_re.eps}}
	% 		\subfloat[$p_{cube}=0.6$]{\label{cube_06_re}	\includegraphics[width=0.33\linewidth]{fig//cube_06_re.eps}}	
	% 		\subfloat[$p_{cube}=0.9$]{\label{cube_09_re}	\includegraphics[width=0.33\linewidth]{fig//cube_09_re.eps}}
	
	% 		\subfloat[$p_{cube}=0.3$]{\label{cube_03_time}	\includegraphics[width=0.33\linewidth]{fig//cube_03_time.eps}}
	% 		\subfloat[$p_{cube}=0.6$]{\label{cube_06_time}	\includegraphics[width=0.33\linewidth]{fig//cube_06_time.eps}}
	% 		\subfloat[$p_{cube}=0.9$]{\label{cube_09_time}	\includegraphics[width=0.33\linewidth]{fig//cube_09_time.eps}}
	% 	\end{minipage}
	% 	\caption{The trend of the number of edges, relative error, and CPU time along with different $\hat{n}$ under $\sigma_{t}^{rel}=0.1$, $\sigma_{r}=0.1$ and $p_{cube}\in\{0.3,0.6,0.9\}$.}
	% 	\label{fig:cube_time}
	% \end{figure*}
	
	\subsection{SLAM Benchmark Datasets}
	The evaluations of this part are performed on the widely-used benchmark datasets \cite{rosen2019se}. We  use chordal initialization technique for all methods. 
	A portion of existing works  \cite{liu2012sdr,fan2019generalized,nasiri2020novel} assesses the precision of solutions by directly comparing rotational angles and translational distances between estimated poses, which quantify error magnitudes across different metrics. The results of the sphere datasets are shown in Table \ref{benchmark_table}, which indicates that PRADMM is efficient in terms of CPU time and translational errors. For more experimental results from other benchmark datasets, please refer to the appendix.
	%converges faster than other algorithms. In terms of accuracy, PRADMM achieves smaller translational errors.
	
	\begin{table}[!tbp]
		\centering
		\begin{tabular}{l|ccc}
			\toprule[1pt]
			Algorithm & Loss$(\theta)$ & Loss$(\bm{t})$ & Time\,(s)\\
			\midrule 
			mGN & 1.51e+06  & 8.94e+03 & 23.739 \\
			SE-sync & 1.49e+06  & 9.23e+03 & 1.142 \\
			RS+PS & 1.49e+06  & 9.23e+03 & 0.987 \\
			PieADMM & 1.63e+06 & 7.63e+03 & 5.697 \\
			PRADMM (ours) & 1.63e+06  & 8.13e+03 & 0.387 \\
			\bottomrule[1pt]
		\end{tabular}
		\caption{The numerical results of errors in rotation angle and translation, along with the CPU time consumed for sphere datasets with $n=2200$, $m=8647$.}
		\label{benchmark_table}
	\end{table}
	
	\section{Conclusion}
	This paper presents PRADMM, a novel algorithm for large-scale PGO. We   decouple   vertex correlations via mathematical reformulation, enabling fully parallelizable vertex updates and remarkably stable computation times regardless of graph size. Crucially, we establish global convergence under significantly weaker conditions (row rank-deficient coefficients) than existing requirements, while permitting extended dual steps $(\tau\in(0,2))$, surpassing conventional limits on the manifolds. 
	Extensive experiments show PRADMM consistently outperforms state-of-the-art graph-SLAM methods.
	%Extensive experimental results on synthetic and real-world 3D SLAM benchmarks demonstrate that PRADMM consistently outperforms state-of-the-art methods, converging faster while maintaining high solution accuracy, making it highly suited for large-scale robotic navigation.
	
	\section*{Acknowledgments}
	The authors are grateful to the anonymous referees for their careful reading of the
	manuscript and their valuable comments. This work was supported by the National Natural Science Foundation of China (Nos. 12471282 and 12131004), Research Center for Intelligent Operations Research, The Hong Kong Polytechnic University (4-ZZT8),
	and the Fundamental Research Funds for the Central Universities (No. YWF-22-T-204).
	
	\bibliography{ref_SLAM,ref_ADMM,ref_opt,ref_probability,ref_others}
	

% \clearpage
% \newpage
% \input{ReproducibilityChecklist}
%
\clearpage
\newpage

\onecolumn
\begin{center}
\fontsize{15}{17}\selectfont \textbf{Appendix}
\end{center}
~\\

\noindent The appendices are organized as follows:\\
In Appendix \ref{app-preliminaries}, we give several  basic definitions, notations, and lemmas for the following analysis.\\
In Appendix \ref{app-subproblems}, we show the detailed derivations for the subproblems.\\
In Appendix \ref{app-inertial Riemannian ADMM}, we prove Theorem \ref{TRS-subproblem3} for the equivalence of $\tilde{\bm{p}}$-subproblem and generalized eigenvalue problem.\\
In Appendix \ref{app-Convergence Analysis}, we prove the main lemmas and theorems for the convergence of PRADMM.\\
In Appendix \ref{app-numerical experiments}, we provide the experimental details and additional experiments.\\

\appendix
\section{Notation and Preliminaries}  \label{app-preliminaries}
\subsection{Basic Notations}
The fields of real numbers, quaternion numbers, and unit quaternion numbers are denoted by $\mathbb{R}$, $\mathbb{Q}$, and $\mathbb{U}$, respectively. 
Throughout this paper, scalars, vectors, matrices, and quaternions are denoted by lowercase letters
(e.g., $x$), boldface lowercase letters (e.g., $\bm{x}$), boldface capital letters (e.g., $X$), and lowercase letters with tilde (e.g., $\tilde{x}$), respectively.

We denote by $\mathbb{S}^{n \times n}$ the set of all symmetric matrices. The notation $\left\| \cdot\right\|$ denotes the $2$-norm of vectors or the Frobenius norm of matrices. Let $M$ be a positive definite linear operator; we use $\|\bm{x}\|_{M}:=\sqrt{\left\langle \bm{x},M\bm{x}\right\rangle }$ to denote its $M$-norm; and $\sigma_{\min}(M)$ and $\sigma_{\max}(M)$  denote the smallest and largest eigenvalue
of $M$, respectively. For symmetric matrices $M_{1}, M_{2} \in \mathbb{R}^{n \times n}$,
$M_{1}\succ M_{2}$ and $M_{1}\succeq M_{2}$ means that $M_{1}- M_{2}$ is positive definite and positive semidefinite, respectively.

\subsection{Quaternion and Pose}
A quaternion number $\tilde{q} \in \mathbb{Q}$, proposed by Hamilton, has the form $\tilde{q}=q_0+q_1 \mathbf{i}+q_2 \mathbf{j}+q_3 \mathbf{k}$, where $q_0, q_1, q_2, q_3 \in \mathbb{R}$ and $\mathbf{i}, \mathbf{j}, \mathbf{k}$ are three imaginary units.
We may also write $\tilde{q}=[q_0,q_1,q_2,q_3]=[q_0, \bm{q}] \in \mathbb{R}^4$ as the vector representation where $\bm{q}=[q_1,q_2,q_3]\in \mathbb{R}^3$ for convenience. We note that we also regard the above representation as a column
vector and its transpose $[q_0, \bm{q}]^\top$ a row vector.
The sum of $\tilde{p}$ and $\tilde{q}$ is defined as
$
\tilde{p}+\tilde{q}=[p_0+q_0,\bm{p}+\bm{q}].
$
The product of $\tilde{p}$ and $\tilde{q}$ is defined by
$$
\tilde{p} \tilde{q}=\left[p_0 q_0-\bm{p} \cdot \bm{q}, p_0 \bm{q}+q_0 \bm{p}+\bm{p} \times \bm{q}\right],
$$
where $\bm{p} \cdot \bm{q}$ is the $\operatorname{dot}$ product, and $\bm{p} \times \bm{q}$ is the cross product of $\bm{p}$ and $\bm{q}$. Thus, in general, $\tilde{p} \tilde{q} \neq \tilde{q} \tilde{p}$, and we have $\tilde{p} \tilde{q}=\tilde{q} \tilde{p}$ if and only if $\bm{p} \times \bm{q}=\vec{0}$, i.e., either $\bm{p}=\vec{0}$ or $\bm{q}=\vec{0}$, or $\bm{p}=\alpha \bm{q}$ for several real number $\alpha$ (see \cite{qi2023Augmented}). The multiplication of quaternions is associative and distributive over vector addition, but is not commutative.

The conjugate of $\tilde{q}$ is the quaternion $\tilde{q}^*=q_0-q_1 \mathbf{i}-q_2 \mathbf{j}-q_3 \mathbf{k}$. Then, $(\tilde{p} \tilde{q})^*=\tilde{q}^* \tilde{p}^*$ for any $\tilde{p}, \tilde{q} \in \mathbb{Q}$.  The magnitude of $\tilde{q}$ is defined by 
$
\|\tilde{q}\|=\sqrt{\tilde{q} \tilde{q}^*}=\sqrt{\tilde{q}^* \tilde{q}}.
$
And $\tilde{q}$ is invertible if and only if $\|\tilde{q}\|$ is positive. In this case, we have
$
\tilde{q}^{-1}=\tilde{q}^{*}/\|\tilde{q}\|.
$

The set of all unit quaternions is $\mathbb{U}:=\{\tilde{q}\in \mathbb{R}^{4}\mid\|\tilde{q}\|=1\}$, which can be regarded as a unit sphere in $\mathbb{R}^{4}$. Equivalently, a unit quaternion has the following form:  
$
\tilde{q}=[\cos(\theta/2),\sin(\theta/2)\bm{n}],
$
where $\bm{n}=(n_{x},n_{y},n_{z})$ is a unit vector and $\theta$ is an angle. Let a vector $\bm{t}_{1}\in \mathbb{R}^{3}$ rotates $\theta$ radians around axis $\bm{n}$ to reach $\bm{t}_{2}\in \mathbb{R}^{3}$. This process can be represented by a quaternion as $$
[0,\bm{t}_{2}]=\tilde{q}[0,\bm{t}_{1}]\tilde{q}^{*}.
$$ 
Using rotation matrix in $SO(3)$, we also have 
$\bm{t}_{2}=R\bm{t}_{1}$, where 
$$
R =\cos(\theta)I_{3}+(1-\cos(\theta))\bm{n}\bm{n}^{\top}+\sin(\theta)\bm{n}^{\land}, \quad
\text{and} \quad
\bm{n}^\land=\left(\begin{array}{ccc}
	0 & -n_z & n_y  \\
	n_z & 0 & -n_x\\
	-n_y & n_x & 0 
\end{array}\right).
$$

\begin{table}[t]
	\centering
	\footnotesize
	\caption{Rotation and pose representations in SLAM. ``Y'' and ``N'' represent ``Yes'' and ``No'', respectively. 
    }
		\begin{tabular}{|c|c|c|c|c|}
			\Xhline{1pt}
			Rotation representation & Definition  & Storage & Constraints & Singularity\\
			\hline
			Euler angles  &   $(\phi,\theta,\varphi)$            & 3      & N       &    Y         \\
			Axis-angle  &   $so(3):=\{\bm{n}\in \mathbb{R}^{3}\}$   & 3              &    N    & Y      \\
			Special orthogonal group  &   $SO(3):=\{R \in \mathbb{R}^{3 \times 3} \mid R^\top R = I_{3}, \operatorname{det}(R)=1\}$      & 9              &    Y    & N       \\
			Unit quaternion  &   $\mathbb{U} :=\{\tilde{q}\in \mathbb{R}^{4}\mid\|\tilde{q}\|=1\}$             & 4         &    Y    & N    \\ 
			\Xhline{1pt}
			Pose representation & Definition & Storage& Constraints & Singularity \\
			\hline
			Lie algebra   &        $se(3):=\{(\bm{\rho},\bm{n})\in \mathbb{R}^{6}\}$    &  6&          N          &  Y           \\
                Special Euclidean group   &        $SE(3):=\{(R,\bm{t}): R\in SO(3),\bm{t}\in\mathbb{R}^{3}\}$    &12&          Y          &  N           \\
			Augmented unit quaternion   &        $\mathbb{AU}:=\{(\tilde{q},\bm{t})\mid \tilde{q}\in\mathbb{U},\bm{t}\in\mathbb{R}^{3}\}$    &  7&          Y          &  N           \\
			Dual quaternion   &        $\hat{\mathbb{U}}:=\{(\tilde{q}_r,\tilde{q}_d)\in\mathbb{R}^{8}\mid\tilde{q}_r\in\mathbb{U}, \left\langle\tilde{q}_r,\tilde{q}_d\right\rangle=0\}$    &  8&          Y          &  N           \\
			\Xhline{1pt}
	\end{tabular}
	\label{app-table-pose-representation}
\end{table}

We list the rotation and pose representations in SLAM in Table \ref{app-table-pose-representation}. If the rotation matrix $R$ is compound motion of two rotations, i.e., $R=R_{2}R_{1}$, then the corresponding quaternion $\tilde{q}$ can be formulated as $\tilde{q}=\tilde{q}_{2}\tilde{q}_{1}$. 
Next, we show a lemma that can simplify the product of two quaternions by multiplication between matrix and vector.
Given any $\tilde{a}=[a_0,a_1,a_2,a_3]$, we define 
$$
M(\tilde{a})=\left(\begin{array}{rrrr}
	a_0 & -a_1 & -a_2 & -a_3 \\
	a_1 & a_0 & -a_3 & a_2 \\
	a_2 & a_3 & a_0 & -a_1 \\
	a_3 & -a_2 & a_1 & a_0
\end{array}\right), \qquad 
W(\tilde{a})=\left(\begin{array}{rrrr}
	a_0 & -a_1 & -a_2 & -a_3 \\
	a_1 & a_0 & a_3 & -a_2 \\
	a_2 & -a_3 & a_0 & a_1 \\
	a_3 & a_2 & -a_1 & a_0
\end{array}\right) .
$$
\begin{lemma}
	\cite{chen2024regularization} For any $\tilde{a} = [a_0,a_1,a_2,a_3]\in \mathbb{Q}$ and $\tilde{b} = [b_0,b_1,b_2,b_3]\in \mathbb{Q}$, the following statements hold
	\begin{enumerate}
		\item[(a)] $M(\tilde{a}^*)=M(\tilde{a})^\top$, $W(\tilde{a}^*)=W(\tilde{a})^\top$.
		\item[(b)] $\tilde{a}\tilde{b}=M(\tilde{a})\tilde{b}=W(\tilde{b})\tilde{a}$.
		\item[(c)] $M(\tilde{a})^\top M(\tilde{a})=W(\tilde{a})^\top W(\tilde{a})=\|\tilde{a}\|^2 I_{4}$, where $I_{4}$ is the identity matrix of size $4 \times 4$.
	\end{enumerate}
	\label{lem-MW}
\end{lemma}

\subsection{Riemannian Optimization} \label{sec-Riemannian optimization}
Suppose $\mathcal{M}$ is a Riemannian manifold \cite{absil2009optimization}. For all $\bm{x} \in \mathcal{M}$, the tangent space is denoted by $T_{\bm{x}}\mathcal{M}$. If $f$ is differentiable at $\bm{x} \in \mathcal{M}$, we represent the Euclidean and Riemannian gradients by $\nabla f(\bm{x})$ and $\operatorname{grad} f(\bm{x})$, respectively. For an $m$-dimensional Riemannian submanifold $\mathcal{M}$, we have 
$
\operatorname{grad} f(\bm{x}) = \operatorname{Proj}_{T_{\bm{x}}\mathcal{M}}(\nabla f(\bm{x})),
$
where $\operatorname{Proj}_{T_{\bm{x}}\mathcal{M}}$ is the Euclidean projection operator onto the subspace $T_{\bm{x}}\mathcal{M}$.

When we regard the unit quaternion $\mathbb{U}$ as a sphere $\mathcal{S}^{3}$ embedded in $\mathbb{R}^{4}$, the tangent space is 
$
T_{\bm{x}}\mathbb{U}=\{\bm{v}\in\mathbb{R}^{4}:\bm{x}^{\top} \bm{v}=0\}, 
$ 
and the Riemannian gradient of $f$ is 
$
\operatorname{grad}f(\bm{x})=(I_{4}-\bm{x}\bm{x}^{\top})\nabla f(\bm{x}).
$

Consider the Riemannian optimization problem 
$\min_{\bm{x}\in \mathcal{M}} f(\bm{x})$ where the first-order optimality condition requires $\operatorname{grad} f(\bm{x}^{*})=0$. When linear equality constraints are present:
$$\min ~f(\bm{x}) \quad \text{s.t.}~\bm{x}\in \mathcal{M},~A\bm{x}=b,$$
and let 
$
\mathcal{L}(\bm{x},\bm{\lambda}) = f(\bm{x})-\bm{\lambda}^{\top}(A\bm{x}-b),
$
then $\{\bm{x}^{*},\bm{\lambda}^{*}\}$ is a first-order stationary point \cite{yang2014optimality} if
\begin{align*}
	\operatorname{grad} \mathcal{L}(\bm{x}^{*},\bm{\lambda}^{*})
	=\operatorname{Proj}_{T_{\bm{x}}\mathcal{M}}(\nabla f(\bm{x}^{*})-A^{\top}\bm{\lambda}^{*})
	=0, \quad
	A\bm{x}^{*}=b, \quad
	\bm{x}^{*}\in \mathcal{M}.
\end{align*}

 Let $\sigma: \mathcal{M} \rightarrow \mathbb{R}$ be a  continuous function. For $-\infty <\eta_1<\eta_2<+\infty$, we define 
$$
\left[ \eta_1<\sigma(u)< \eta_2 \right]:=\{u\in\mathcal{M}:\eta_1<\sigma(u)< \eta_2\}. 
$$
\begin{definition} \cite{huang2022riemannian} \label{def-KLproperty}
	(Riemannian Kurdyka--\L ojasiewicz property) The function $\sigma$ is said to have the Riemannian K\L property at $\bar{u}$ if and only if there exist $\eta \in \left( \left. 0,+\infty\right] \right. $, a neighborhood $U\subset \mathcal{M}$ of $\bar{u}$, and a continuous
	concave function $\varphi:\left[ \left. 0,\eta\right) \right.\rightarrow \mathbb{R}_{+} $ such that 
	\begin{itemize}
		\item[(i)] $\varphi(0)=0$;
		\item[(ii)] $\varphi$ is $C^{1}$ on $(0,\eta)$;
		\item[(iii)] for all $s\in (0,\eta)$, $\varphi^{\prime}(s)>0$;
		\item[(iv)] for all $u \in U \cap \left[ \sigma(\bar{u})<\sigma(u)<\sigma(\bar{u})+\eta\right] $, the following inequality holds
		\begin{align*}
			\varphi^{\prime}(\sigma(u)-\sigma(\bar{u}))\operatorname{dist}(0,\hat{\partial} \sigma(u))\geq 1,
		\end{align*}
		where $\operatorname{dist}(0,\hat{\partial} \sigma(u))=\inf\{\|v\|_{u}:v\in\hat{\partial} \sigma(u)\}$ and $\hat{\partial}$ denotes the Riemannian generalized subdifferential.
	\end{itemize}
\end{definition}
The function $\varphi$ is called the desingularising function. We denote by $\mit\Phi_{\eta}$ the class of $\varphi$ which satisfies the above definitions (i), (ii) and (iii).
\begin{lemma} \cite{huang2022riemannian} \label{lem-uniformizedKL}
	(Uniformized Riemannian K\L property)
	Let $\Omega$ be a compact set and let $\sigma: \mathcal{M} \rightarrow \mathbb{R} \cup \{+\infty\}$ be a continuous function. Assume that $\sigma$ is constant on $\Omega$ and satisfies the Riemannian K\L property at each point of $\Omega$. Then, there exist $\varepsilon>0, \eta>0$ and $\varphi \in \mit\Phi_\eta$ such that for all $\bar{u}$ in $\Omega$ and all $u$ in the following intersection
	\begin{align*}
		\left\{u \in \mathcal{M}: \operatorname{dist}(u, \Omega)<\varepsilon\right\} \cap[\sigma(\bar{u})<\sigma(u)<\sigma(\bar{u})+\eta]
	\end{align*}
	one has,	
	\begin{align*}
		\varphi^{\prime}(\sigma(u)-\sigma(\bar{u})) \operatorname{dist}(0, \hat{\partial} \sigma(u)) \geq 1.
	\end{align*}
\end{lemma}

\subsection{Distance Metrics for Rotations} \label{sec-Distance Metrics}

Assume that $R_{1}, R_{2}$ are two rotation matrices, and $\tilde{q}_{1}, \tilde{q}_{2}$ are corresponding quaternion representations. Define the rotation angle $\theta_{ij}$ of two different rotations as 
$$\theta=\operatorname{dist}_{a}(R_{1},R_{2})=\|\operatorname{Log}(R_{1}^{\top}R_{2})\|
=\arccos\frac{\operatorname{tr}(R_{1}^{\top}R_{2})-1}{2}, $$
where $\operatorname{Log}(R)$ denotes the logarithm map for $SO(3)$ and return its axis-angle representation $\bm{n}$.
It can be shown that this distance is a geodesic distance \cite{Hartley2012RA}, i.e., it is the length of the minimum path between $R_{1}$ and $R_{2}$ on the manifold $SO(3)$. Then the
chordal distance and the quaternion distance between two rotations can be interpreted as 
$$
\operatorname{dist}_{c}(R_{1},R_{2})
=\|R_{1}-R_{2}\|_{F}
=2\sqrt{2}\sin(\theta/2),
$$
$$
\operatorname{dist}_{q}(R_{1},R_{2})
=\|\tilde{q}_{1}-\tilde{q}_{2}\|_{F}
=2\sin(\theta/4).
$$
\section{Detailed Derivations for the Subproblems in Section \ref{sec-algorithm}}  \label{app-subproblems}

\subsection{$\tilde{\bm{p}}$--subproblems}
For the $\tilde{\bm{p}}$-subproblem, according to Lemma \ref{lem-MW}, the problem \eqref{subproblem_p_1} can be reformulated as 
\begin{align*}
	\tilde{\bm{p}}^{k+1}  =& \mathop{\arg\min}\limits_{\tilde{\bm{p}} \in \mathbb{U}^{n}}  
	\sum_{i=1}^{n}  \left\lbrace  \sum_{(i,j)\in \mathcal{E}_{i}^{out}}
	\|M(\tilde{q}_{i}^{k})M(\tilde{t}_{ij})D\tilde{p}_{i}-(\tilde{t}_j^{k} - \tilde{s}_i^{k}) \|_{\Sigma_{1}}^2\right. \\
	& \left. + \sum_{(l,i)\in \mathcal{E}_{i}^{in}}\|W(\tilde{q}_{li})W(\tilde{q}_{l}^{k})\tilde{p}_{i}-1\|_{\Sigma_{2}}^2
	+ \frac{\beta_1}{2}\|\tilde{p}_{i}-(\tilde{q}_{i}^{k}+\frac{1}{\beta_1}\bm{\lambda}_{i}^{k})\|^{2}
	+\frac{1}{2}\|\tilde{p}_{i}-\tilde{p}_{i}^{k+\frac{1}{2}}\|_{H_{1,i}}^2\right\rbrace. 
\end{align*}
where the matrix $D=\text{diag}(1,-1,-1,-1)$ is a diagonal matrix of size $4 \times 4$. Since $\tilde{p}_{i}$ are fully separable, we can update them in parallel, i.e.,
\begin{align}
	\tilde{p}_{i}^{k+1}  =& \mathop{\arg\min}\limits_{\tilde{p}_{i} \in \mathbb{U} } 
	\sum_{(i,j)\in \mathcal{E}_{i}^{out}}
	\|G_{1,ij}^k\tilde{p}_{i}-(\tilde{t}_j^{k} - \tilde{s}_i^{k}) \|_{\Sigma_{1}}^2 
	+\sum_{(l,i)\in \mathcal{E}_{i}^{in}}\|G_{2,ij}^k\tilde{p}_{i}-1\|_{\Sigma_{2}}^2\nonumber\\
	&  + \frac{\beta_1}{2}\|\tilde{p}_{i}-(\tilde{q}_{i}^{k}+\frac{1}{\beta_1}\bm{\lambda}_{i}^{k})\|^{2}
	+\frac{1}{2}\|\tilde{p}_{i}-\tilde{p}_{i}^{k+\frac{1}{2}}\|_{H_{1,i}}^2. \nonumber\\
	= & \mathop{\arg\min}\limits_{\tilde{p}_{i} \in \mathbb{U} } ~ \frac{1}{2}\tilde{p}_{i}^{\top}A_{1,i}^{k}\tilde{p}_{i} + (b_{1,i}^{k})^{\top}\tilde{p}_{i} \nonumber
\end{align}
where $G_{1,ij}^k = M(\tilde{q}_{i}^{k})M(\tilde{t}_{ij})D$,  $G_{2,ij}^k =W(\tilde{q}_{li})W(\tilde{q}_{l}^{k})$, and
\begin{align*}
A_{1,i}^{k}=&2\sum_{(i,j)\in \mathcal{E}_{i}^{out}}
\left\lbrace(G_{1,ij}^k)^{\top}\Sigma_{1}G_{1,ij}^k\right\rbrace 
+2\sum_{(l,i)\in \mathcal{E}_{i}^{in}}
\left\lbrace(G_{2,ij}^k)^{\top}\Sigma_{2}G_{2,ij}^k\right\rbrace
+\beta_{1}I_{4}+H_{1,i},
    \nonumber\\
b_{1,i}^{k}=& b_{1,i}^{k} = 
-2\sum_{(i,j)\in \mathcal{E}_{i}^{out}}
(G_{1,ij}^k)^{\top}\Sigma_{1}(\tilde{t}_j^{k} - \tilde{s}_i^{k})
-2\sum_{(l,i)\in \mathcal{E}_{i}^{in}}
(G_{2,ij}^k)^{\top}\Sigma_{2}\tilde{1}
-\beta_1\tilde{q}_{i}^{k}-\bm{\lambda}_{i}^{k}-H_{1,i}\tilde{p}_{i}^{k+\frac{1}{2}}.
\end{align*}

\subsection{$\tilde{\bm{q}}$--subproblems}
For the $\tilde{\bm{q}}$-subproblem, \eqref{subproblem_q_1} can be written as
\begin{align}
	\tilde{\bm{q}}^{k+1} &= 
	\mathop{\arg\min}\limits_{\tilde{\bm{q}} \in \mathbb{R}^{4n}}
	\sum_{(i,j)\in \mathcal{E}} 
	\|W(\tilde{p}_{i}^{k+1})^{\top}W(\tilde{t}_{ij})\tilde{q}_{i}-(\tilde{t}_j^{k} - \tilde{s}_i^{k}) \|_{\Sigma_{1}}^2
	+\|W(\tilde{q}_{ij})M(\tilde{p}_{j}^{k+1})^{\top}\tilde{q}_{i}-1\|_{\Sigma_{2}}^2.
	\nonumber\\	
	&\quad +\sum_{i=1}^{n} \left\lbrace 
	\frac{\beta_1}{2}\|\tilde{q}_{i}-(\tilde{p}_{i}^{k+1}-\frac{1}{\beta_1}\bm{\lambda}_{i}^{k})\|^{2}
	+\frac{1}{2}\|\tilde{q}_{i}-\tilde{q}_{i}^{k+\frac{1}{2}}\|_{H_{2,i}}^2 \right\rbrace.
	\nonumber
\end{align}
Then we can update $\tilde{q}_{i}$ in parallel. For $i=1,2,\dots,n$, we have
\begin{align*}
	\tilde{q}_{i}^{k+1}  =
	\mathop{\arg\min}\limits_{\tilde{q}_{i}} ~ \frac{1}{2}\tilde{q}_{i}^{\top}A_{2,i}^{k}\tilde{q}_{i} + (b_{2,i}^{k})^{\top}\tilde{q}_{i} = (A_{2,i}^{k})^{-1}b_{2,i}^{k},
\end{align*}
where $G_{3,ij}^{k}=W(\tilde{p}_{i}^{k+1})^{\top}W(\tilde{t}_{ij})$, $G_{4,ij}^{k}=W(\tilde{q}_{ij})M(\tilde{p}_{j}^{k+1})^{\top}$, and
\begin{align*}
A_{2,i}^{k}=&\sum_{(i,j)\in \mathcal{E}_{i}^{out}}
	2\left\lbrace
	(G_{3,ij}^k)^{\top}\Sigma_{1}G_{3,ij}^k
	+(G_{4,ij}^k)^{\top}\Sigma_{2}G_{4,ij}^k
	\right\rbrace +\beta_1 I_{4}+H_{2,i},
    \nonumber\\
b_{2,i}^{k}=& \sum_{(i,j)\in \mathcal{E}_{i}^{out}} 
	2\left\lbrace 
	(G_{3,ij}^k)^{\top}\Sigma_{1}(\tilde{t}_j^{k} - \tilde{s}_i^{k})
	+(G_{4,ij}^k)^{\top}\Sigma_{2}1		
	\right\rbrace +\beta_1\tilde{p}_{i}^{k+1}-\bm{\lambda}_{i}^{k}+H_{2,i}\tilde{q}_{i}^{k+\frac{1}{2}}.
\end{align*}

\subsection{$\bm{t}$--subproblems}
For the $\bm{t}$-subproblem, denote 
$$
\Sigma_{1}=\left( \begin{array}{*{20}{c}}
	\sigma_{11}& \bm{\sigma}_{12}^{\top}\\
	\bm{\sigma}_{21}& \hat{\Sigma}_{1}
\end{array}
\right). 
$$
We can reformulate \eqref{subproblem_t_1} as 
%\left\right
\begin{align*}
	\bm{t}_{i}^{k+1} =& \mathop{\arg\min}\limits_{\bm{t}_{i}}
	\sum_{(l,i)\in \mathcal{E}_{i}^{in}}
	\|\bm{t}_{i}-\bm{s}_{l}^{k}-\bm{c}_{li}^{k}\|_{\hat{\Sigma}_{1}}^{2}
	+\frac{\beta_2}{2}\|\bm{t}_{i}-\bm{s}_{i}^{k}-\frac{1}{\beta_2}\bm{z}_{i}^{k}\|^{2}
	+\frac{1}{2}\|\bm{t}_{i}-\bm{t}_{i}^{k+\frac{1}{2}}\|_{H_{3,i}}^2=(A_{3,i}^{k})^{-1}b_{3,i}^{k},
\end{align*}
where $\tilde{c}_{li}^{k}=\tilde{q}_{i}^{k+1}\tilde{t}_{ij}(\tilde{p}_{i}^{k+1})^{*}$; $\bm{c}_{li}^{k}$ is the imaginary part of $\tilde{c}_{li}^{k}$, and
\begin{align*}
A_{3,i}^{k}=&\sum_{(l,i)\in \mathcal{E}_{i}^{in}}2\hat{\Sigma}_{1}+\beta_2 I_3 + H_{3,i},
    \nonumber\\
b_{3,i}^{k}=& \sum_{(l,i)\in \mathcal{E}_{i}^{in}}2\hat{\Sigma}_{1}(\bm{s}_{l}^{k}+\bm{c}_{li}^{k})
	+\beta_2\bm{s}_{i}^{k}+\bm{z}_{i}^{k}+H_{3,i}\bm{t}_{i}^{k+\frac{1}{2}}.
\end{align*}

\subsection{$\bm{s}$--subproblems}
For the $\bm{s}$-subproblem, the update scheme is similar as $\bm{t}_{i}^{k+1}$, we have
\begin{align*}
	\bm{t}_{i}^{k+1}=& \mathop{\arg\min}\limits_{\bm{t}_{i}}
	 \sum_{(i,j)\in \mathcal{E}_{i}^{out}}
	\|\bm{s}_{i}-\bm{t}_{j}^{k+1}+\bm{c}_{ij}^{k}\|_{\hat{\Sigma}_{1}}^{2}
	+\frac{\beta_2}{2}\|\bm{s}_{i}-\bm{t}_{i}^{k+1}+\frac{1}{\beta_2}\bm{z}_{i}^{k}\|^{2}
	+\frac{1}{2}\|\bm{s}_{i}-\bm{s}_{i}^{k+\frac{1}{2}}\|_{H_{4,i}}^2
    =(A_{4,i}^{k})^{-1}b_{4,i}^{k},
\end{align*}
where
\begin{align*}
A_{4,i}^{k}=&\sum_{(i,j)\in \mathcal{E}_{i}^{out}}2\hat{\Sigma}_{1}+\beta_2 I_3 + H_{4,i},
    \nonumber\\
    b_{4,i}^{k}=& \sum_{(i,j)\in \mathcal{E}_{i}^{out}}
	2\hat{\Sigma}_{1}(\bm{t}_{j}^{k+1}-\bm{c}_{ij}^{k})
	+\beta_2\bm{t}_{i}^{k+1}-\bm{z}_{i}^{k}+H_{4,i}\bm{s}_{i}^{k+\frac{1}{2}}.
\end{align*}

\section{Discussions and Proofs for the Generalized Eigenvalue Problem in Section \ref{sec-algorithm}}  \label{app-inertial Riemannian ADMM}
\subsection{A Lemma for Proving Theorem \ref{TRS-subproblem3}}
Denote
$$
\tilde{Q}(\lambda)=\left(\begin{array}{cc}
	-I & A+\lambda I\\
	A+\lambda I & \frac{-gg^{\top}}{\Delta^{2}}.
\end{array}\right),
\quad
\bar{Q}(\lambda)=\left(\begin{array}{ccc}
	\Delta^{2} &0& g^{\top}\\
	0& -I& A+\lambda I\\
	g& A+\lambda I & O_{n}
\end{array}\right).
$$
Then $\tilde{Q}(\lambda)$ is the Schur complement of $\bar{Q}(\lambda)$, which implies $\det(\tilde{Q}(\lambda))=\Delta^{2}\det(\bar{Q}(\lambda))$.
\begin{lemma}{\label{TRS-subproblem2}}
	Consider the problem \eqref{eq-lemma4.2-1}. Then we have
	\begin{enumerate}
		\item[(\rmnum{1})] The vector $x^{*}$ is a global solution of problem \eqref{eq-lemma4.2-1} if and only if there is a real number $\lambda^{*}$ such that the following conditions are satisfied:
		\begin{align}\label{eq-lemma4.2-2}
			\|x^{*}\|= \Delta, \quad 
			(A+\lambda^{*} I)x^{*} = -g, \quad
			A+\lambda^{*} I \succeq 0.
		\end{align}
		
		\item[(\rmnum{2})] Suppose $(x^{*},\lambda^{*})$ is a global solution, which satisfies \eqref{eq-lemma4.2-2}. Let $\bar{W}=[w_{1},\dots,w_{n}]$ be the matrix that satisfies $\bar{W}^{\top}A\bar{W}=\Lambda:=\operatorname{diag}(-\lambda_{1},\dots,-\lambda_{n})$, and $\lambda_{1}\leq\dots\leq\lambda_{n-t}<\lambda_{n-t+1}=\dots=\lambda_{n}$.
		Then we have either $\lambda^{*}=\lambda_{n}$, or $\lambda^{*}>\lambda_{n}$ and $h(\lambda^{*})=0$, where $h(\lambda)=\sum_{i=1}^{n} \frac{(w_{i}^{\top}g)^{2}}{(\lambda-\lambda_{i})^{2}}-\Delta^{2}$. 
		
		\item[(\rmnum{3})] Suppose $(x^{*},\lambda^{*})$ is a global solution that satisfies \eqref{eq-lemma4.2-2}. Then we have $\operatorname{det}(\bar{Q}(\lambda^{*}))=\operatorname{det}(\tilde{Q}(\lambda^{*}))=0$.
		Furthermore, there always exists $\lambda^{*}\in\left. \left[ \lambda_{n},+\infty\right.\right) $ which is a solution.
	\end{enumerate}
\end{lemma}
\begin{proof}
(\rmnum{1}) 
(\rmnum{1}) Define a quadratic function 
\begin{equation*}
	\hat{m}(x)=g^{\top}x+\frac{1}{2}x^{\top}(A+\lambda I)x
	=m(x)+\frac{\lambda}{2}x^{\top}x.
\end{equation*}
When there is $\lambda^{*}$ such that the conditions \eqref{eq-lemma4.2-2} are satisfied, we have $x^{*}$ is a global minimum of $\hat{m}(x)$ and
\begin{align*}
	m(x) \geq m(x^{*})+\frac{\lambda}{2}\left((x^{*})^{\top}x^{*}-x^{\top}x\right)
	=m(x^{*})+\frac{\lambda}{2}\left(\Delta^{2}-x^{\top}x\right).
\end{align*}
For any $x$ with $\|x\|=\Delta$, we have $m(x) \geq m(x^{*})$ and $x^{*}$ is a global solution of problem \eqref{eq-lemma4.2-1}. 

For the converse, if $x^{*}$ is a global solution of problem \eqref{eq-lemma4.2-1}, there is a scalar $\lambda^{*}$ such that
\begin{align*}
	\|x^{*}\|= \Delta, \quad 
	(A+\lambda^{*} I)x^{*} = -g.
\end{align*}
from the KKT conditions. Since $m(x) \geq m(x^{*})$ for any $x$ with $\|x\|=\Delta$, we have 
\begin{align}
	m(x) \geq m(x^{*})+\frac{\lambda}{2}\left((x^{*})^{\top}x^{*}-x^{\top}x\right),
	\label{eq-lem4.2-L1}
\end{align}
Combining $(A+\lambda^{*} I)x^{*} = -g$, we can reformulate \eqref{eq-lem4.2-L1} as 
\begin{align*}
	\frac{1}{2}(x-x^{*})^{\top}(A+\lambda^{*} I)(x-x^{*}) \geq 0.
\end{align*}
Since all directions $d=x-x^{*}$ consist of a half space and can be expanded to the whole space by changing the sign, we obtain $A+\lambda^{*} I \succeq 0$.
\noindent(\rmnum{2}) Since $A+\lambda^{*}I \succeq 0$, it is not hard to conclude $\lambda^{*}\geq \lambda_{n}$. If $\lambda^{*}> \lambda_{n}$, the matrix $A+\lambda^{*}I$ is nonsingular and we can represent $x^{*}=-(A+\lambda^{*}I)^{-1}g$. Together with  $\|x^{*}\|=\Delta$, we have 
	\begin{align*}
		g^{\top}\bar{W}(\Lambda+\lambda^{*}I)^{-2}\bar{W}^{\top}g=\Delta^{2},
	\end{align*}
	which implies $h(\lambda^{*})=0$.
	
	\noindent(\rmnum{3}) If $A+\lambda I$ is nonsingular, i.e., $\lambda^{*}>\lambda_{n}$, we denote $x(\lambda)=-(A+\lambda I)^{-1}g$ and $X(\lambda)=\left(\begin{array}{ccc}
		1 &0& 0\\
		x(\lambda)& I& O_{n}\\
		0& O_{n} & I
	\end{array}\right)$.
	Then $\det X(\lambda)=1$, and 
	\begin{align}
		\det\bar{Q}(\lambda)
		&=\det X(\lambda)^{\top}\bar{Q}(\lambda)X(\lambda)\nonumber\\
		&=\det \left(\begin{array}{ccc}
			\Delta^{2}-x(\lambda)^{\top}x(\lambda) &-x(\lambda)^{\top}& 0\\
			-x(\lambda)& -I & A+\lambda I\\
			0& A+\lambda I & O_{n}
		\end{array}\right)\nonumber\\
		&=(-1)^{n}\det(A+\lambda I)^{2}(\Delta^{2}-x(\lambda)^{\top}x(\lambda))
		\label{eq-lemma4.2-3}\\
		&=(-1)^{n+1}\det(A+\lambda I)^{2}h(\lambda).
		\label{eq-lemma4.2-4}
	\end{align}  
	
	Suppose $(x^{*},\lambda^{*})$ is a global solution, and $\lambda^{*}>\lambda_{n}$. Then we have $x^{*}=x(\lambda^{*})$ such that $\|x^{*}\|=\Delta$, and $\det\bar{Q}(\lambda^{*})=0$ by \eqref{eq-lemma4.2-3}.
	
	When $\lambda^{*}=\lambda_{n}$, we have $g\perp \mathcal{N}(A+\lambda^{*} I)$, which means $w_{i}^{\top}g=0$, $\forall~ n-t+1\leq i \leq n$. Define
	\begin{align}\label{eq-lemma4.2-5}
		x(\lambda^{*},\bm{\alpha})=-(A+\lambda^{*} I)^{\dagger}g + w(\bm{\alpha})
	\end{align}
	where $w(\bm{\alpha})=\sum_{i=1}^{t}\alpha_{i}w_{n-t+i} \in \operatorname{span}\{w_{n-t+1},\dots,w_{n}\}$. We have  $(A+\lambda^{*} I)x(\lambda^{*},\bm{\alpha})=-g$ for any $\bm{\alpha}$. There exists an $\bm{\alpha}^{*}$ such that $x(\lambda^{*},\bm{\alpha}^{*})=x^{*}$, and $\|x(\lambda^{*},\bm{\alpha}^{*})\|=\Delta$. Replacing $x(\lambda)$ with $x(\lambda^{*},\bm{\alpha}^{*})$ in \eqref{eq-lemma4.2-3}, we have
	\begin{align*}
		\det\bar{Q}(\lambda^{*})
		%&=\det \left(\begin{array}{ccc} \Delta^{2}-x(\lambda^{*},\bm{\alpha}^{*})^{\top}x(\lambda^{*},\bm{\alpha}^{*})
            &=\det \left(\begin{array}{ccc}
			\Delta^{2}-x(\lambda^{*},\bm{\alpha}^{*})^{\top}x(\lambda^{*},\bm{\alpha}^{*})
            &-x(\lambda^{*},\bm{\alpha}^{*})^{\top}
			& 0\\
			-x(\lambda^{*},\bm{\alpha}^{*})& -I & A+\lambda^{*} I\\
			0& A+\lambda^{*} I & O_{n}
		\end{array}\right)\nonumber\\
		&=(-1)^{n}\det(A+\lambda^{*} I)^{2}(\Delta^{2}-\sum_{i=1}^{n-t}
		\frac{(w_{i}^{\top}g)^{2}}{(\lambda^{*}-\lambda_{i})^{2}}-\|\bm{\alpha}^{*}\|^{2}).
		\label{eq-lemma4.2-6}
	\end{align*}
    It follows from $\det(A+\lambda^{*} I)=0$ that $\det\bar{Q}(\lambda^{*})=0$. Hence, $\det(\tilde{Q}(\lambda^{*}))=\Delta^{2}\det(\bar{Q}(\lambda^{*}))=0$. 
	
	Note that the objective function in \eqref{eq-lemma4.2-1} is continuous, and the constraint set is bounded. Then the solution set is nonempty, and there always exists $\lambda\in\left. \left[ \lambda_{n},+\infty\right.\right) $ such that
	$\operatorname{det}(\bar{Q}(\lambda))=\operatorname{det}(\tilde{Q}(\lambda))=0$. 
\end{proof}
\subsection{Proof of Theorem \ref{TRS-subproblem3}}
\begin{proof}
	Since $\operatorname{det}(\tilde{Q}(\lambda))=\Delta^{2}\operatorname{det}(\bar{Q}(\lambda))$, we only discuss $\operatorname{det}(\tilde{Q}(\lambda))$ in the proof.
	
	Firstly, suppose that $(x^{*},\lambda^{*})$ is a global solution of \eqref{eq-lemma4.2-1}. We have $\lambda^{*}\geq \lambda_{n}$ by \eqref{eq-lemma4.2-2} and $\operatorname{det}(\tilde{Q}(\lambda^{*}))=0$ by (\rmnum{3}). When $\lambda^{*}>\lambda_{n}$, since $h(\lambda)$ is strictly decreasing on $(\lambda_{n},+\infty)$, we have $h(\lambda)$ has exactly one zero point on $(\lambda_{n},+\infty)$. Hence, by \eqref{eq-lemma4.2-4}, $\tilde{Q}(\lambda)$ has exactly one real eigenvalue larger than $\lambda_{n}$, which must be $\lambda^{*}$.
	
	When $\lambda^{*}=\lambda_{n}$, we have $w_{i}^{\top}g=0$, $\forall~ n-t+1\leq i \leq n$, and $h(\lambda)=\sum_{i=1}^{n-t}
	\frac{(w_{i}^{\top}g)^{2}}{(\lambda-\lambda_{i})^{2}}-\Delta^{2}$. If $h(\lambda^{*})\leq 0$, there is no solution for $h(\lambda)= 0$ from monotonicity. Hence, $\lambda^{*}=\lambda_{n}$ is the largest real eigenvalue of $\tilde{Q}(\lambda)$.
    If $h(\lambda^{*})> 0$, there is no $\bm{\alpha}^{*}$ such that \eqref{eq-lemma4.2-5} and 
	\begin{align*}
		\|x(\lambda^{*},\bm{\alpha}^{*})\|^{2}=
		\sum_{i=1}^{n-t}
		\frac{(w_{i}^{\top}g)^{2}}{(\lambda^{*}-\lambda_{i})^{2}}+\|\bm{\alpha}^{*}\|^{2}
		-\Delta^{2}
		=0,
	\end{align*}
	which is a contradiction with the existence of $x^{*}$.
	
	For the converse, suppose that $\lambda^{*}$ is the largest real eigenvalue of $\tilde{Q}(\lambda)$. By Lemma \ref{TRS-subproblem2}(\rmnum{3}), we know $\lambda^{*}\in \left. \left[ \lambda_{n},+\infty\right.\right)$. If $\lambda^{*}>\lambda_{n}$, then we have $h(\lambda^{*})=0$ and $x^{*}=-(A+\lambda^{*} I)^{-1}g$ with $\|x^{*}\|=\Delta$, which implies $(x^{*},\lambda^{*})$ is a global solution.
	
	Then we discuss the case that $\lambda^{*}=\lambda_{n}$. If $\exists~ i\geq n-t+1$, such that $w_{i}^{\top}g\neq0$, there is no solution $x^{*}$ satisfying $(A+\lambda^{*}I)x^{*}=-g$. Hence, $\lambda^{*}$ is not a solution. There is no solution $\lambda \in\left. \left[ \lambda_{n},+\infty\right.\right)$ such that $\det \tilde{Q}(\lambda)=0$, which is a contradiction. Otherwise, $\forall~ i\geq n-t+1$, $w_{i}^{\top}g=0$. When $h(\lambda^{*})=\sum_{i=1}^{n-t}
	\frac{(w_{i}^{\top}g)^{2}}{(\lambda^{*}-\lambda_{i})^{2}}-\Delta^{2}\leq0$, we can find $\bm{\alpha}^{*}$ such that $(A+\lambda^{*} I)x(\lambda^{*},\bm{\alpha}^{*})=-g$, and $\|x(\lambda^{*},\bm{\alpha}^{*})\|=\Delta$, where $x(\lambda^{*},\bm{\alpha})$ is defined in \eqref{eq-lemma4.2-5}. When $h(\lambda^{*})>0$, there is also no solution $x^{*}$ satisfying $(A+\lambda^{*}I)x^{*}=-g$, which is a contradiction. 
%	\cx{Can I say $\lambda^{*}=\lambda_{n}$ must be a solution by Lemma \ref{TRS-subproblem2}(\rmnum{3}) without a proof (Instead of the whole paragraph)?}
\end{proof}

\begin{remark}
	In fact, the problem \eqref{eq-lemma4.2-1} with spherical constraint is related to the boundary solution of a trust region problem \cite{nocedal1999numerical,adachi2017solving} 
	\begin{equation}
		\min_{x\in\mathbb{R}^{n}} m(x)=g^{\top}x+\frac{1}{2}x^{\top}Ax,  \qquad 
		\operatorname{s.t.}~\|x\|\leq \Delta, 
		\label{eq-lemma4.1-1}
	\end{equation}
	where $A \in \mathbb{S}^{n \times n}$. The optimality conditions is that there is a real number $\lambda^{*}\geq0$ such that the following conditions are satisfied:
				\begin{align}\label{eq-lemma4.1-2}
					\|x^{*}\|\leq \Delta, \quad 
					(A+\lambda^{*} I)x^{*} = -g, \quad
					\lambda^{*} (\Delta-\|x^{*}\|)=0, \quad
					A+\lambda^{*} I \succeq 0.
				\end{align}
	When $A\succ0$ and $\|A^{-1}g\|<\Delta$, we find a global solution of \eqref{eq-lemma4.1-1} directly with $\lambda^{*}=0$. But the solution of \eqref{eq-lemma4.2-1} always stays on the boundary. 
	
	In addition, when $A\succ0$, $A+\lambda^{*} I \succeq 0$ always holds for the trust region problem, and any KKT point is a global solution (Problem \eqref{eq-lemma4.1-1} is a convex optimization problem). However, problem \eqref{eq-lemma4.2-1} is non-convex. Since $\lambda^{*}\in \mathbb{R}$, we must select the global solutions from all KKT points by the condition $A+\lambda^{*} I \succeq 0$. It is worth noting that the unique difference between \eqref{eq-lemma4.2-2} and \eqref{eq-lemma4.1-2} with $\|x^{*}\|=\Delta$ is $\lambda^{*}\in\mathbb{R}$ or $\lambda^{*}\geq 0$. 
\end{remark}

\section{Proofs for Section \ref{sec-Convergence Analysis}}  \label{app-Convergence Analysis}

First, we give the first-order optimality condition, which follows from Section \ref{sec-Riemannian optimization} directly.
\begin{lemma} (Optimality Conditions) 
	If there exists a Lagrange multiplier $\lambda^{*}$ such that $x_{i}^{*} \in \mathcal{X}_{i}$ and 
	\begin{equation*}
		\left\lbrace 
		\begin{aligned}
			&0=\operatorname{Proj}_{T_{x_{1}^{*}}\mathcal{M}}
			\left( \nabla_{1}f(x_{1}^{*},x_{2}^{*},x_{3}^{*},x_{4}^{*})
			+\nabla_{1}g(x_{1}^{*},x_{2}^{*})-A_1^{\top}\lambda^{*}
			\right),\\
%			+N_{\mathcal{M}}(x_{1}^{*})
			&0=\nabla_{2}f(x_{1}^{*},x_{2}^{*},x_{3}^{*},x_{4}^{*})
			+\nabla_{2}g(x_{1}^{*},x_{2}^{*})-A_2^{\top}\lambda^{*},\\
			&0=\nabla_{3}f(x_{1}^{*},x_{2}^{*},x_{3}^{*},x_{4}^{*})-A_3^{\top}\lambda^{*},\\
			&0=\nabla_{4}f(x_{1}^{*},x_{2}^{*},x_{3}^{*},x_{4}^{*})-A_4^{\top}\lambda^{*},\\
			&0=A_{1}x_{1}^{*} + A_{2}x_{2}^{*} +  A_{3}x_{3}^{*} + A_{4}x_{4}^{*},
		\end{aligned}
		\right.
        		\label{eq-optimality of PGO}
	\end{equation*}
	then $(x_{1}^{*},x_{2}^{*},x_{3}^{*},x_{4}^{*})$ is a stationary point of the problem \eqref{sec5-problem}.
	\label{thm-optimality of PGO}
\end{lemma}
Before presenting the main results, we show the first-order optimality conditions of each
subproblem, which is fundamental to the following analysis:
\begin{small}
\begin{subequations}
	\begin{numcases}{}
		0 =
		\operatorname{Proj}_{T_{x_{\!1}^{\!k\!+\!1\!}}\mathcal{M}}
		\left( \nabla_{1}f(x_{1}^{k+1},x_{2}^{k},x_{3}^{k},x_{4}^{k})
		+\nabla_{1}g(x_{1}^{k+1},x_{2}^{k})-A_{1}^{\top}\lambda^{k}
		+\beta (x_{1}^{k+1}-x_{2}^{k})
		+H_1(x_{1}^{k+1}-x_{1}^{k})\right),
		\label{eq-subsubgradient-1}
		\\
		0 = \nabla_{2}f(x_{1}^{k+1},x_{2}^{k+1},x_{3}^{k},x_{4}^{k})
		+\nabla_{2}g(x_{1}^{k+1},x_{2}^{k+1}) -A_{2}^{\top}\lambda^{k}
		-\beta (x_{1}^{k+1}-x_{2}^{k+1})
		+H_2(x_{2}^{k+1}-x_{2}^{k}),\label{eq-subsubgradient-2}
		\\
		0 = \nabla_{3}f(x_{1}^{k+1},x_{2}^{k+1},x_{3}^{k+1},x_{4}^{k})
		-A_{3}^{\top}\lambda^{k} 
		+\beta (x_{3}^{k+1}-x_{4}^{k})
		+H_3(x_{3}^{k+1}-x_{3}^{k}),\label{eq-subsubgradient-3}
		\\
		0 = \nabla_{4}f(x_{1}^{k+1},x_{2}^{k+1},x_{3}^{k+1},x_{4}^{k+1})
		-A_{4}^{\top}\lambda^{k}
		-\beta (x_{3}^{k+1}-x_{4}^{k+1})
		+H_4(x_{4}^{k+1}-x_{4}^{k}),\label{eq-subsubgradient-4}
		\\
		\lambda^{k+1} = \lambda^{k} - \tau\beta(A_{1}x_{1}^{k+1} + A_{2}x_{2}^{k+1} +  A_{3}x_{3}^{k+1} + A_{4}x_{4}^{k+1}).
		\label{eq-subsubgradient-5}
	\end{numcases} 
\end{subequations}
\end{small}

\subsection{Proof of Lemma \ref{lem-5.1}}
\begin{proof}
    From \eqref{eq-subsubgradient-2}, \eqref{eq-subsubgradient-4}, \eqref{eq-subsubgradient-5}, and the definition of $A_i$, we have
	\begin{small}
	\begin{align*}
		&A_{2}^{\top}\lambda^{k}-\frac{1}{\tau}A_{2}^{\top}(\lambda^{k}-\lambda^{k+1})
		=z_{1}^{k+1}
		:=\nabla_{2}f(x_{1}^{k+1},x_{2}^{k+1},x_{3}^{k},x_{4}^{k})
		+\nabla_{2}g(x_{1}^{k+1},x_{2}^{k+1})
		+H_2(x_{2}^{k+1}-x_{2}^{k})
		,\\
		&A_{4}^{\top}\lambda^{k}-\frac{1}{\tau}A_{4}^{\top}(\lambda^{k}-\lambda^{k+1})
		=z_{2}^{k+1}
		:=\nabla_{4}f(x_{1}^{k+1},x_{2}^{k+1},x_{3}^{k+1},x_{4}^{k+1})
		+H_4(x_{4}^{k+1}-x_{4}^{k}).
	\end{align*}
	\end{small}
	Hence,
	\begin{align}\label{eq-lemma5.1-0}
		\left( \begin{array}{c}
			A_{2}^{\top}\lambda^{k+1}\\
			A_{4}^{\top}\lambda^{k+1}
		\end{array}	\right)
		=(1-\tau)\left( \begin{array}{c}
			A_{2}^{\top}\lambda^{k}\\
			A_{4}^{\top}\lambda^{k}
		\end{array}	\right) 
		+\tau\left( \begin{array}{c}
			z_{1}^{k+1}\\
			z_{2}^{k+1}
		\end{array}	\right),
	\end{align}
	which implies that 
	\begin{align}\label{eq-lemma5.1-1}
		\Delta \lambda^{k+1} = (1-\tau)\Delta \lambda^{k}
		-\tau\Delta z^{k+1},
	\end{align}
	where $z=(z_{1},z_{2})$, and $\Delta z^{k+1}=z^{k+1}-z^{k}$. We now consider two cases.
	    
         Case $\Rmnum{1}$: $0<\tau\leq1$. The convexity of $\|\cdot\|^2$ gives 
		\begin{align}\label{eq-lemma5.1-2}
			\|\Delta \lambda^{k+1}\|^2 \leq  (1-\tau)\|\Delta \lambda^{k}\|^2
			+\tau\|\Delta z^{k+1}\|^2.
		\end{align}
		
        Case $\Rmnum{2}$: $1<\tau<2$. We rewrite \eqref{eq-lemma5.1-1} as 
		\begin{align*}
			\Delta \lambda^{k+1} = -(\tau-1)\Delta \lambda^{k}
			-\frac{\tau}{2-\tau}(2-\tau)\Delta z^{k+1},
		\end{align*}
		which implies that 
		\begin{align}\label{eq-lemma5.1-3}
			\|\Delta \lambda^{k+1}\|^{2} \leq (\tau-1)\|\Delta \lambda^{k}\|^{2}
			+\frac{\tau^{2}}{2-\tau}\|\Delta z^{k+1}\|^{2}.
		\end{align}
	Combine \eqref{eq-lemma5.1-2} and \eqref{eq-lemma5.1-3}, we have
	\begin{align*}
		\|\Delta \lambda^{k+1}\|^{2} \leq |1-\tau|\|\Delta \lambda^{k}\|^{2}
		+\frac{\tau^{2}}{1-|1-\tau|}\|\Delta z^{k+1}\|^{2},
	\end{align*}
	or equivalently,
	\begin{align}\label{eq-lemma5.1-4}
		(1-|1-\tau|)\|\Delta \lambda^{k+1}\|^{2} \leq 
		|1-\tau|\left( \|\Delta \lambda^{k}\|^{2} -\|\Delta \lambda^{k+1}\|^{2}\right) 
		+\frac{\tau^{2}}{1-|1-\tau|}\|\Delta z^{k+1}\|^{2}.
	\end{align}
	By using the Lipschitz continuity of $\nabla f$ and $\nabla g$ defined in Lemma \ref{assumption1} (ii), we have
	\begin{align}\label{eq-lemma5.1-5}
		\|\Delta z_{1}^{k+1}\|^{2}=&
		\|\nabla_{2}f(x_{1}^{k+1},x_{2}^{k+1},x_{3}^{k},x_{4}^{k})
		-\nabla_{2}f(x_{1}^{k},x_{2}^{k},x_{3}^{k-1},x_{4}^{k-1})
		\nonumber\\
		&+\nabla_{2}g(x_{1}^{k+1},x_{2}^{k+1})-\nabla_{2}g(x_{1}^{k},x_{2}^{k})
		+H_2\Delta x_{2}^{k+1}
		-H_2\Delta x_{2}^{k}\|^{2}\nonumber\\
		\leq&4(L_{f}^{2}+L_{g}^{2})\|\Delta x_{1}^{k+1}\|^2
		+4\|\Delta x_{2}^{k+1}\|_{(L_{f}^{2}+L_{g}^{2})I+H_{2}^{\top}H_{2}}^2
        +4L_{f}^{2}\|\Delta x_{3}^{k}\|^2
		+4L_{f}^{2}\|\Delta x_{4}^{k}\|^2
		+4\|\Delta x_{2}^{k}\|_{H_{2}^{\top}H_{2}}^{2}.
	\end{align}
	Similarly, we can bound $\Delta z_{2}^{k+1}$ by
	\begin{align}\label{eq-lemma5.1-7}
		\|\Delta z_{2}^{k+1}\|^{2}
		\leq&3L_{f}^{2}(\|\Delta x_{1}^{k+1}\|^2+\|\Delta x_{2}^{k+1}\|^2+\|\Delta x_{3}^{k+1}\|^2)
		+3\|\Delta x_{4}^{k+1}\|_{L_{f}^{2}I+H_{4}^{\top}H_{4}}^2
		+3\|\Delta x_{4}^{k}\|_{H_{4}^{\top}H_{4}}^2.
	\end{align}
	Equation \eqref{eq-bounded lambda} is obtained directly by \eqref{eq-lemma5.1-4}--\eqref{eq-lemma5.1-7}.
\end{proof}

\subsection{Lemma \ref{lem-5.3} and Lemma \ref{lem-5.4}}
Then we summarize the recursion of $\mathcal{L}_{\beta}$ in the next lemma. To simplify the notations in our analysis, we denote
\begin{align}
	M_1&=\frac{1}{2}H_1-\frac{\alpha_{2}}{\beta}(7L_{f}^{2}+4L_{g}^{2})I, 
	\quad N_{1}=0,
	\label{eq-lem5.3-1}\\
	M_2&= H_2-\frac{\alpha_{2}}{\beta}(7L_{f}^{2}I+4L_{g}^{2}I+4H_{2}^{\top}H_{2}), 
	\quad N_{2}=\frac{4\alpha_{2}}{\beta}H_{2}^{\top}H_{2},
	\label{eq-lem5.3-2}\\
	M_3&= H_3-\frac{3\alpha_{2}}{\beta}L_{f}^{2}I,
	\quad N_{3}=\frac{4\alpha_{2}}{\beta}L_{f}^{2}I,
	\label{eq-lem5.3-3}\\
	M_4&= H_4-\frac{\alpha_{2}}{\beta}(3L_{f}^{2}I+3H_{4}^{\top}H_{4}),
	\quad N_{4}=\frac{\alpha_{2}}{\beta}(4L_{f}^{2}I+3H_{4}^{\top}H_{4}).
	\label{eq-lem5.3-4}
\end{align}
\begin{lemma}\label{lem-5.3}
	Suppose that Assumption \ref{assumption1} holds. Let $\{(\bm{x}^{k},\lambda^{k})\}$ be the sequence generated by \eqref{sec5-ADMM} which is assumed to be bounded, then
	we have
	\begin{align}\label{eq-lem5.3-5}
		&\mathcal{L}_{\beta}(\bm{x}^{k+1},\lambda^{k+1})  
		+\frac{\alpha_{1}}{\tau\beta}\|\Delta\lambda^{k+1}\|^{2}
		+\sum_{i=1}^4 \|\Delta x_{i}^{k+1}\|_{M_{i}}^{2}
		\nonumber\\
		\leq
		&\mathcal{L}_{\beta}(\bm{x}^{k},\lambda^{k})
		+\frac{\alpha_{1}}{\tau\beta}\|\Delta\lambda^{k}\|^{2}
		+\sum_{i=1}^4 \|\Delta x_{i}^{k}\|_{N_{i}}^{2}, 
	\end{align}
	where $M_{i}$ and $N_{i}^{k}$ are defined in \eqref{eq-lem5.3-1}--\eqref{eq-lem5.3-4}.
\end{lemma}
\begin{proof}
    Denote $\bm{x}_{-1}^{k}=(x_{2}^{k},x_{3}^{k},x_{4}^{k})$. From the $x_{1}$-subproblem in \eqref{sec5-ADMM}, we have
	\begin{align*}
		\mathcal{L}_{\beta}(\bm{x}^{k,1},\lambda^{k})+\frac{1}{2}\|x_{1}^{k+1}-x_{1}^{k}\|_{H_{1}}^{2}
		\leq
		\mathcal{L}_{\beta}(x_{1},\bm{x}_{-1}^{k},\lambda^{k})+\frac{1}{2}\|x_{1}-x_{1}^{k}\|_{H_{1}}^{2},~ \forall~x_{1}\in \mathcal{M}.
	\end{align*}
	Let $x_{1}=x_{1}^{k}$, and we have
	\begin{align}
		\mathcal{L}_{\beta}(\bm{x}^{k,0},\lambda^{k})-\mathcal{L}_{\beta}(\bm{x}^{k,1},\lambda^{k})
		&\geq\frac{1}{2}\|x_{1}^{k}-x_{1}^{k+1}\|_{H_{1}}^{2}.
		\label{eq-lem5.2-2}
	\end{align}
	Note that $f$ and $g$ are block multi-convex defined in Lemma \ref{assumption1} (v), which implies that
	$\mathcal{L}_{\beta}(\cdot)$ is also a block multi-convex function about $x_{i}$ while all the other blocks are fixed. Then we have
	\begin{align*}
		\mathcal{L}_{\beta}(\bm{x}^{k,i-1},\lambda^{k})
		\geq
		\mathcal{L}_{\beta}(\bm{x}^{k,i},\lambda^{k})
		+\left\langle \nabla_{i} \mathcal{L}_{\beta}(\bm{x}^{k,i},\lambda^{k}),x_{i}^{k}-x_{i}^{k+1} \right\rangle, 
        \quad \forall ~i = 2,3,4. 
	\end{align*}
	According to the optimality conditions in \eqref{eq-subsubgradient-2}--\eqref{eq-subsubgradient-4}, we obtain
	\begin{align*}
		\mathcal{L}_{\beta}(\bm{x}^{k,i-1},\lambda^{k})
		-\mathcal{L}_{\beta}(\bm{x}^{k,i},\lambda^{k})
		+\left\langle H_{i}(x_{i}^{k+1}-x_{i}^{k}),x_{i}^{k}-x_{i}^{k+1} \right\rangle
		\geq 0. 
	\end{align*}
	By rearranging the terms,
	\begin{align}\label{eq-lemma5.2-3}
		\mathcal{L}_{\beta}(\bm{x}^{k,i-1},\lambda^{k})
		-\mathcal{L}_{\beta}(\bm{x}^{k,i},\lambda^{k})
		&\geq \|x_{i}^{k+1}-x_{i}^{k}\|_{H_{i}}^{2}.
	\end{align}
	Summing \eqref{eq-lem5.2-2} and \eqref{eq-lemma5.2-3} from $i=2$ to $4$, we obtain
    \begin{align*}
    &\mathcal{L}_{\beta}(\bm{x}^{k},\lambda^{k})
    -\mathcal{L}_{\beta}(\bm{x}^{k+1},\lambda^{k})  
    \geq \frac{1}{2}\|x_{1}^{k+1}-x_{1}^{k}\|_{H_{1}}^{2}
    +\sum_{i=2}^4 \|x_{i}^{k+1}-x_{i}^{k}\|_{H_{i}}^{2}.
     \end{align*}
     Note that
	\begin{align}\label{eq-lemma5.2-1}
		\mathcal{L}_{\beta}(\bm{x}^{k+1},\lambda^{k+1}) = 
		\mathcal{L}_{\beta}(\bm{x}^{k+1},\lambda^{k})
		+\frac{1}{\tau\beta}\|\lambda^{k+1}-\lambda^{k}\|^{2}. 
	\end{align}
    Combining Lemma \ref{lem-5.1}, Eq. \eqref{eq-lemma5.2-1}, and rearranging the terms, the result is derived directly.
\end{proof}

Let us denote
\begin{align}\label{def-beta_s1}
\beta^{\prime}=\max
		&\left\lbrace 
		2(L_{f}+L_{g,2}),
		\frac{2\alpha_{2}}{\underline{\kappa}_{1}}(7L_{f}^{2}+4L_{g}^{2}),
		\frac{\alpha_{2}}{\underline{\kappa}_{2}}(7L_{f}^{2}+4L_{g}^{2})
        +\frac{(4\alpha_{2}+2)\overline{\kappa}_{2}^{2}}{\underline{\kappa}_{2}}, 
		\frac{3\alpha_{2}+2}{\underline{\kappa}_{3}}L_{f}^{2},
        \frac{(3\alpha_{2}+2)}{\underline{\kappa}_{4}}L_{f}^{2}
        +\frac{(3\alpha_{2}+1)\overline{\kappa}_{4}^{2}}{\underline{\kappa}_{4}}
		\right\rbrace,
\end{align}
\begin{align}\label{def-beta_s2}
\beta^{\prime\prime}=\max\left\lbrace
\frac{2\alpha_{2}}{\underline{\kappa}_{1}}(7L_{f}^{2}+4L_{g}^{2}),
\frac{\alpha_{2}}{\underline{\kappa}_{2}}
(7L_{f}^{2}+4L_{g}^{2}+8\overline{\kappa}_{2}^{2}),
\frac{7\alpha_{2}L_{f}^{2}}{\underline{\kappa}_{3}},
\frac{\alpha_{2}}{\underline{\kappa}_{4}}(7L_{f}^{2}+6\overline{\kappa}_{4}^{2})
\right\rbrace.
\end{align}
Next, we prove that $\{\Psi^{k}\}$ is bounded from below, which is the key to prove the convergence.
\begin{lemma}\label{lem-5.4}
Suppose that Assumption \ref{assumption1} holds. Let $\{(\bm{x}^{k},\lambda^{k})\}$ be the sequence generated by \eqref{sec5-ADMM} which is assumed to be bounded. $\underline{\kappa}_{i}$, $\overline{\kappa}_{i}>0$ are the smallest and largest eigenvalue of $H_{i}$, respectively. If $\beta > \beta^{\prime}$, then $\{\Psi^{k}\}$ is bounded from below, i.e.,
	\begin{align*}
		\Psi^{k+1} \geq f^{*}+g^{*}.
	\end{align*}
\end{lemma}

\begin{proof}
	It follows from \eqref{eq-subsubgradient-5} and \eqref{eq-lemma5.1-0}, and we have $\tau\lambda^{k+1}=-(1-\tau)\Delta\lambda^{k+1}-\tau z^{k+1}$, and
	\begin{align*}
		\left\langle \lambda^{k+1},\sum_{i=1}^{4}A_{i}x_{i}^{k+1} \right\rangle
		&=\frac{1-\tau}{\tau^{2}\beta}\|\Delta\lambda^{k+1}\|^{2}
		-\left\langle z^{k+1}, \sum_{i=1}^{4}A_{i}x_{i}^{k+1}\right\rangle.
	\end{align*}
	Combining the definition of $z^{k+1}$ in Lemma \ref{lem-5.1}, we can reformulate $\mathcal{L}_{\beta}(\bm{x}^{k+1},\lambda^{k+1})$ as
	\begin{align}
		&\mathcal{L}_{\beta}(\bm{x}^{k+1},\lambda^{k+1}) 
		= f(\bm{x}^{k+1})
		+g(x_{1}^{k+1},x_{2}^{k+1}) 
		-\left\langle \lambda^{k+1},\sum_{i=1}^{4}A_{i}x_{i}^{k+1} \right\rangle 
		+\frac{\beta}{2}\left\| \sum_{i=1}^{4}A_{i}x_{i}^{k+1}\right\|^{2}
		\nonumber\\
		&\qquad= f(\bm{x}^{k+1})
		+g(x_{1}^{k+1},x_{2}^{k+1})+\frac{\beta}{2}\left\| \sum_{i=1}^{4}A_{i}x_{i}^{k+1}\right\|^{2}
		-\frac{1-\tau}{\tau^{2}\beta}\|\Delta\lambda^{k+1}\|^{2}
		\nonumber\\
		&\qquad+\left\langle \nabla_{2}f(x_{1}^{k+1},x_{2}^{k+1},x_{3}^{k},x_{4}^{k})
		+\nabla_{2}g(x_{1}^{k+1},x_{2}^{k+1})
		+H_2(x_{2}^{k+1}-x_{2}^{k}), 
		x_{1}^{k+1}-x_{2}^{k+1}\right\rangle
		\nonumber\\
		&\qquad+\left\langle
		\nabla_{4}f(x_{1}^{k+1},x_{2}^{k+1},x_{3}^{k+1},x_{4}^{k+1})
		+H_4(x_{4}^{k+1}-x_{4}^{k}),
		x_{3}^{k+1}-x_{4}^{k+1}
		\right\rangle
		\label{eq-lem-5.4-1}
	\end{align}
	Since the gradient $\nabla f(\bm{x})$ and the partial gradient $\nabla_{x_{2}}g(x_{1},x_{2})$ are Lipschitz continuous in a bounded subset in Lemma \ref{assumption1}, we obtain
	\begin{align}
		&f(\bm{x}^{k+1})
		+\left\langle \nabla_{2}f(\bm{x}^{k+1}),
		x_{1}^{k+1}-x_{2}^{k+1}\right\rangle
		+\left\langle
		\nabla_{4}f(\bm{x}^{k+1}),
		x_{3}^{k+1}-x_{4}^{k+1}
		\right\rangle
		\nonumber\\
		 \geq&
		f(x_{1}^{k+1},x_{1}^{k+1},x_{3}^{k+1},x_{3}^{k+1})
		-\frac{L_{f}}{2}(\|x_{1}^{k+1}-x_{2}^{k+1}\|^{2}
		+\|x_{3}^{k+1}-x_{4}^{k+1}\|^{2}),
		\label{eq-lem-5.4-2}
	\end{align}
	and
	\begin{align}
		&g(x_{1}^{k+1},x_{2}^{k+1}) 
		+\left\langle
		\nabla_{2}g(x_{1}^{k+1},x_{2}^{k+1}), 
		x_{1}^{k+1}-x_{2}^{k+1}\right\rangle
		\geq g(x_{1}^{k+1},x_{1}^{k+1}) 
		-\frac{L_{g,2}}{2}\|x_{1}^{k+1}-x_{2}^{k+1}\|^{2}.
		\label{eq-lem-5.4-3}
	\end{align}
	In addition, by using Cauchy inequality, we also have 
	\begin{align}
		&\left\langle \nabla_{2}f(x_{1}^{k+1},x_{2}^{k+1},x_{3}^{k},x_{4}^{k})
		-\nabla_{2}f(x_{1}^{k+1},x_{2}^{k+1},x_{3}^{k+1},x_{4}^{k+1}), 
		x_{1}^{k+1}-x_{2}^{k+1}\right\rangle
		\nonumber\\
		&\qquad \geq -\frac{2L_{f}^{2}}{\beta}
		(\|\Delta x_{3}^{k+1}\|^{2}+\|\Delta x_{4}^{k+1}\|^{2})
		-\frac{\beta}{8}\|x_{1}^{k+1}-x_{2}^{k+1}\|^{2},
		\label{eq-lem-5.4-4}
	\end{align}
	\begin{align}
		&\left\langle H_2(x_{2}^{k+1}-x_{2}^{k}),
		x_{1}^{k+1}-x_{2}^{k+1}
		\right\rangle
		\geq -\frac{2}{\beta}\|\Delta x_{2}^{k+1}\|_{H_{2}^{\top}H_{2}}^{2}
		-\frac{\beta}{8}\|x_{1}^{k+1}-x_{2}^{k+1}\|^{2},
		\label{eq-lem-5.4-5}
	\end{align}
	and
	\begin{align}
		&\left\langle H_4(x_{4}^{k+1}-x_{4}^{k}),
		x_{3}^{k+1}-x_{4}^{k+1}
		\right\rangle
		\geq -\frac{1}{\beta}
		\|\Delta x_{4}^{k+1}\|_{H_{4}^{\top}H_{4}}^{2}
		-\frac{\beta}{4}
		\|x_{3}^{k+1}-x_{4}^{k+1}\|^{2}.
		\label{eq-lem-5.4-6}
	\end{align}
%	Let 
%	\begin{align}
%		\frac{1}{2\nu_{1}}=\frac{24\alpha_{2}}{\beta}\\
%		\frac{1}{2\nu_{2}}=\frac{(1-\delta_{2})24\alpha_{2}}{\beta}\\
%		\frac{1}{2\nu_{3}}=\frac{18\alpha_{2}}{\beta}\\
%		\frac{1}{2\nu_{4}}=\frac{(1-\delta_{4})18\alpha_{2}}{\beta}\\
%		\frac{L_{f}^{2}}{2\nu_{5}}=\frac{4\alpha_{2}L_{f}^{2}}{\beta}
%	\end{align}
	Combining \eqref{eq-lem-5.4-1} with \eqref{eq-lem-5.4-2}--\eqref{eq-lem-5.4-6}, we have
	\begin{align*}
		\Psi^{k+1}
		\geq&
		f(x_{1}^{k+1},x_{1}^{k+1},x_{3}^{k+1},x_{3}^{k+1})
		+g(x_{1}^{k+1},x_{1}^{k+1})
		+(\frac{\alpha_{1}}{\tau\beta}-\frac{1-\tau}{\tau^{2}\beta})
		\|\Delta\lambda^{k+1}\|^{2}
		\nonumber\\		
		&+\left( \frac{\beta}{4}
		-\frac{L_{f}}{2}-\frac{L_{g,2}}{2}
		\right) 
		\|x_{1}^{k+1}-x_{2}^{k+1}\|^{2}
		+\left( \frac{\beta}{4}
		-\frac{L_{f}}{2}\right) 
		\|x_{3}^{k+1}-x_{4}^{k+1}\|^{2}
		\nonumber\\ 
		&+\|\Delta x_{1}^{k+1}\|_{M_{1}}^{2}
		+\|\Delta x_{2}^{k+1}\|_{\bar{M}_{2}}^{2}
		+\|\Delta x_{3}^{k+1}\|_{\bar{M}_{3}}^{2}
		+\|\Delta x_{4}^{k+1}\|_{\bar{M}_{4}}^{2}
		,
	\end{align*}
	where 
	$\bar{M}_{2}=M_{2}-\frac{2}{\beta}H_{2}^{\top}H_{2}$,~ 
	$\bar{M}_{3}=M_{3}-\frac{2L_{f}^{2}}{\beta}I$,~
	$\bar{M}_{4}=M_{4}-\frac{1}{\beta}H_{4}^{\top}H_{4}-\frac{2L_{f}^{2}}{\beta}I$.~
	When $\tau\in (0,2)$,
	\begin{equation*}
		\alpha_{1}-\frac{1-\tau}{\tau}
		=\left\lbrace 
		\begin{array}{l}
			\displaystyle0, \qquad\qquad\qquad\qquad\operatorname{if}~ \tau\in \left. \left( 0,1\right. \right] ,\\
			\displaystyle\frac{\tau-1}{2-\tau}+\frac{\tau-1}{\tau}>0, ~\operatorname{if}~ \tau\in (1,2).
		\end{array}
		\right. 
	\end{equation*}
	Therefore, if $\beta > \beta^{\prime}$, we have 
	\begin{align*}
		\Psi^{k+1}\geq &~ f(x_{1}^{k+1},x_{1}^{k+1},x_{3}^{k+1},x_{3}^{k+1})
		+g(x_{1}^{k+1},x_{1}^{k+1})\\
		\geq &~
		f^{*}+g^{*}.
	\end{align*}
\end{proof}
\subsection{Proof of Theorem \ref{thm-decrease}}
\begin{proof}
	(i) Let $\underline{\beta}=\max\{\beta^{\prime},\beta^{\prime\prime}\}$. When $\beta>\underline{\beta}$, there exist $0<\delta_{i}<1$ such that  $\delta_{i}M_{i}\succeq N_{i}$. Then we can rewrite the right-hand side of \eqref{eq-lem5.3-5} as
	\begin{align*}
		&\mathcal{L}_{\beta}(\bm{x}^{k+1},\lambda^{k+1})  
		+\frac{\alpha_{1}}{\tau\beta}\|\Delta\lambda^{k+1}\|^{2}
		+\sum_{i=1}^4 \|\Delta x_{i}^{k+1}\|_{M_{i}}^{2}
		\nonumber\\
		\leq
		&\mathcal{L}_{\beta}(\bm{x}^{k},\lambda^{k})
		+\frac{\alpha_{1}}{\tau\beta}\|\Delta\lambda^{k}\|^{2}
		+\sum_{i=1}^{4} \delta_{i}\|\Delta x_{i}^{k}\|_{M_{i}}^{2}.
	\end{align*}
	Rearranging the terms, we get the first conclusion.
%	\begin{align*}
%		&\mathcal{L}_{\beta}(\bm{x}^{k+1},\lambda^{k+1})  
%		+\frac{\alpha_{1}}{\tau\beta}\|\Delta\lambda^{k+1}\|^{2}
%		+\sum_{i=1}^4 \|\Delta x_{i}^{k+1}\|_{M_{i}}^{2}
%		+\sum_{i=2,4} \hat{\delta}_{i}\|\Delta x_{i}^{k}\|_{M_{i}}^{2}
%		\nonumber\\
%		\leq
%		&\mathcal{L}_{\beta}(\bm{x}^{k},\lambda^{k})
%		+\frac{\alpha_{1}}{\tau\beta}\|\Delta\lambda^{k}\|^{2}
%		+\sum_{i=1}^{4} \delta_{i}\|\Delta x_{i}^{k}\|_{M_{i}}^{2}
%		+\sum_{i=2,4} \hat{\delta}_{i}\|\Delta x_{i}^{k-1}\|_{M_{i}}^{2}.
%	\end{align*}

	(ii) Summing \eqref{eq-thm5.2-L1} from $k=1$ to $K$, we obtain	
	\begin{small}
		\begin{align*}
			\Psi^{K+1}
			+\sum_{k=1}^{K}\sum_{i=1}^4 (1-\delta_{i})\|\Delta x_{i}^{k}\|_{M_{i}}^{2}
			\leq\Psi^{1}. 
%			\label{eq-thm5.2-1}
		\end{align*}
	\end{small}
%	\begin{small}
%		\begin{align}\label{eq-thm5.2-1}
%			&\mathcal{L}_{\beta}(\bm{x}^{K+1},\lambda^{K+1})  
%			+\frac{\alpha_{1}}{\tau\beta}\|\Delta\lambda^{K+1}\|^{2}
%			+\sum_{i=1}^4 \|\Delta x_{i}^{K+1}\|_{M_{i}}^{2}
%			+\sum_{i=1}^{4} (1-\delta_{i})\|\Delta x_{i}^{K}\|_{M_{i}}^{2}
%			+\sum_{k=1}^{K-1}\sum_{i=1}^4 (1-\delta_{i})\|\Delta x_{i}^{k}\|_{M_{i}}^{2}
%			\nonumber\\
%			&\leq
%			\mathcal{L}_{\beta}(\bm{x}^{1},\lambda^{1})
%			+\frac{\alpha_{1}}{\tau\beta}\|\Delta\lambda^{1}\|^{2}
%			+\sum_{i=1}^4 \|\Delta x_{i}^{1}\|_{M_{i}}^{2}
%			+\sum_{i=2,4} \hat{\delta}_{i}\|\Delta x_{i}^{0}\|_{M_{i}}^{2}
%			. 
%		\end{align}
%	\end{small}
	Since $\{\Psi^{k}\}$ is bounded from below in Lemma \ref{lem-5.4}, we derive that $\sum_{k=1}^{+\infty}\sum_{i=1}^4 (1-\delta_{i})\|\Delta x_{i}^{k}\|_{M_{i}}^{2} < +\infty$, which implies $\{\Delta x_{i}^{k}\}$, $i=1,\dots,4$, converge to $0$. It follows from \eqref{eq-lemma5.1-5} and \eqref{eq-lemma5.1-7} that $\sum_{k=1}^{+\infty}\|\Delta z^{k}\|<+\infty$. Then summing \eqref{eq-lemma5.1-4} from $k=1$ to $K$, we obtain
	\begin{align*}
		(1-|1-\tau|)\sum_{k=1}^{K}\|\Delta \lambda^{k+1}\|^{2} \leq 
		|1-\tau| \|\Delta \lambda^{1}\|^{2} 
		+\tau\alpha_{2}\sum_{k=1}^{K}\|\Delta z^{k+1}\|^{2}< +\infty,
	\end{align*}
	which implies that $\{\Delta \lambda^{k}\}$ converges to $0$. 
	
	(iii) Let $(\hat{\bm{x}},\hat{\lambda})$
	be an arbitrary cluster point of $\{(\bm{x}^{k},\lambda^{k})\}$, and $(\bm{x}^{k_j},\lambda^{k_j})$ be the subsequence converging to $(\hat{\bm{x}},\hat{\lambda})$. Next we will show that $(\hat{\bm{x}},\hat{\lambda})$ is a stationary point of $\mathcal{L}_{\beta}$.
	Firstly, by \eqref{eq-subsubgradient-2} and \eqref{eq-subsubgradient-5}, we have that
	\begin{align}
		0=&\nabla_{2}f(x_{1}^{k_j},x_{2}^{k_j},x_{3}^{k_j-1},x_{4}^{k_j-1})
		+\nabla_{2}g(x_{1}^{k_j},x_{2}^{k_j})\nonumber\\
		&-A_{2}^{\top}\lambda^{k_j}
		-(1-\tau)\beta A_{2}^{\top}(A_{1}x_{1}^{k_j}+A_{2}x_{2}^{k_j}+A_{3}x_{3}^{k_j}+A_{4}x_{4}^{k_j})
		+H_2(x_2^{k_j}-x_2^{k_j-1}), ~\forall j.
		\label{eq-thm5.2-2}
	\end{align}
	Since $\{\Delta \lambda^{k}\}$ converges to $0$, it follows from \eqref{eq-subsubgradient-5} that 
	\begin{align}
		\lim_{j\rightarrow +\infty} \sum_{i=1}^{4}A_{i}x_{i}^{k_j}=\sum_{i=1}^{4}A_{i}\hat{x}_{i}=0.
		\label{eq-thm5.2-3}
	\end{align}
	In addition, since $\{\Delta x_{2}^{k}\}$ converges to $0$, we have
	\begin{align*}
		\lim_{j\rightarrow +\infty}H_2(x_2^{k_j}-x_2^{k_j-1})
		=0.
	\end{align*}
	Taking limit along the subsequence in \eqref{eq-thm5.2-2} and using the continuity of $\nabla f$ and $\nabla g$, we have
	\begin{align}
		\nabla_{2}f(\hat{x}_{1},\hat{x}_{2},\hat{x}_{3},\hat{x}_{4})
		+\nabla_{2}g(\hat{x}_{1},\hat{x}_{2})
		-A_{2}^{\top}\hat{\lambda}=0.
		\label{eq-thm5.2-4}
	\end{align}
	Similarly, taking limit in \eqref{eq-subsubgradient-3} and \eqref{eq-subsubgradient-4}, we can obtain
	\begin{align}
		\nabla_{3}f(\hat{x}_{1},\hat{x}_{2},\hat{x}_{3},\hat{x}_{4})	-A_{3}^{\top}\hat{\lambda}=0,\label{eq-thm5.2-5}\\
		\nabla_{4}f(\hat{x}_{1},\hat{x}_{2},\hat{x}_{3},\hat{x}_{4})	-A_{4}^{\top}\hat{\lambda}=0.\label{eq-thm5.2-6}
	\end{align}
	Next, we analyze the optimality condition of the first block. 
	We define 
	$$s^{k+1}= 
	\nabla_{1}f(x_{1}^{k+1},x_{2}^{k},x_{3}^{k},x_{4}^{k})
	+\nabla_{1}g(x_{1}^{k+1},x_{2}^{k})-A_{1}^{\top}\lambda^{k}
	+\beta (x_{1}^{k+1}-x_{2}^{k})
	+H_1(x_{1}^{k+1}-x_{1}^{k}).$$
	From the first-order optimality condition of $x_1$-subproblem \eqref{eq-subsubgradient-1}, we know that the projection of $s^{k+1}$ on $T_{x_{1}^{k+1}}\mathcal{M}$ is $0$. 
	By using the same technique as \eqref{eq-lem-5.4-4}, we have 
	\begin{align}
	\lim_{j\rightarrow +\infty} s^{k_j}=\nabla_{1}f(\hat{x}_{1},\hat{x}_{2},\hat{x}_{3},\hat{x}_{4})
	+\nabla_{1}g(\hat{x}_{1},\hat{x}_{2})
	-A_{1}^{\top}\hat{\lambda}:=\hat{s}.
	\label{eq-thm5.2-7}
	\end{align}
	Then we derive
%\begin{align}
%    &\|\lim_{j\rightarrow +\infty} \operatorname{Proj}_{T_{x_{1}^{k_j}}\mathcal{M}}s^{k_j}-
%    \operatorname{Proj}_{T_{\hat{x}_{1}}\mathcal{M}}\hat{s}\|\nonumber\\
%    \leq & \lim_{j\rightarrow +\infty}
%    \|\operatorname{Proj}_{T_{x_{1}^{k_j}}\mathcal{M}}s^{k_j}
%    -\operatorname{Proj}_{T_{x_{1}^{k_j}}\mathcal{M}}\hat{s}\|
%    +\lim_{j\rightarrow +\infty}
%    \|\operatorname{Proj}_{T_{x_{1}^{k_j}}\mathcal{M}}\hat{s}
%    -\operatorname{Proj}_{T_{\hat{x}_{1}}\mathcal{M}}\hat{s}\|\nonumber\\
%    \leq & \lim_{j\rightarrow +\infty}\|s^{k_j}-\hat{s}\|
%    +\lim_{j\rightarrow +\infty}
%    \|\operatorname{Proj}_{T_{x_{1}^{k_j}}\mathcal{M}}\hat{s}
%    -\operatorname{Proj}_{T_{\hat{x}_{1}}\mathcal{M}}\hat{s}\|,
%    \label{eq-thm5.2-8}
%\end{align}
    \begin{align}
		&\|\lim_{j\rightarrow +\infty} \operatorname{Proj}_{T_{y}\mathcal{M}}s^{k_j}-
		\operatorname{Proj}_{T_{\hat{x}_{1}}\mathcal{M}}\hat{s}\|\nonumber\\
		\leq & \lim_{j\rightarrow +\infty}
		\|\operatorname{Proj}_{T_{y}\mathcal{M}}s^{k_j}
		-\operatorname{Proj}_{T_{y}\mathcal{M}}\hat{s}\|
		+\lim_{j\rightarrow +\infty}
		\|\operatorname{Proj}_{T_{y}\mathcal{M}}\hat{s}
		-\operatorname{Proj}_{T_{\hat{x}_{1}}\mathcal{M}}\hat{s}\|\nonumber\\
		\leq & \lim_{j\rightarrow +\infty}\|s^{k_j}-\hat{s}\|
		+\lim_{j\rightarrow +\infty}
		\|\operatorname{Proj}_{T_{y}\mathcal{M}}\hat{s}
		-\operatorname{Proj}_{T_{\hat{x}_{1}}\mathcal{M}}\hat{s}\|,
		\label{eq-thm5.2-8}
	\end{align}
	where $y = x_{1}^{k_{j}}$, and the second inequality follows from the nonexpansiveness of the projection operator. By \eqref{eq-thm5.2-7}, the first term of right-hand side satisfies $\lim_{j\rightarrow +\infty}\|s^{k_j}-\hat{s}\|=0$. Next, we analyze the second term of \eqref{eq-thm5.2-8}.
	Note that the Riemannian submanifold $\mathcal{M}$ is smooth. The projection on the tangent space $\operatorname{Proj}_{T_{u}\mathcal{M}}(v)$ is also smooth respect to $u$ (Exercise 3.66 in \cite{boumal2023introduction}). Then we have 
	$$
	\lim_{j\rightarrow +\infty}
	\|\operatorname{Proj}_{T_{y}\mathcal{M}}\hat{s}
	-\operatorname{Proj}_{T_{\hat{x}_{1}}\mathcal{M}}\hat{s}\|=0.
	$$
	Finally, we have 
	\begin{align}
		\operatorname{Proj}_{T_{\hat{x}_{1}}\mathcal{M}}\hat{s} = 0.
		\label{eq-thm5.2-9}
	\end{align}
	The subsequential convergence is obtained from \eqref{eq-thm5.2-3}--\eqref{eq-thm5.2-6} and \eqref{eq-thm5.2-9}.
\end{proof}

\subsection{Proof of Lemma \ref{lem-subgradient-bound}}
\begin{proof}
Denote $M_{-1}=\operatorname{diag}(M_{2},M_{3},M_{4})$. By the definition of $\Psi$ (see \eqref{def-Psi}), we have that
	\begin{align*}
		\Psi^{k+1}:
		=\Psi(\bm{x}^{k+1},\bm{x}^{k},\lambda^{k+1},\lambda^{k})
		=\mathcal{L}_{\beta}(\bm{x}^{k+1},\lambda^{k+1})
		+\frac{\alpha_{1}}{\tau\beta}\|\Delta\lambda^{k+1}\|^{2}
		+\|\Delta \bm{x}^{k+1}\|_{M}^{2}.
	\end{align*}
	where $M=\operatorname{diag}(M_{1},M_{2},M_{3},M_{4})$, and
	\begin{align}
		\operatorname{grad} \Psi^{k+1}=\left( \begin{array}{c}
			\operatorname{grad}_{x_{1}} \mathcal{L}_{\beta}(\bm{x}^{k+1},\lambda^{k+1})
			+2\operatorname{Proj}_{T_{x_{1}^{k+1}}\mathcal{M}}(M_{1}(x_{1}^{k+1}-x_{1}^{k})) \\
			\nabla_{\bm{x}_{-1}}\mathcal{L}_{\beta}(\bm{x}^{k+1},\lambda^{k+1})
			+2M_{-1}(\bm{x}_{-1}^{k+1}-\bm{x}_{-1}^{k})\\
			2\operatorname{Proj}_{T_{x_{1}^{k}}\mathcal{M}}
			(-M_{1}(x_{1}^{k+1}-x_{1}^{k}))\\
			-2M_{-1}(\bm{x}_{-1}^{k+1}-\bm{x}_{-1}^{k})\\
            \nabla_{\lambda}\mathcal{L}_{\beta}(\bm{x}^{k+1},\lambda^{k+1})+\frac{2\alpha_{1}}{\tau\beta}(\lambda^{k+1}-\lambda^{k})\\
			-\frac{2\alpha_{1}}{\tau\beta}(\lambda^{k+1}-\lambda^{k})
		\end{array}\right). 
		\label{eq-lem5.6-1}
	\end{align}
	By using \eqref{eq-subsubgradient-1}, \eqref{eq-subsubgradient-5}, and the definition of $\mathcal{L}_{\beta}$, we have
	\begin{align}
		&\left\|\operatorname{grad}_{x_{1}} \mathcal{L}_{\beta}(\bm{x}^{k+1},\lambda^{k+1})+2\operatorname{Proj}_{T_{x_{1}^{k+1}}\mathcal{M}}(M_{1}(x_{1}^{k+1}-x_{1}^{k}))\right\|\nonumber\\
		=&\left\|\operatorname{Proj}_{T_{x_{1}^{k+1}}\mathcal{M}}\left( \nabla_{1}f(\bm{x}^{k+1})+\nabla_{1}g(x_{1}^{k+1},x_{2}^{k+1})-A_{1}^{\top}\lambda^{k+1}+\beta(x_{1}^{k+1}-x_{2}^{k+1})+2M_{1}(x_{1}^{k+1}-x_{1}^{k}) 
		\right) \right\|\nonumber\\
		=&\left\|\operatorname{Proj}_{T_{x_{1}^{k+1}}\mathcal{M}}\left(
		\nabla_{1}f(\bm{x}^{k+1})-\nabla_{1}f(x_{1}^{k+1},x_{2}^{k},x_{3}^{k},x_{4}^{k})+\nabla_{1}g(x_{1}^{k+1},x_{2}^{k+1})-\nabla_{1}g(x_{1}^{k+1},x_{2}^{k})
		\right. \right.\nonumber\\
		&-\left. \left. A_{1}^{\top}\Delta \lambda^{k+1}-\beta(x_{2}^{k+1}-x_{2}^{k})-H_{1}(x_{1}^{k+1}-x_{1}^{k})+2M_{1}(x_{1}^{k+1}-x_{1}^{k}) 
		\right)\right\|\nonumber\\
		\leq& 12\sigma_{\max}(M_{1})\|\Delta x_{1}^{k+1}\|
		+6(L_{f}+L_{g}+\beta)\|\Delta x_{2}^{k+1}\|
		+6L_{f}\|\Delta x_{3}^{k+1}\|\nonumber\\
		&+6L_{f}\|\Delta x_{4}^{k+1}\|
		+6\sigma_{\max}(A_{1})\|\Delta \lambda^{k+1}\|
		+6\overline{\kappa}_{1}\|\Delta x_{1}^{k}\|,
		\label{eq-lem5.6-2}
	\end{align}
	where $\sigma_{\max}(\cdot)$ represents the largest eigenvalue, and the first inequality follows from the nonexpansive of projection operator and the Lipschitz continuity of $\nabla f$ and  $\nabla g$. Similarly, from \eqref{eq-subsubgradient-2}--\eqref{eq-subsubgradient-4} and \eqref{eq-subsubgradient-5}, we have that
	\begin{align}
		&\left\|\nabla_{x_{2}}\mathcal{L}_{\beta}(\bm{x}^{k+1},\lambda^{k+1})+2M_{2}(x_{2}^{k+1}-x_{2}^{k})\right\|\nonumber\\
		=&\left\|\nabla_{2}f(\bm{x}^{k+1})+\nabla_{2}g(x_{1}^{k+1},x_{2}^{k+1})-A_{2}^{\top}\lambda^{k+1}-\beta(x_{1}^{k+1}-x_{2}^{k+1})+2M_{2}(x_{2}^{k+1}-x_{2}^{k}) \right\|\nonumber\\
		=&\left\|\nabla_{2}f(\bm{x}^{k+1})-\nabla_{2}f(x_{1}^{k+1},x_{2}^{k+1},x_{3}^{k},x_{4}^{k})
		-A_{2}^{\top}\Delta \lambda^{k+1}
		-H_{2}(x_{2}^{k+1}-x_{2}^{k})+2M_{2}(x_{2}^{k+1}-x_{2}^{k}) \right\|\nonumber\\
		\leq&8\sigma_{\max}(M_{2})\|\Delta x_{2}^{k+1}\|
		+4L_{f}\|\Delta x_{3}^{k+1}\|
		+4L_{f}\|\Delta x_{4}^{k+1}\|\nonumber\\
		&+4\sigma_{\max}(A_{2})\|\Delta \lambda^{k+1}\|
		+4\overline{\kappa}_{2}\|\Delta x_{2}^{k}\|,
		\label{eq-lem5.6-3}
	\end{align}
	\begin{align}
		&\left\|\nabla_{x_{3}}\mathcal{L}_{\beta}(\bm{x}^{k+1},\lambda^{k+1})+2M_{3}(x_{3}^{k+1}-x_{3}^{k})\right\|\nonumber\\
		=&\left\|\nabla_{3}f(\bm{x}^{k+1})-A_{3}^{\top}\lambda^{k+1}+\beta(x_{3}^{k+1}-x_{4}^{k+1})+2M_{3}(x_{3}^{k+1}-x_{3}^{k}) \right\|\nonumber\\
		=&\left\|\nabla_{3}f(\bm{x}^{k+1})-\nabla_{3}f(x_{1}^{k+1},x_{2}^{k+1},x_{3}^{k+1},x_{4}^{k})
		-A_{3}^{\top}\Delta \lambda^{k+1}\right.\nonumber\\
		&\left.-\beta(x_{4}^{k+1}-x_{4}^{k})
		-H_{3}(x_{3}^{k+1}-x_{3}^{k})+2M_{3}(x_{3}^{k+1}-x_{3}^{k}) \right\|\nonumber\\
		\leq&10\sigma_{\max}(M_{3})\|\Delta x_{3}^{k+1}\|
		+5(L_{f}+\beta)\|\Delta x_{4}^{k+1}\|
		+5\sigma_{\max}(A_{3})\|\Delta \lambda^{k+1}\|
		+5\overline{\kappa}_{3}\|\Delta x_{3}^{k}\|,
        \label{eq-lem5.6-4}
	\end{align}
	and 
	\begin{align}
		&\left\|\nabla_{x_{4}}\mathcal{L}_{\beta}(\bm{x}^{k+1},\lambda^{k+1})+2M_{4}(x_{4}^{k+1}-x_{4}^{k})\right\|\nonumber\\
		=&\left\|\nabla_{4}f(\bm{x}^{k+1})-A_{4}^{\top}\lambda^{k+1}-\beta(x_{3}^{k+1}-x_{4}^{k+1})+2M_{4}(x_{4}^{k+1}-x_{4}^{k}) \right\|\nonumber\\
		=&\left\|-A_{4}^{\top}\Delta \lambda^{k+1}
		-H_{4}(x_{4}^{k+1}-x_{4}^{k})+2M_{4}(x_{4}^{k+1}-x_{4}^{k}) \right\|\nonumber\\
		\leq&6\sigma_{\max}(M_{4})\|\Delta x_{4}^{k+1}\|
		+3\sigma_{\max}(A_{4})\|\Delta \lambda^{k+1}\|
		+3\overline{\kappa}_{4}\|\Delta x_{4}^{k}\|.
		\label{eq-lem5.6-5}
	\end{align}
	We can also bound the gradient of $\Psi^{k}$ about the dual variable $\lambda$, i.e., 
	\begin{align}
		\left\|\nabla_{\lambda}\mathcal{L}_{\beta}(\bm{x}^{k+1},\lambda^{k+1})+\frac{2\alpha_{1}}{\tau\beta}(\lambda^{k+1}-\lambda^{k})\right\|
		=\left\|\sum_{i=1}^{4}A_{i}x_{i}^{k+1}+\frac{2\alpha_{1}}{\tau\beta}(\lambda^{k+1}-\lambda^{k})\right\|
		\leq\frac{1+2\alpha_{1}}{\tau\beta}\|\Delta \lambda^{k+1}\|.
		\label{eq-lem5.6-6}
	\end{align} 
	Using an analysis technique similar to the above, the remaining three parts of \eqref{eq-lem5.6-1} can also be bounded by the residual error of the generated sequence of iterates. We omit the analysis. Finally, combining \eqref{eq-lem5.6-1}--\eqref{eq-lem5.6-6}, there exist $\varrho_1,\varrho_2, \varrho_3>0$ such that 
	\begin{align*}
		\|\operatorname{grad} \Psi^{k+1}\|\leq \varrho_1\|\Delta \bm{x}^{k+1}\|+\varrho_2\|\Delta \bm{x}^{k}\|+\varrho_3\|\Delta \lambda^{k+1}\|.
	\end{align*}
	The proof is completed.
\end{proof}

\subsection{Proof of Theorem \ref{thm-global-convergence}}
\begin{proof}
	(i) Since $\{(\bm{x}^{k},\lambda^{k})\}$ is bounded, $\Omega^{*}$ is compact. By Lemma \ref{lem-5.4} and Theorem \ref{thm-decrease} that $\{\Psi^{k}\}$ is nonincreasing and bounded from below, we know that the limit of $\Psi^{k}$ exists. Note that $\Psi$ is continuous and $\{\Delta \lambda^{k}\}$ and $\{\Delta \bm{x}^{k}\}$ converge to $0$, we have 
%	that $\Psi$ has the same value at all the points in $\omega(\bm{x}^{0},\bm{x}^{1},\lambda^{1})$, i.e.,
	\begin{align}
		\lim_{k\rightarrow +\infty} \Psi^{k}
		=\Psi^{*}:=\Psi(\hat{z})
		\quad\forall \hat{z}\in \Omega^{*}.
		\label{eq-thm-gc-2}
	\end{align}
	
	On the one hand, if there exists an integer $\bar{k}$ such that $\Psi^{\bar{k}}=\Psi^{*}$, then the nonincreasing property \eqref{eq-thm5.2-L1} would imply that $\bm{x}^{k+1}=\bm{x}^{k}$ for any $k\geq\bar{k}$. Associated with \eqref{eq-thm5.2-3}, we have that $\sum_{i=1}^{4}A_{i}x_{i}^{k}=0,~\forall k\geq\bar{k}$.
	Hence, for any $k\geq\bar{k}$, it follows from \eqref{sec5-ADMM} that $\lambda^{k+1}=\lambda^{k}$ and the assertion \eqref{eq-thm-gc-1} holds.
	
	On the other hand, assume $\Psi^{k}>\Psi^{*}$ for all $k$. From \eqref{eq-thm-gc-2}, for any $\eta_{1},\eta_{2}>0$, there exists a nonnegative integer $k_{0}$ such that for any $k>k_{0}$,
	$$
	\Psi^{k}<\Psi^{*}+\eta_{1}, \quad \operatorname{dist}\left( z^{k},\Omega^{*}\right) <\eta_{2}.
	$$
	Applying Lemma \ref{lem-uniformizedKL} with $\Omega=\Omega^{*}$, there exists $\varphi \in \mit\Phi_\eta$ such that for any $k>k_{0}$,
	\begin{align}
		\varphi^{\prime}(\Psi^{k}-\Psi^{*}) \|\operatorname{grad}\Psi^{k}\| \geq 1.
		\label{eq-thm-gc-3}
	\end{align}
	Using of the concavity of $\varphi$, we get
	\begin{align}
		\varphi(\Psi^{k}-\Psi^{*}) -\varphi(\Psi^{k+1}-\Psi^{*})
		\geq\varphi^{\prime}(\Psi^{k}-\Psi^{*})(\Psi^{k}-\Psi^{k+1}).
		\label{eq-thm-gc-4}
	\end{align}
	Then combining \eqref{eq-thm-gc-3}, \eqref{eq-thm-gc-4} with Theorem \ref{thm-decrease}, we obtain
	\begin{align}
		\varphi(\Psi^{k}-\Psi^{*}) -\varphi(\Psi^{k+1}-\Psi^{*})
		\geq C_{1}\|\operatorname{grad}\Psi^{k}\|^{-1}
		\|\Delta \bm{x}^{k}\|^{2},
		\label{eq-thm-gc-5}
	\end{align}
	where $C_{1}=\min_{i}\{(1-\delta_{i})\sigma_{\min}(M_{i})\}$.
	For convenience, we define for all $p,q\in \mathbb{N}$ the following quantities
	$$
	\Delta_{p,q}=\varphi(\Psi^{p}-\Psi^{*}) -\varphi(\Psi^{q}-\Psi^{*}).
	$$
	Combining Lemma \ref{lem-subgradient-bound} with \eqref{eq-thm-gc-5} yields for any $k>k_{0}$ that 
	\begin{align}
	 	\|\Delta \bm{x}^{k}\|^{2}
	 	\leq &\frac{1}{C_{1}}\Delta_{k,k+1} \|\operatorname{grad}\Psi^{k}\|\nonumber\\
	 	\leq&\frac{1}{C_{1}}\Delta_{k,k+1}
	 	\left( \varrho_1\|\Delta \bm{x}^{k}\|+\varrho_2\|\Delta \bm{x}^{k-1}\|+\varrho_3\|\Delta \lambda^{k}\|\right).
	 	\label{eq-thm-gc-6}
	\end{align}
 	Next, we will give a upper bound of $\|\Delta \lambda^{k}\|$. By using \eqref{eq-lemma5.1-1}, it holds that 
 	\begin{align*}
 	\Delta \lambda^{k+1} = (1-\tau)\Delta \lambda^{k}-\tau\Delta z^{k+1}.
 	\end{align*}
 	Then we consider two cases.
 	\begin{itemize}
 		\item[$\bullet$] When $0<\tau\leq1$, we have
 		\begin{align*}
 			\tau\|\Delta \lambda^{k+1}\| \leq (1-\tau)(\|\Delta \lambda^{k}\|-\|\Delta \lambda^{k+1}\|)+\tau\|\Delta z^{k+1}\|.
 		\end{align*}
 		\item[$\bullet$] When $1<\tau<2$, we have
 		\begin{align*}
 			(2-\tau)\|\Delta \lambda^{k+1}\|\leq (\tau-1)(\|\Delta \lambda^{k}\|-\|\Delta \lambda^{k+1}\|)+\tau\|\Delta z^{k+1}\|.
 		\end{align*}
 	\end{itemize}
 	In addition, using the same technique as the proof of Lemma \ref{lem-5.1}, there exist $\varrho_4,\varrho_5>0$ such that
 	\begin{align*}
 		\|\Delta z^{k+1}\|\leq\varrho_4\|\Delta\bm{x}^{k+1}\| + \varrho_5\|\Delta\bm{x}^{k}\|.
 	\end{align*}
	Hence, we obtain
	\begin{align}
		\|\Delta \lambda^{k+1}\|\leq \alpha_{1}(\|\Delta \lambda^{k}\|-\|\Delta \lambda^{k+1}\|)
		+\alpha_{3}
		(\varrho_4\|\Delta\bm{x}^{k+1}\| + \varrho_5\|\Delta\bm{x}^{k}\|),
		\label{eq-thm-gc-7}
	\end{align}
	where $\alpha_{3}=\frac{\tau}{1-|1-\tau|}$.
% 	By Lemma \ref{lem-5.1}, there exist $\varrho_4,\varrho_5,\varrho_6>0$ such that
% 	\begin{align}
% 	 	\|\Delta \lambda^{k+1}\|^{2} 
% 	 	\leq& \alpha_1\left( \|\Delta \lambda^{k}\|^{2} -\|\Delta \lambda^{k+1}\|^{2}\right) + \varrho_4\|\Delta\bm{x}^{k+1}\|^{2} + \varrho_5\|\Delta\bm{x}^{k}\|^{2} +\varrho_6\|\Delta\bm{x}^{k-1}\|^{2}\nonumber\\
% 	 	\leq&
% 	\end{align}
  	Denote $C_{2}=\max\{\varrho_1+\alpha_{3}\varrho_3\varrho_4,\varrho_2+\alpha_{3}\varrho_3\varrho_5\}$.
  	Combining \eqref{eq-thm-gc-6} with \eqref{eq-thm-gc-7} and using the fact that $2\sqrt{\alpha\beta}\leq\alpha+\beta$, we have that 
  		\begin{align}
  			\|\Delta \bm{x}^{k}\|
  			\leq&\sqrt{\frac{C_{2}}{C_{1}}\Delta_{k,k+1}\left( 
  				\|\Delta \bm{x}^{k}\|
  				+\|\Delta \bm{x}^{k-1}\|
  				+\frac{\alpha_{1}\varrho_3}{C_{2}}(\|\Delta \lambda^{k-1}\| 
  				-\|\Delta \lambda^{k}\|)
  				\right) }\nonumber\\
  			\leq& \frac{C_{2}}{C_{1}}\Delta_{k,k+1}
  			+\frac{1}{4}(\|\Delta \bm{x}^{k}\|
  			+\|\Delta \bm{x}^{k-1}\|)
  			+\frac{\alpha_{1}\varrho_3}{4C_{2}}(\|\Delta \lambda^{k-1}\|
  			-\|\Delta \lambda^{k}\|). 
  			\label{eq-thm-gc-8}
  		\end{align}	
  	Summing the inequality \eqref{eq-thm-gc-8} over $k = k_{0}+2,\dots,m$ yields
  	\begin{align*}
  		\frac{1}{2}\sum_{k=k_{0}+2}^{m}\|\Delta \bm{x}^{k}\|
  		\leq& \frac{C_{2}}{C_{1}}\Delta_{k_{0}+3,m}
  		+ \frac{1}{4}\|\Delta \bm{x}^{k_{0}+1}\|
  		+\frac{\alpha_{1}\varrho_3}{4C_{2}}\|\Delta \lambda^{k_{0}+2}\|
  		\nonumber\\
  		\leq& \frac{C_{2}}{C_{1}}\varphi(\Psi^{k_{0}+3}-\Psi^{*})
  		+ \frac{1}{4}\|\Delta \bm{x}^{k_{0}+1}\|
  		+\frac{\alpha_{1}\varrho_3}{4C_{2}}\|\Delta \lambda^{k_{0}+2}\|.
  	\end{align*}
  	Letting $m\rightarrow +\infty$,  we get
  	\begin{align*}
  		\sum_{k=1}^{+\infty}\|\bm{x}^{k+1}-\bm{x}^{k}\|< +\infty.
  	\end{align*}
  	Further, \eqref{eq-thm-gc-7} implies that 
  	\begin{align*}
  		\sum_{k=1}^{+\infty}\|\lambda^{k+1}-\lambda^{k}\|< +\infty.
  	\end{align*}
    Hence, the sequence $\{(\bm{x}^{k},\lambda^{k})\}$  is a Cauchy sequence, which
    converges. The
    assertion then follows immediately from Theorem \ref{thm-decrease}.
\end{proof}

\section{Supplemental Numerical Experiments}  \label{app-numerical experiments}
\subsection{Quantitative Indicators and Stopping Criterion}
We measure the quality of restoration by the relative error (Rel.Err.) and Normalized Root Mean Square Error (NRMSE), which are respectively defined as
\begin{align*}
	\operatorname{Rel.Err.}  & = \frac{\|\tilde{\bm{q}}-\tilde{\bm{q}}_{0}\|+\|\bm{t}-\bm{t}_{0}\|}{\|\tilde{\bm{q}}_{0}\|+\|\bm{t}_{0}\|},\\
	\operatorname{NRMSE} & = \frac{\|\tilde{\bm{q}}-\tilde{\bm{q}}_{0}\|+\|\bm{t}-\bm{t}_{0}\|}{(\max(\bm{t})-\min(\bm{t}))\sqrt{n}},
\end{align*}
where $(\tilde{\bm{q}},\bm{t})$ is the restored pose and $(\tilde{\bm{q}}_{0},\bm{t}_{0})$ is the true pose. In practical experiments, we use $(\tilde{\bm{p}}^{k},\bm{t}^{k})$ for computation since $\|\tilde{q}_{i}^{k}\| \neq 1$. We adopt the residual
\begin{align*}
	e^{k+1}&=\frac{1}{\beta}\|\bm{\lambda}^{k+1}-\bm{\lambda}^{k}\|^{2}
	+\frac{1}{\beta}\|\bm{z}^{k+1}-\bm{z}^{k}\|^{2}
	+\beta\left( \|\tilde{\bm{q}}^{k+1}-\tilde{\bm{q}}^{k}\|^{2}
	+\|\bm{t}^{k+1}-\bm{t}^{k}\|^{2}\right),
\end{align*}
to quantify the accuracy of PRADMM. As $e^{k} \rightarrow 0$, it implies that  $\|\tilde{\bm{p}}^{k+1}-\tilde{\bm{p}}^{k}\|$ also converges to zero. We terminate the solvers when the iteration residual $e^{k}<tol$ or the maximum number of iterations $MaxIter$ is reached.

\subsection{More Results of Ring Datasets}
We test the circular ring datasets with different sizes and list the numerical results about relative error, NRMSE, and CPU time in Table \ref{app-tab:circular ring_results}. Our experiments demonstrate that PRADMM achieves both faster computation and higher accuracy on small-scale circular ring datasets. As dataset size increases, PRADMM maintains comparable or occasionally superior precision to RS+PS despite marginally incremental runtime requirements.
Compared with SE-sync, PRADMM achieves significantly faster computation times across all scales while maintaining superior solution accuracy in most scenarios. 

\begin{table}[tbp]
	\centering
	\caption{Numerical results of different sizes and noise levels of circular ring datasets.}
	\label{app-tab:circular ring_results}
        \resizebox{0.8\linewidth}{!}{
	\begin{tabular}{l|ccc|ccc|ccc}
		\toprule[1pt]
		\multicolumn{1}{c}{}& 
		\multicolumn{3}{c}{$\sigma_{r}=0.01$, $\sigma_{t}=0.01$} & 
		\multicolumn{3}{c}{$\sigma_{r}=0.03$, $\sigma_{t}=0.05$} & 
		\multicolumn{3}{c}{$\sigma_{r}=0.05$, $\sigma_{t}=0.1$} \\
		\cmidrule(lr){2-4}\cmidrule(lr){5-7}\cmidrule(lr){8-10}
		\multicolumn{1}{c}{Algorithm} & Rel. Err. & NRMSE & \multicolumn{1}{c}{Time (s)} & Rel. Err. & NRMSE & \multicolumn{1}{c}{Time (s)} & Rel. Err. & NRMSE & Time (s) \\
		\midrule[1pt]
		\multicolumn{10}{c}{$m=n=100$} \\
		\midrule[0.5pt]
		mG-N & 0.0711 & 0.0354 & 0.407 & 0.3683 & 0.1834 & 0.401 & 0.5024 & 0.2502 & 0.407 \\
		SE-sync & 0.0711 & 0.0354 & 0.179 & 0.3683 & 0.1834 & 0.166 & 0.5025 & 0.2502 & 0.161 \\
		RS+PS & 0.0699 & 0.0348 & 0.069 & 0.3457 & 0.1721 & 0.065 & 0.4861 & 0.2420 & 0.072 \\
		RGD & 0.0692 & 0.0344 & 0.322 & 0.3087 & 0.1537 & 0.355 & 0.4729 & 0.2355 & 0.354 \\
		SOC & 0.0691 & 0.0344 & 0.647 & 0.3085 & 0.1536 & 0.425 & 0.4729 & 0.2354 & 0.436 \\
		PieADMM & 0.0689 & 0.0343 & 0.123 & 0.3085 & 0.1536 & 0.131 & 0.4729 & 0.2354 & 0.215 \\
		PRADMM & \textbf{0.0689} & \textbf{0.0343} & \textbf{0.065} & \textbf{0.3085} & \textbf{0.1536} & \textbf{0.034} & \textbf{0.4724} & \textbf{0.2352} & \textbf{0.042} \\
		\midrule[0.5pt]
		\multicolumn{10}{c}{$m=n=500$} \\
		\midrule[0.5pt]
		mG-N & 0.0892 & 0.0445 & 1.257 & 0.2801 & 0.1399 & 1.621 & 0.4766 & 0.2381 & 1.610 \\
		SE-sync & \textbf{0.0892} & \textbf{0.0445} & 0.257 & 0.2801 & 0.1399 & 0.289 & \textbf{0.4766} & \textbf{0.2381} & 0.204 \\
		RS+PS & 0.0911 & 0.0455 & \textbf{0.081} & 0.2856 & 0.1427 & \textbf{0.083} & 0.4851 & 0.2423 & \textbf{0.084} \\
		RGD & 0.1002 & 0.0500 & 0.886 & 0.2987 & 0.1487 & 0.450 & 0.4881 & 0.2438 & 0.780 \\
		SOC & 0.0916 & 0.0458 & 5.640 & 0.2997 & 0.1497 & 5.512 & 0.4851 & 0.2423 & 4.285 \\
		PieADMM & 0.0903 & 0.0451 & 0.445 & 0.2993 & 0.1495 & 0.457 & 0.4853 & 0.2424 & 0.463 \\
		PRADMM & 0.0908 & 0.0454 & 0.108 & \textbf{0.2856} & \textbf{0.1427} & 0.151 & 0.4851 & 0.2423 & 0.141 \\
		\midrule[0.5pt]
		\multicolumn{10}{c}{$m=n=1000$} \\
		\midrule[0.5pt]
		mG-N & 0.0463 & 0.0232 & 3.152 & 0.1557 & 0.0778 & 3.084 & 0.2729 & 0.1364 & 3.085 
		\\
		SE-sync & 0.0463 & 0.0232 & 0.381 & 0.1557 & 0.0778 & 0.405 & \textbf{0.2729} & \textbf{0.1364} & 0.428 \\
		RS+PS & 0.0457 & 0.0229 & \textbf{0.107} & 0.1560 & 0.0780 & \textbf{0.101} & 0.2747 & 0.1373 & \textbf{0.101} \\
		RGD & 0.0460 & 0.0230 & 4.476 & \textbf{0.1546} & \textbf{0.0772} & 4.245 & 0.2823 & 0.1411 & 4.243 \\
		SOC & 0.0459 & 0.0229 & 7.641 & 0.1559 & 0.0779 & 6.257 & 0.2746 & 0.1373 & 5.167 \\
		PieADMM & 0.0459 & 0.0230 & 1.080 & 0.1565 & 0.0783 & 1.267 & 0.2743 & 0.1371 & 1.054 \\
		PRADMM & \textbf{0.0457} & \textbf{0.0229} & 0.224 & 0.1557 & 0.0778 & 0.341 & 0.2741 & 0.1370 & 0.337 \\
		\midrule[0.5pt]
		\multicolumn{10}{c}{$m=n=5000$} \\
		\midrule[0.5pt]
		mG-N & 0.0451 & 0.0225 & 4.649 & 0.1448 & 0.0724 & 6.115 & 0.2477 & 0.1238 & 5.965 \\
		SE-sync & 0.0451 & 0.0225 & 0.716 & 0.1448 & 0.0724 & 0.745 & 0.2477 & 0.1238 & 0.756 \\
		RS+PS & 0.0439 & 0.0219 & \textbf{0.174} & 0.1403 & 0.0702 & \textbf{0.174} & 0.2414 & 0.1207 & \textbf{0.186} \\
		RGD & 0.0452 & 0.0226 & 9.715 & 0.1435 & 0.0717 & 9.654 & 0.2463 & 0.1231 & 9.612 \\
		SOC & 0.0442 & 0.0221 & 17.301 & 0.1402 & 0.0701 & 17.538 & 0.2414 & 0.1207 & 17.699 \\
		PieADMM & 0.0449 & 0.0224 & 4.291 & 0.1434 & 0.0717 & 5.150 & 0.2462 & 0.1231 & 5.076 \\
		PRADMM & \textbf{0.0439} & \textbf{0.0219} & 0.264 & \textbf{0.1403} & \textbf{0.0702} & 0.236 & \textbf{0.2414} & \textbf{0.1207} & 0.237 \\
		\bottomrule[1pt]
	\end{tabular}}
\end{table}

\subsection{More Results of Cube Datasets}
\begin{figure*}[!tbp]
	\begin{minipage}{.26\linewidth}
		\centering
		\subfloat[$n$ and $m$]{\label{cube_numedge2}	
			\includegraphics[width=1\linewidth]{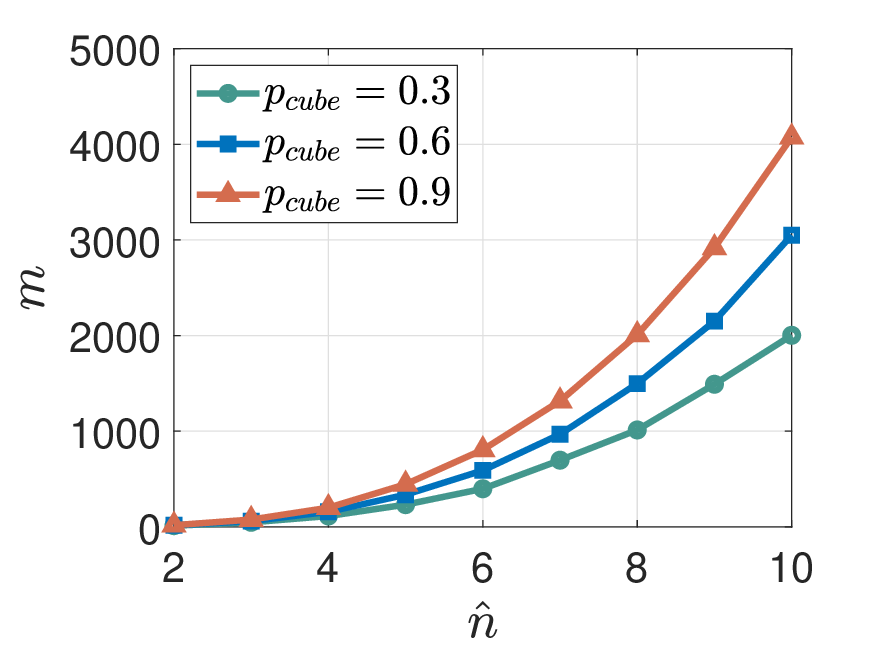}}
	\end{minipage}
	\begin{minipage}{.74\linewidth}
		\centering
		\subfloat[$p_{cube}=0.3$]{\label{cube_03_re2}	\includegraphics[width=0.33\linewidth]{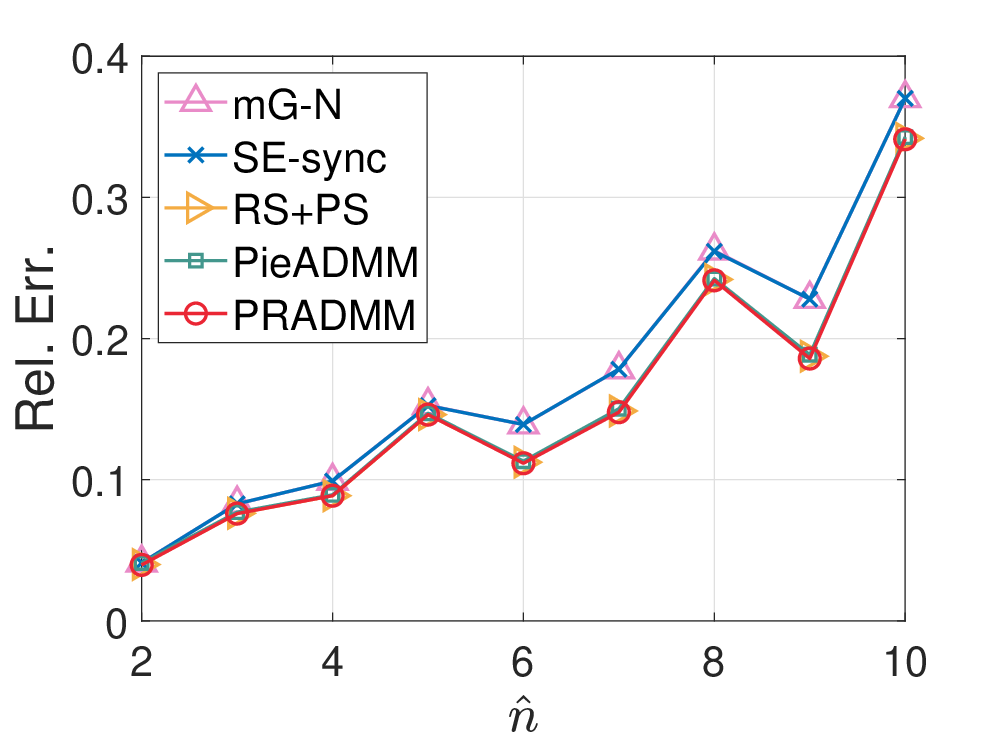}}
		\subfloat[$p_{cube}=0.6$]{\label{cube_06_re2}	\includegraphics[width=0.33\linewidth]{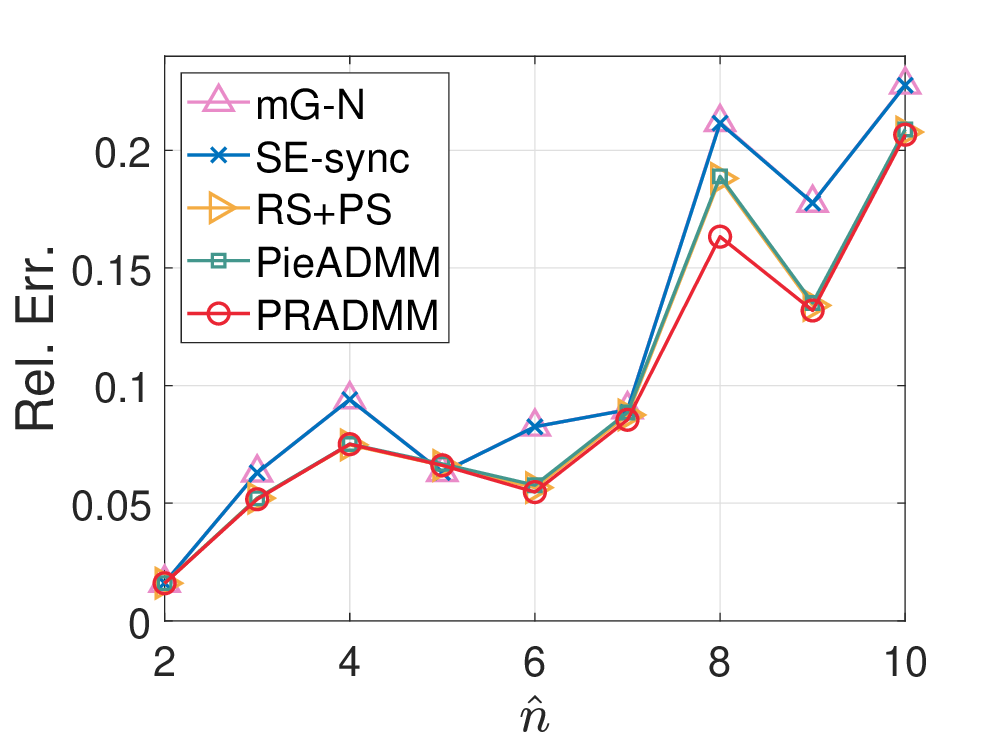}}	
		\subfloat[$p_{cube}=0.9$]{\label{cube_09_re2}	\includegraphics[width=0.33\linewidth]{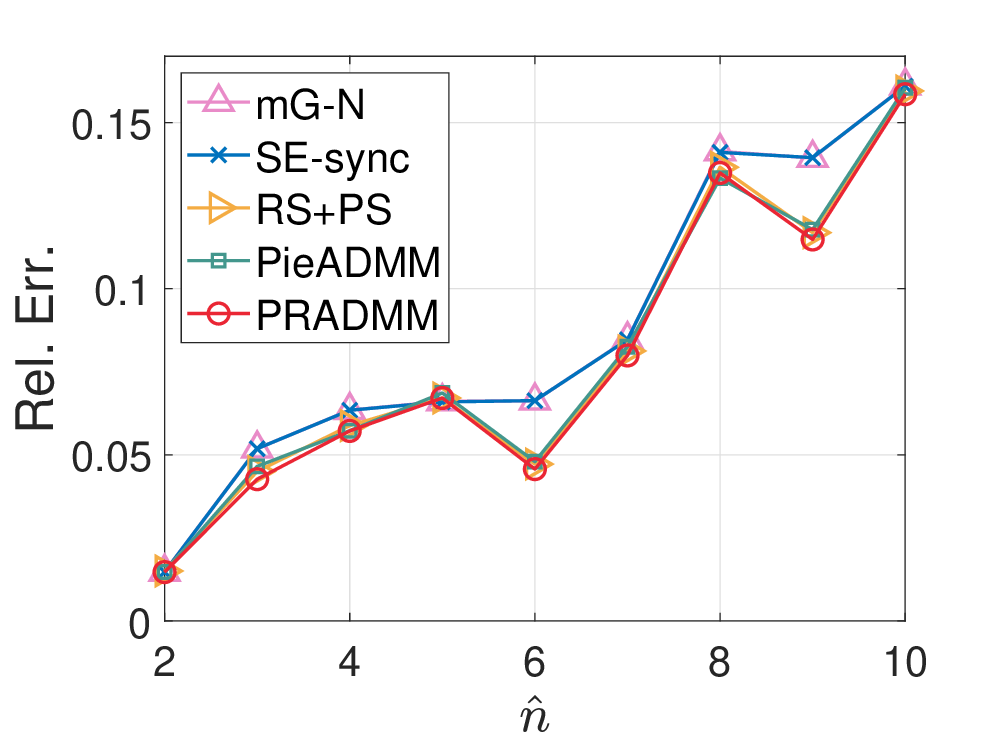}}
		
		\subfloat[$p_{cube}=0.3$]{\label{cube_03_time2}	\includegraphics[width=0.33\linewidth]{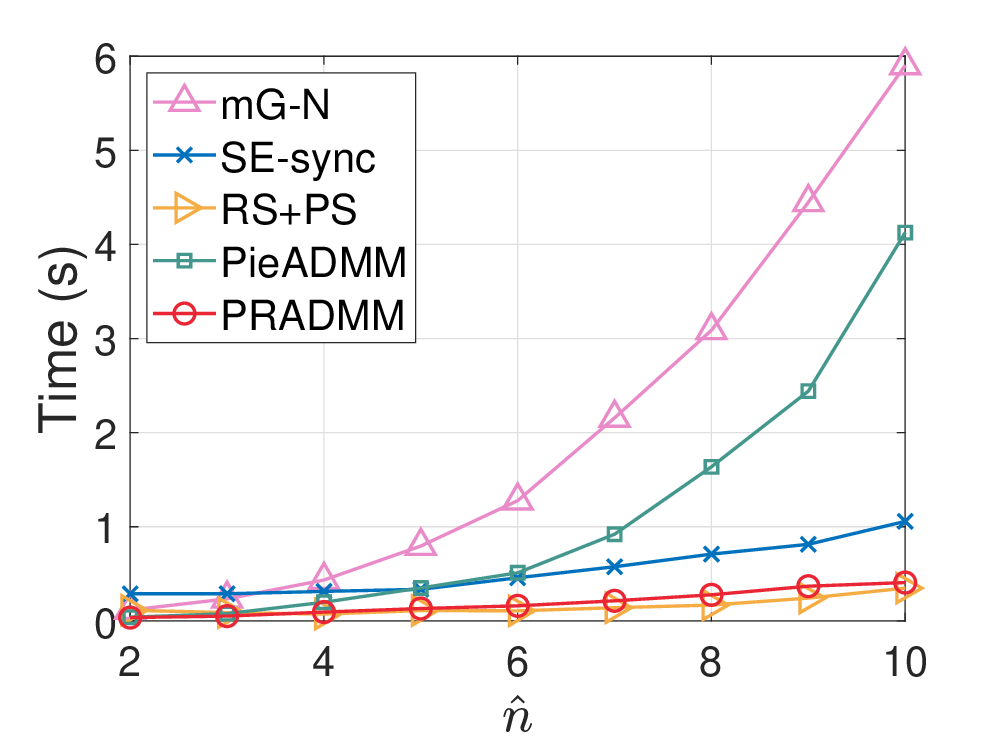}}
		\subfloat[$p_{cube}=0.6$]{\label{cube_06_time2}	\includegraphics[width=0.33\linewidth]{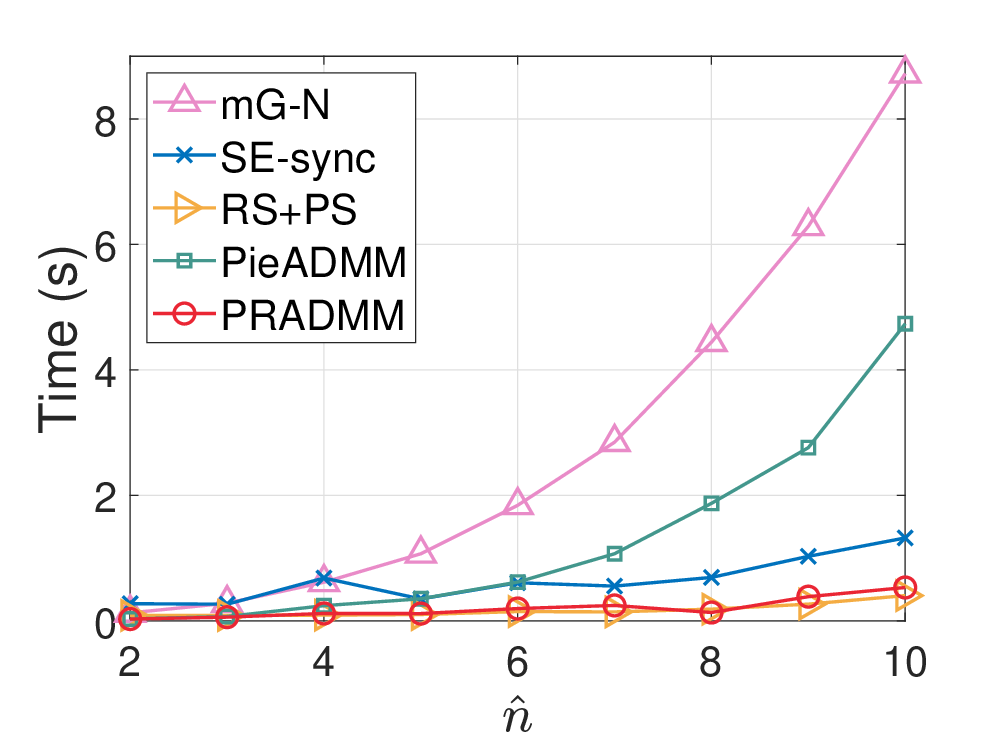}}
		\subfloat[$p_{cube}=0.9$]{\label{cube_09_time2}	\includegraphics[width=0.33\linewidth]{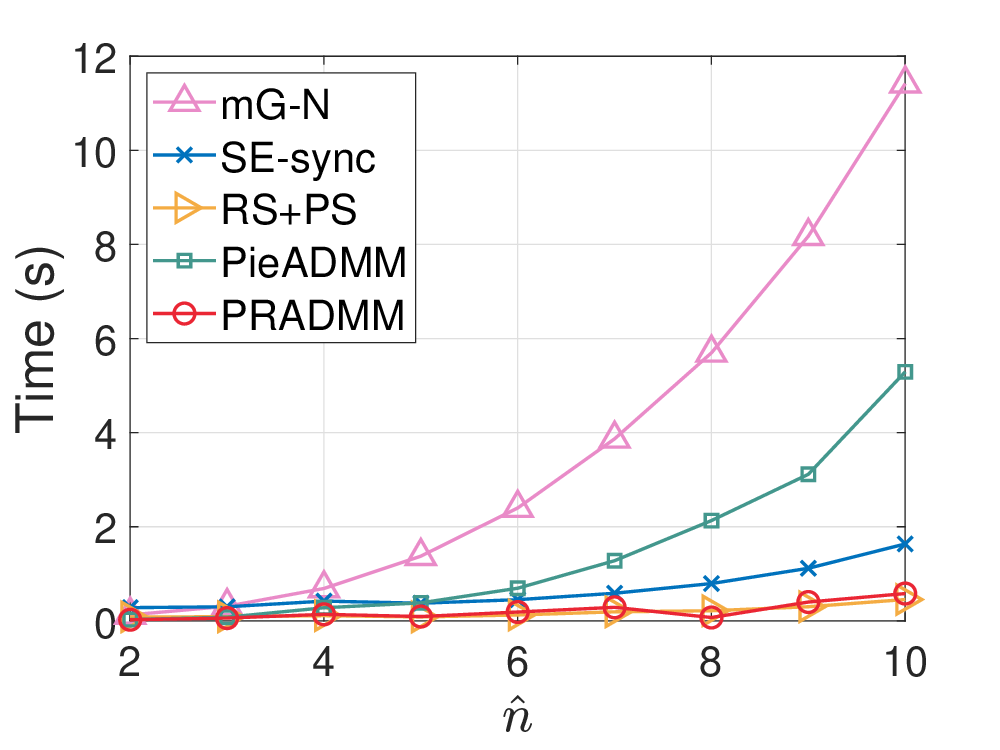}}
	\end{minipage}
	\caption{The trend of the number of edges, relative error, and CPU time along with different $\hat{n}$ under $\sigma_{t}^{rel}=0.1$, $\sigma_{r}=0.1$ and $p_{cube}\in\{0.3,0.6,0.9\}$.}
	\label{fig:cube_time2}
\end{figure*}

We define the translational noise parameter as $\sigma_{t}=\sigma_{t}^{rel}/\hat{n}$, where $\sigma_{t}^{rel}$ represents the relative noise level of translation and is used in our tests. We evaluate $\hat{n}$ values ranging from $2$ to $10$ under fixed noise conditions $\sigma_{t}^{rel}=0.1$ and $\sigma_{r}=0.1$. We list nine scenarios combining parameters $\hat{n}\in\{2,6,10\}$ and $p_{cube}\in\{0.3,0.6,0.9\}$, with corresponding observation counts $m$ listed in Table \ref{app-tab:cube_results}. The trends in the number of edges, relative error, and CPU time for different $\hat{n}$ values are shown in Figure \ref{fig:cube_time2}.
PRADMM achieves superior reconstruction accuracy compared to alternative algorithms while exhibiting faster convergence on small-scale datasets. For large-scale PGO problems, PRADMM maintains competitive efficiency -- operating marginally slower than RS+PS yet substantially outperforming all other methods in computation speed.
\begin{figure}[!tbp]
	\centering
	\subfloat{\label{cube_03_4_trajectory}	\includegraphics[width=1\linewidth]{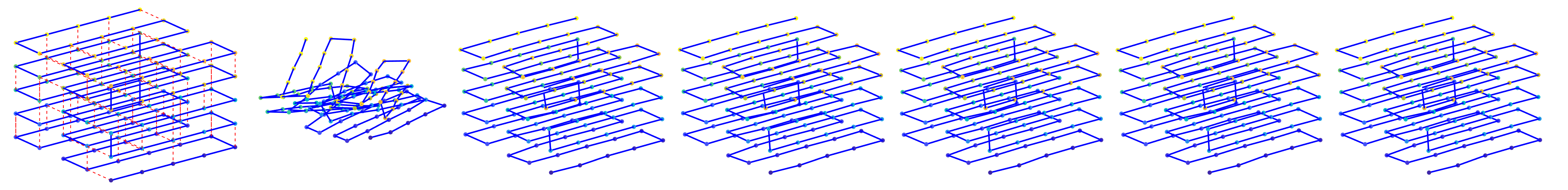}}

    {\raggedright \quad~~ Ground truth\qquad\quad Noisy \qquad\qquad RE:0.1527 \qquad~ RE:0.1525 \qquad~ RE:0.1462 \qquad~ RE:0.1470 \qquad~ RE:0.1461 \qquad\\}
    \caption{The comparison of cube trajectory with $\hat{n}=5$. From left to right are the real trajectory, the corrupted trajectory and the recovered results by mG-N, SE-sync, RS+PS, PieADMM and PRADMM, respectively.}
	\label{fig:cube_trajectory2}
\end{figure}

\begin{table}[!tbp]
	\centering
	\footnotesize
	\renewcommand{\arraystretch}{1.0}
	\caption{Numerical results of different $\hat{n}$ and $p_{cube}$ of cube datasets with $\sigma_{r}=0.1$, $\sigma^{rel}_{t}=0.1$.}
	\label{app-tab:cube_results}
	\resizebox{0.8\linewidth}{!}{
	\begin{tabular}{l|ccc|ccc|ccc}
		\toprule[1pt]
		\multicolumn{1}{c}{} & \multicolumn{3}{c}{$\hat{n} = 2$} & \multicolumn{3}{c}{$\hat{n} = 6$} & \multicolumn{3}{c}{$\hat{n} = 10$} \\
		\cmidrule(lr){2-4}\cmidrule(lr){5-7}\cmidrule(lr){8-10}
		\multicolumn{1}{c}{Algorithm} & Rel. Err. & NRMSE & \multicolumn{1}{c}{Time (s)} & Rel. Err. & NRMSE & \multicolumn{1}{c}{Time (s)} & Rel. Err. & NRMSE & \multicolumn{1}{c}{Time (s)} \\
		\midrule[0.5pt]
		\multicolumn{10}{c}{$p_{cube}=0.3$, $m=12,~	398,~	2003$.} \\
		\midrule[0.5pt]
		mG-N & 0.0411 & 0.0544 & 0.114 & 0.1391 & 0.1616 & 1.277 & 0.3702 & 0.4230 & 5.903 \\
		SE-sync & 0.0411 & 0.0543 & 0.214 & 0.1391 & 0.1616 & 0.272 & 0.3702 & 0.4229 & 0.434 \\
		RS+PS & 0.0399 & 0.0528 & 0.113 & 0.1124 & 0.1306 & \textbf{0.106} & 0.3419 & 0.3906 & \textbf{0.347} \\
		PieADMM & 0.0410 & 0.0542 & 0.038 & 0.1131 & 0.1314 & 0.510 & 0.3427 & 0.3916 & 4.126 \\
		PRADMM & \textbf{0.0396} & \textbf{0.0524} & \textbf{0.035} & \textbf{0.1116} &\textbf{0.1297} & 0.159 & \textbf{0.3412} & \textbf{0.3899} & 0.407 \\
		\midrule[0.5pt]
		\multicolumn{10}{c}{$p_{cube}=0.6$, $m=16,	~590, ~ 3050$.} \\
		\midrule[0.5pt]
		mG-N & 0.0161 & 0.0213 & 0.127 & 0.0825 & 0.0958 & 1.842 & 0.2279 & 0.2604 & 8.723 \\
		SE-sync & 0.0161 & 0.0213 & 0.192 & 0.0825 & 0.0959 & 0.392 & 0.2276 & 0.2601 & 0.582 \\
		RS+PS & 0.0160 & 0.0211 & 0.087 & 0.0566 & 0.0658 & \textbf{0.148} & 0.2078 & 0.2374 & \textbf{0.402} \\
		PieADMM & 0.0160 & 0.0212 & 0.030 & 0.0576 & 0.0670 & 0.618 & 0.2089 & 0.2387 & 4.735 \\
		PRADMM & \textbf{0.0160} & \textbf{0.0211} & \textbf{0.029} & \textbf{0.0547} & \textbf{0.0636} & 0.196 & \textbf{0.2066} & \textbf{0.2361} & 0.531 \\
		\midrule[0.5pt]	
		\multicolumn{10}{c}{$p_{cube}=0.9$, $m=17,~	806,~	4075$.} \\
		\midrule[0.5pt]
		mG-N & 0.0147 & 0.0195 & 0.126 & 0.0663 & 0.0770 & 2.399 & 0.1611 & 0.1841 & 11.417 \\
		SE-sync & 0.0147 & 0.0195 & 0.206 & 0.0663 & 0.0770 & 0.256 & 0.1610 & 0.1839 & 0.652 \\
		RS+PS & 0.0150 & 0.0199 & 0.085 & 0.0472 & 0.0548 & \textbf{0.124} & 0.1596 & 0.1823 & \textbf{0.453} \\
		PieADMM & 0.0147 & 0.0195 & 0.027 & 0.0479 & 0.0557 & 0.693 & 0.1605 & 0.1834 & 5.292 \\
		PRADMM & \textbf{0.0147} & \textbf{0.0194} & \textbf{0.023} & \textbf{0.0457} & \textbf{0.0531} & 0.187 & \textbf{0.1586} & \textbf{0.1812} & 0.580 \\
		\bottomrule[1pt]
	\end{tabular}}
\end{table}

\subsection{More Results of Benchmark Datasets}
\begin{figure*}[!tbp]
	\centering
	\subfloat{\label{tinyGrid3D_iRADMM}	\includegraphics[width=0.8\linewidth]{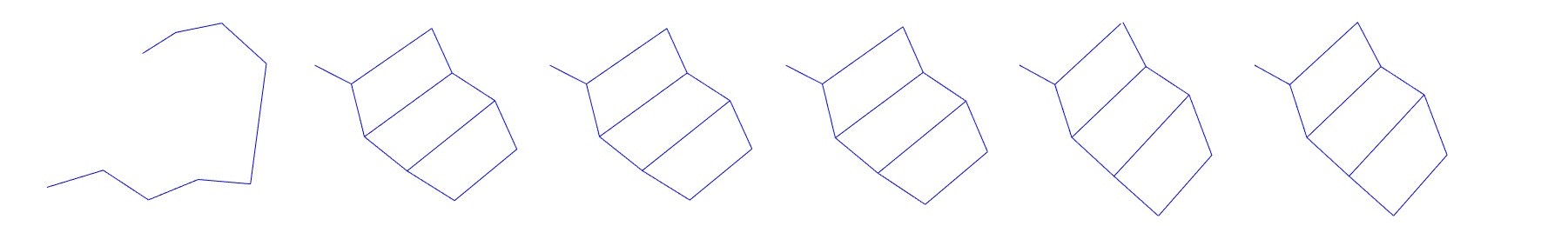}}
	
	% \vspace{-1em}
	\subfloat{\label{garage_iRADMM}	
    \includegraphics[width=0.8\linewidth]{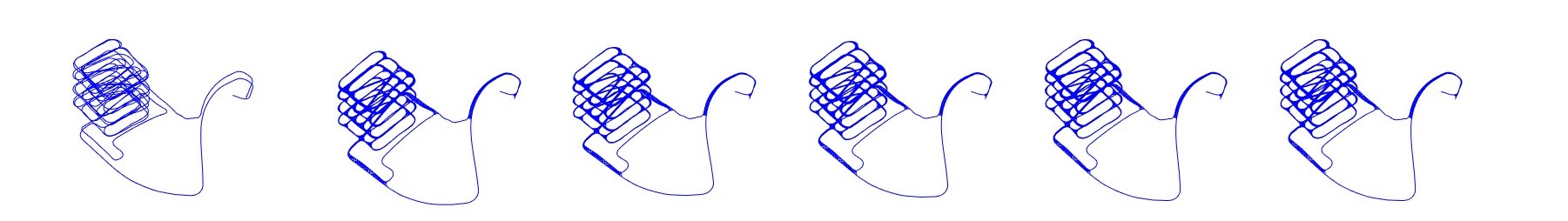}}
	
	% \vspace{-1.4em}
	\subfloat{\label{sphere1_iRADMM}	\includegraphics[width=0.8\linewidth]{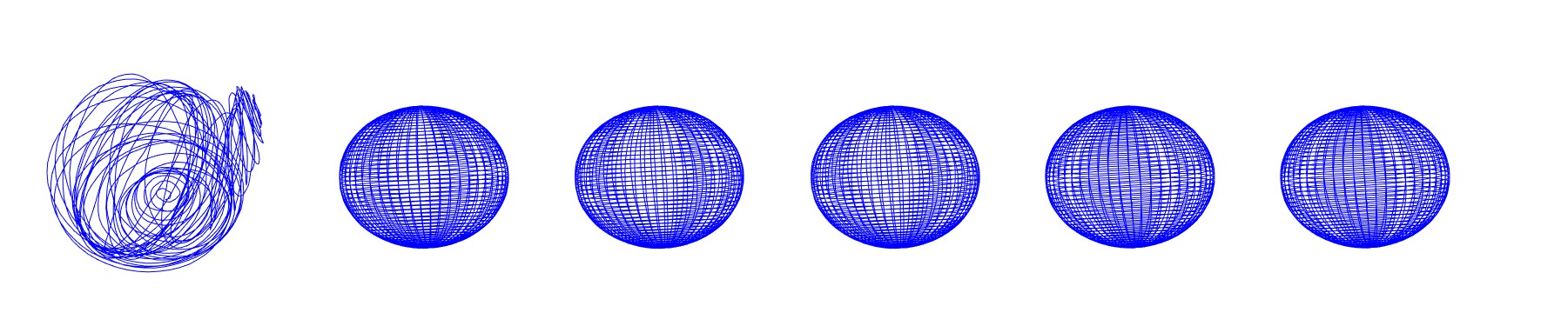}}
	
	% \vspace{-2em}
	\subfloat{\label{sphere2_iRADMM}	
	~\includegraphics[width=0.8\linewidth,height=2.5cm,keepaspectratio=false]{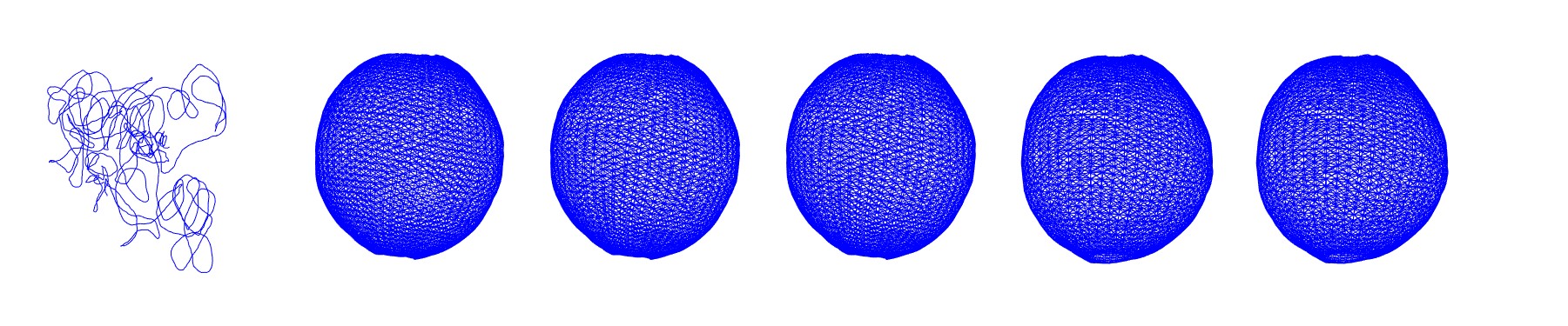}}
	
	% \vspace{-2em}
	\subfloat{\label{torus3D_iRADMM}	\includegraphics[width=0.8\linewidth]{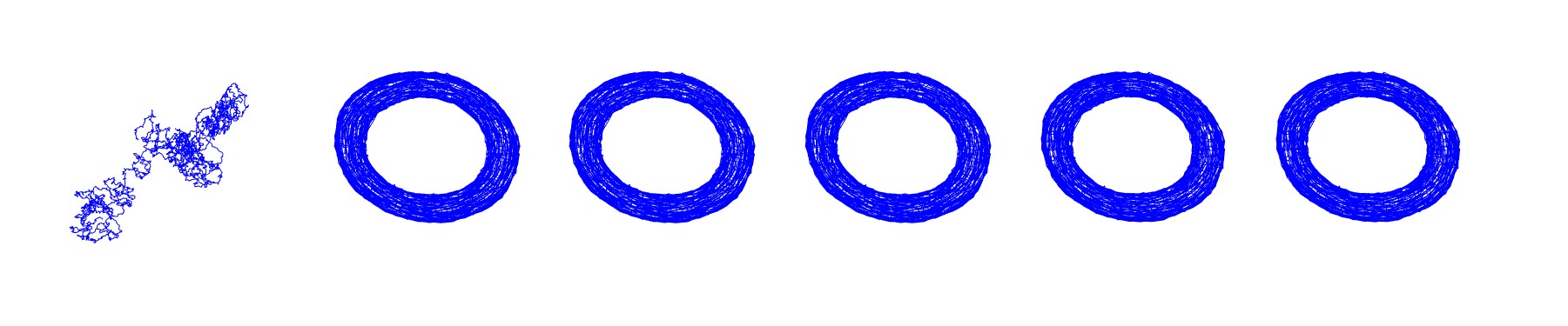}}
	
	~ Noised data  ~~\quad\quad mG-N ~\quad\qquad SE-sync ~~~\quad\quad RS+PS ~\quad\quad PieADMM ~~\quad PRADMM~~
	\caption{The results of SLAM benchmark datasets in visual. From top to bottom are the tinyGrid, garage, torus, sphere $1$, and sphere $2$ datasets, respectively. From left to right are the corrupted data and the recovered results by mG-N, SE-sync, RS+PS, PieADMM, and PRADMM, respectively.}
	\label{benchmark trajectory}
\end{figure*}

Figure \ref{benchmark trajectory} shows the results of the trajectory in visual, and the corresponding numerical results are listed in Table \ref{app-benchmark_table}.

\begin{table}[ht]
	\centering
	\renewcommand{\arraystretch}{1.0}
	\caption{The numerical results of errors in rotation angle and translation, along with the CPU time consumed of SLAM benchmark datasets.}
        \resizebox{0.7\linewidth}{!}{
	\begin{tabular}{ccccccc}
		\toprule[1pt]
		Datasets & Size & Algorithm & Loss$(\theta)$ & Loss$(\tilde{q})$ & Loss$(\bm{t})$ & Time\,(s)\\
		\midrule
		\multirow{5}{*}{tinyGrid} 
		& \multirow{5}{*}{$\begin{array}{c}
				n=9\\m=11
		\end{array}$} 
		& mGN & 7.04 & 1.76 & 4.54 & 0.122 \\
		& & SE-sync & 7.06 & 1.76 & 4.49 & 0.093 \\
		& & RS+PS & 7.01 & 1.75 & 4.61 & 0.002 \\
		& & PieADMM & 8.04 & 2.01 & 4.23 & 0.046 \\
		& & PRADMM & 8.06 & 2.01 & 4.22 & 0.024 \\
		\midrule
		\multirow{5}{*}{garage} 
		& \multirow{5}{*}{$\begin{array}{c}
				n=1661\\m=6275
			\end{array}$}
		& mGN & 1.23e-02 & 3.20e-03 & 1.24 & 16.988 \\
		& & SE-sync & 1.25e-02 & 3.21e-03 & 1.24 & 4.890 \\
		& & RS+PS & 1.28e-02 & 3.20e-03 & 1.24 & 0.387 \\
		& & PieADMM & 8.94e-03 & 2.32e-03 & 1.30 & 10.371 \\
		& & PRADMM & 1.03e-02 & 2.67e-03 & 1.30 & 0.180 \\
		\midrule
		\multirow{5}{*}{sphere 1} 
		& \multirow{5}{*}{$\begin{array}{c}
				n=2500\\m=4949
			\end{array}$} 
		& mGN & 5.43e+02 & 1.36e+02 & 6.45e+02 & 13.656 \\
		& & SE-sync & 5.03e+02 & 1.26e+02 & 6.82e+02 & 0.812 \\
		& & RS+PS & 5.03e+02 & 1.26e+02 & 6.82e+02 & 0.690 \\
		& & PieADMM & 6.96e+02 & 1.74e+02 & 5.72e+02 & 5.817 \\
		& & PRADMM & 7.95e+02 & 1.99e+02 & 5.66e+02 & 0.444 \\
		\midrule
		\multirow{5}{*}{sphere 2} 
		& \multirow{5}{*}{$\begin{array}{c}
				n=2200\\m=8647
			\end{array}$}
		& mGN & 1.51e+06 & 3.75e+05 & 8.94e+03 & 23.739 \\
		& & SE-sync & 1.49e+06 & 3.73e+05 & 9.23e+03 & 1.142 \\
		& & RS+PS & 1.49e+06 & 3.73e+05 & 9.23e+03 & 0.987 \\
		& & PieADMM & 1.63e+06 & 4.06e+05 & 7.63e+03 & 5.697 \\
		& & PRADMM & 1.63e+06 & 4.05e+05 & 8.13e+03 & 0.387 \\
		\midrule
		\multirow{5}{*}{torus} 
		& \multirow{5}{*}{$\begin{array}{c}
				n=5000\\m=9048
			\end{array}$} 
		& mGN & 6.21e+03 & 1.55e+03 & 1.18e+04 & 25.306 \\
		& & SE-sync & 6.21e+03 & 1.55e+03 & 1.18e+04 & 1.525 \\
		& & RS+PS & 6.21e+03 & 1.55e+03 & 1.18e+04 & 4.210 \\
		& & PieADMM & 7.92e+03 & 1.98e+03 & 1.12e+04 & 22.194 \\
		& & PRADMM & 7.06e+03 & 1.76e+03 & 1.19e+04 & 1.037 \\
		\bottomrule[1pt]
	\end{tabular}}
	\label{app-benchmark_table}
\end{table}
\end{document}